\documentclass{article}

\usepackage{amsfonts}
\usepackage{amssymb}
\usepackage{textcomp}
\usepackage{latexsym,amsfonts,amssymb,amsmath,epsfig}

\usepackage{hyperref}

\begin{document}

\newcommand{\spanz}{\mathop{\rm span}}

\renewcommand{\phi}{\varphi}
\newcommand{\be}{\begin{equation}}
\newcommand{\ee}{\end{equation}}
\newcommand{\ba}{\begin{eqnarray}}
\newcommand{\ea}{\end{eqnarray}}
\newcommand{\ban}{\begin{eqnarray*}}
\newcommand{\ean}{\end{eqnarray*}}

\newcommand{\too}{\mathop{\longrightarrow}}
\newcommand{\esssup}{\mathop{\rm esssup}}
\newcommand{\nul}{{\bf0}}
\newcommand{\odin}{{\bf1}}
\newcommand{\cupm}{\mathop{\cup}}
\newcommand{\capm}{\mathop{\cap}}
\renewcommand{\t}{{\mathbb T}}
\newcommand{\td}{{\mathbb T}^d}
\newcommand{\rd}{{\mathbb R}^d}
\newcommand{\zdd}{{\mathbb Z}^{d-1}}
\newcommand{\zddp}{{\mathbb Z}^{d-1}_+}
\newcommand{\rdd}{{\mathbb R}^{d-1}}
\newcommand{\zd}{{\mathbb Z}^{d}}
\newcommand{\tdd}{{\mathbb T}^{d-1}}
\newcommand{\ld}{L(\td)}
\renewcommand{\r}{{\mathbb R}}
\newcommand{\z} {{\mathbb Z}}
\newcommand{\cn} {{\mathbb C}}
\newcommand{\n} {{\mathbb N}}
\newcommand{\nd} {{``````}^d}

\newcommand{\calt}{{\cal T}}
\newcommand{\call}{{\cal L}}
\newcommand{\cons}{{C\,}}
\newcommand{\ddd}{,\dots,}
\renewcommand{\lll}{\left(}
\newcommand{\rrr}{\right)}
\newcommand{\ex}[1]{e^{2\pi i{#1}}}
\newcommand{\exm}[1]{e^{-2\pi i{#1}}}
\newcommand{\sml}[3]{\sum\limits_{{#1}={#2}}^{#3}}
\newcommand{\Ldvad}{L_2(\rd)}
\newcommand{\h}{\widehat}
\newcommand{\smzd}[1]{\sum\limits_{#1\in\zd}}
\newcommand{\inti}{\int\limits_{-\infty}^\infty}
\newcommand{\w}{\widetilde}

\newtheorem{theo}{Theorem}
\newtheorem{lem}[theo]{Lemma}
\newtheorem {prop} [theo] {Proposition}
\newtheorem {coro} [theo] {Corollary}
\newtheorem {defi} [theo] {Definition}
\newtheorem {rem} [theo] {Remark}

\newcommand{\tocsecindent}{\hspace{0mm}}

\centerline{\Large\bf{On construction of multivariate symmetric }}

\centerline{\Large\bf{MRA-based
wavelets\footnote {
This research was supported by Grant 09-01-00162 of RFBR.
}
}
}
\bigskip \centerline{\large A. Krivoshein} \bigskip

\centerline{St.Petersburg State University}
\centerline{e-mail: KrivosheinAV@gmail.com}

\vspace{1cm}
\centerline{\bf Abstract}
For an arbitrary matrix dilation, any integer $n$ and any  integer/semi-integer $c$,
    we describe all masks that are symmetric with respect to the point $c$
    and have sum rule of order $n$. For each such mask, we give explicit formulas
    for wavelet functions that are point symmetric/antisymmetric and generate frame-like wavelet system
    providing approximation order $n.$
For any matrix dilations (which are appropriate for axial symmetry group on $\r^2$
    in some natural sense) and given integer~$n$,
    axial symmetric/antisymmetric frame-like wavelet systems providing
    approximation order $n$ are constructed.
Also, for several matrix dilations the explicit construction
    of highly symmetric
    frame-like wavelet systems providing
    approximation order $n$ is presented.

\vspace{.5cm}
{\bf Keywords:} MRA-based  wavelet systems, frame-type
expansion, approximation order, symmetry

\vspace{.5cm}
{\bf Mathematics Subject Classifications:} 42C40

\section{Introduction}
%\pdfbookmark[1]{Introduction}{Introduction}

In recent years, methods for the construction of wavelet systems
    with different types of symmetry
    are actively developed (see~\cite{Bel}, \cite{ChHe}, \cite{ChLian}, \cite{DH},
    \cite{GLS}, \cite{HanM=4}, \cite{Han09}, \cite{Han10ComplSym}, \cite{HanMo},
    \cite{HanZh}, \cite{Han11lpm}, \cite{Jiang4fold}, \cite{Jiang6fold},
    \cite{KST}, \cite{Pet}, \cite{Pet03}
    and the references therein).
    In this paper, we study the construction of symmetric
    frame-like wavelet systems.
    Such systems were introduced in~\cite{KrSk}. They are obtained by the general scheme for the
    construction of MRA-based dual wavelet frames that was developed by A.~Ron and Z.~Shen~\cite{RS}.
    The realization of this scheme leads to a
    dual wavelet systems but not necessary to dual frames. To obtain frames we have to satisfy additional condition, namely
    all wavelet functions should have vanishing moments. However, perfect
    reconstruction property, frame-type decomposition and
    desired approximation properties for MRA-based dual wavelet systems
    can be achieved without frame requirements~\cite{KrSk} (also see~\cite{HanFBased}).
    This simplification of the
    construction scheme
    allows us to provide %different types of
    symmetry in general case.
%    In this paper we present a method
%    that allow to construct frame-like wavelet systems (see~\cite{KrSk})
%    with different types of symmetry and desirable approximation order.
    Let us, firstly, introduce some basic notations and definitions.

\subsection{Notations}

$\n$ is the set of positive integers,
    $\r$  is the set of real numbers,
    $\cn$ is the set of complex numbers.
    For $x~\in~\r$ we denote by $\lfloor x \rfloor$
    %an integer part of $x$.
    the largest integer not greater than $x.$
    $\rd$ denotes the
    $d$-dimensional Euclidean space,  $x = (x_1\ddd x_d)$, $y =
    (y_1\ddd y_d)$ are its elements (vectors),
    $(x)_j=x_j$ for $j~=~1,\dots,d,$
    $(x, y)~=~x_1y_1+~\dots~+x_dy_d$,
    $|x| = \sqrt {(x, x)}$, ${\bf0}=(0\ddd 0)\in\rd$.
    For $x, y \in\rd$ we write
    $x > y$ if $x_j > y_j,$ $j= 1\ddd d;$
    $x \geq y$ if $x_j \geq y_j,$ $j= 1\ddd d.$
    $\zd$ is the integer lattice
    in $\rd$, $\z_+^d:=\{x\in\zd:~x\geq~{\bf0}\}.$
    For $x\in\rd$ vector$\lfloor x \rfloor\in\zd$ and $(\lfloor x \rfloor)_j=\lfloor x_j \rfloor,$
    $j= 1\ddd d.$
    For $\,\alpha\in\zd_+\,$ we denote by $\,o(\alpha)\,$ the set of
    odd coordinates of $\,\alpha$,
    $\,\Pi_{\alpha}(\xi):=
    \prod\limits_{j=1}^d\, (1-e^{2\pi i \xi_j})^{\alpha_j}.$
    If $\alpha,\beta\in\zd_+$, $a,b\in\rd$, we set
    $[\alpha]=\sum\limits_{j=1}^d \alpha_j$,
    $\alpha!=\prod\limits_{j=1}^d\alpha_j!$,
    $\lll{\alpha}\atop{\beta}\rrr=\frac{\alpha!}{\beta!(\alpha-\beta)!}$,
    $a^b=\prod\limits_{j=1}^d a_j^{b_j}$,
    $D^{\alpha}f=\frac{\partial^{[\alpha]} f}{\partial^{\alpha_1}x_1\dots
    \partial^{\alpha_d}x_d}$, $\delta_{ab}$~is the Kronecker delta.
    For $i=1,\dots,d,$ vectors $e_i\in\zd$ are the standard basis of $\rd,$ $(e_i)_j=\delta_{ij},$ $j=1,\dots,d.$
    For $n\in\n,$ we set $\Delta_n~:=~\{\alpha:~\alpha\in\zd_+,\,\,[\alpha]<n\},$
    $\Delta^e_n~:=~\{\alpha:~\alpha\in\Delta_n, [\alpha] \texttt{\,is even}\}.$

Let $M$ be a dilation matrix, i.e. an integer $d\times d$ matrix whose
    eigenvalues are bigger than 1 in module, $M^*$ denote the conjugate
    matrix to $M$, $I_d$ denote the $d\times d$ identity matrix, $m=|\det M|.$ We say
    that vectors $k,\,n\in\zd$ are congruent modulo $M$ if $k-n=Ml,\,\, l\in\zd$.
    The integer lattice $\zd\,$ is
    split into cosets with respect to the congruence. It is known that the number of cosets
    is equal to $m$ (see, e.g., \cite[\S~2.2]{NPS}).
    %All cosets are enumerated from $0$ to $m-1$.
    Let us choose
    an arbitrary representative from each coset, call them all digits and
    denote the set of digits by $D(M)=\{s_0\ddd s_{m-1}\}$.
    %Throughout the paper we
    %assume that $m=|\det M|$,
    %$D(M)=\{s_0\ddd s_{m-1}\}$.
    Each coset we denote by $\Omega_M^i:= \{Mk+s_i, k\in\zd\},$
    $i=0\ddd m-1.$
    %digit $s_i$ belongs to the coset with the number $i$.

If $f$ is a function defined on $\rd$, we set
        $$
        f_{jk}(x):=m^{j/2}f(M^jx+k),\quad j\in\z, k\in\zd
        $$
    and $\{ f_{jk}\}:=\{f_{jk}(x):  j\in\z,\,  k\in\zd\}.$
    If function $f\in L_1(\rd)$ then its Fourier transform $\h f$
    is defined to be $\widehat
    f(\xi)=\int\limits_{\rd} f(x)e^{-2\pi i (x, \xi)}\,dx,$ $\xi\in\rd.$
    This notion can be naturally extended to $L_2(\rd)$.

The Schwartz class of functions defined on $\rd$ is denoted by $S$.
    The dual space of $S$ is $S'$, i.e.
    the space of tempered distributions.
    The basic notion and facts from distribution theory
    can be found in, e.g.,~\cite{Vladimirov-1}.
    If  $f\in S$, $g\in S'$, then
    $\overline{\langle f, g\rangle}:=\langle g, f\rangle:= g(f)$.
    If  $f\in L_p(\rd)$, $g\in L_q(\rd)$, $\frac1p+\frac1q=1$, then
    ${\langle f, g\rangle}:= \int_{\rd}f\overline g$.
    If  $f\in S'$  then $\h f$ denotes its  Fourier transform
    defined by $\langle \h f, \h g\rangle=\langle f, g\rangle$, $g\in S$.
    %Note that if  $f\in S^\prime$ is compactly supported, then,
%    due to the Paley-Wiener theorem for tempered distributions,
%    $\h f$ is an entire function.
    If  $f\in S'$, $j\in\z, k\in\zd$, we define $f_{jk}$ by
        $$
        \langle f_{jk},g\rangle= \langle f, g_{-j,-M^{-j}k}\rangle\quad \forall g\in S.
        $$
    For a function $f\in S',$ its
    Sobolev smoothness exponent is defined to be

        $$\nu_2(f)=\sup \left\{ \nu\in\r\bigcup\{-\infty\}\bigcup\{+\infty\}: \int\limits_{\rd} |\h f (\xi)|^2 (1+ |\xi|^2)^{\nu} d \xi < \infty\right\}.$$

We say that $t(\xi),$ $\xi\in\rd,$ is a trigonometric polynomial of semi-integer degrees
    associated with $\sigma\in \frac 12\zd \bigcap\left[0,\frac 12\right]^d$ %$t$
    if $t$ is a function of
    the following type
    %\be
    $$t(\xi)=\sum\limits_{k\in\zd}h_k e^{2\pi i \left(k+\sigma,\xi\right)},\quad h_k\in\cn.$$
    %\label{fTrigSemi}
    %\ee
    We say that $t(\xi),$ $\xi\in\rd,$ is a trigonometric polynomial of integer degrees
    (or just trigonometric polynomial) if $t$ is a
    trigonometric polynomial of semi-integer degrees associated with $\sigma=\nul.$

A function/distribution $\phi$ is called refinable if there exists a $1$-periodic
    function $m_0\in L_2([0,1]^d)$ (mask, also refinable mask, low-pass filter) such that
        \be
        \widehat\phi(\xi)=m_0(M^{*-1}\xi) \widehat\phi(M^{*-1}\xi).\label{RE}
        \ee
    This condition is called refinement equation.
    It is well
    known (see, e.g.,~\cite[\S~2.4]{NPS}) that for any  trigonometric polynomial
    $m_0$ satisfying $m_0({\bf0})=1$ there exists a unique solution (up to a factor) of
    the refinement equation~(\ref{RE}) in $S'$. The solution is compactly supported
    and given by its Fourier transform
        \be
        \widehat\phi(\xi):=\prod_{j=1}^\infty m_0(M^{*-j}\xi).
        \label{Prod}
        \ee
    Throughout the paper we assume that any refinable mask $m_0$ is a trigonometric
    polynomial with real or complex coefficients and $m_0({\bf0})=1$.

For any trigonometric polynomial $t$ there exists a unique set of
    trigonometric polynomials $\tau_k$, $k=0\ddd m-1$, such that
        \be
        t(\xi)= \frac1{\sqrt m}\sum\limits_{k=0}^{m-1} e^{2\pi i(s_k,\xi)}\tau_k(M^*\xi),
        \label{PR}
        \ee
    where $s_k$ are digits of the matrix $M$.
    This is so called polyphase representation of trigonometric polynomial $t$.
    For $k=0,...,m-1,$ $\tau_k$ is the polyphase component of $t$
    corresponding to the digit $s_k.$
    If we change the set of digits, the polyphase components of trigonometric polynomial $t$
    will also change.
    Any polyphase component $\tau_k$ can be expressed by
        $$
        \tau_k(\xi)=
        \frac1{\sqrt m}\sum\limits_{s\in D(M)}
        e^{-2\pi i(M^{-1}s_k,\xi+s)} t(M^{*-1}(\xi+s)),\quad k=0,...,m-1.
        $$
    If $t(\xi)=\sum_{n\in\zd}h_n e^{2\pi i (n,\xi)}, h_n\in\cn,$ then
    $$\tau_k(\xi)=\frac1{\sqrt m}\sum_{p\in\zd} h_{Mp+s_k} e^{2\pi i (p,\xi)}, \quad k=0,...,m-1.$$

For trigonometric polynomial of semi-integer degrees
     $t(\xi)=\sum\limits_{k\in\zd}h_k e^{2\pi i \left(k+\sigma,\xi\right)},$ $h_k\in\cn,$
    associated with $\sigma\in \frac 12\zd \bigcap\left[0,\frac 12\right]^d,$
    we define its coefficient support to be $\texttt{csupp}(t)=\{k\in\zd : h_k\neq 0\}$
    and the number of its non-zero coefficients to be $|\texttt{csupp}(t)|.$

For $\beta\in\Delta_n,$ $n\in\n,$ we denote by $\Theta_{\beta, n}$ the class of trigonometric polynomials of
    integer or semi-integer degrees associated with $\sigma\in \frac 12\zd \bigcap\left[0,\frac 12\right]^d$ and
    defined on $\rd$ such that $\forall t\in\Theta_{\beta, n}$%$t\in\Theta_{\beta, n}$
        \be
        D^{\alpha} t(\nul)=(2\pi i)^{\alpha}\delta_{\alpha\beta}, \quad\forall \alpha\in\Delta_n.
        \label{fG_beta}
        \ee
    %$$\Theta_{\beta, n}=\{t: D^{\alpha} t(\nul)=(2\pi i)^{\alpha}\delta_{\alpha\beta}, \alpha\in\Delta_n\}.$$
    Also, we will use the following notations:
        \begin{description}
            \item[ ] $\Theta^A_{\beta, n}:=\{t\in\Theta_{\beta, n}: t(\xi)=t(-\xi)\}$ for $\beta\in\Delta^e_n.$
            \item[ ] $\Theta^B_{\beta, n}:=\{t\in\Theta_{\beta, n}: t(\xi)=(-1)^{[\beta]}\overline{t(\xi)}\}$ for $\beta\in\Delta_n.$
        \end{description}
        How to construct trigonometric polynomials from the class $\Theta_{\beta, n}$ ($\Theta^A_{\beta, n}$ or $\Theta^B_{\beta, n}$)?
    %Those classes of trigonometric polynomials choice of functions
    %$G_{\beta}(\xi), G^{symA}_{\beta}(\xi), G^{symB}_{\beta}(\xi), \beta\in\Delta_n, \xi\in \rd$
    %can be done in different ways.
    It can be done in different ways. For example, we can take
        \be
        G_{\beta}(\xi)=\prod\limits_{r=1}^d g_{\beta_r}(\xi_r), \quad\beta\in\Delta_n,
        \label{fG_bProd}
        \ee
    %where $g_{\eta}(u)=\sum\limits_{l=0}^{n}h^{(\eta)}_l e^{ 2\pi i l u}$, $ \eta=0,\dots,n-1$, $u\in\r$
    where $g_{k}(u)\in\Theta_{k, n}$ $k=0,\dots,n-1,$ $u\in\r$. It is clear that $G_\beta(\xi)\in\Theta_{\beta, n}.$
    The coefficients of trigonometric polynomials
    $g_{k}(u)=\sum\limits_{l=0}^{n-1}h^{(k)}_l e^{ 2\pi i l u}$
    %The coefficients $h^{(\eta)}_l$
    are easy to find from the linear system of equations with Vandermonde matrix according with conditions~(\ref{fG_beta}).
     For example, if $n=2:$

        $$g_{0}(u)=1, \quad g_{1}(u)=-1+e^{2\pi i u}.$$
    If $n=3:$

        $$g_{0}(u)=1, \quad g_{1}(u)=-\frac 32 +2 e^{2\pi i u}- \frac 12 e^{4\pi i u}, \quad g_{2}(u)=\frac 12 - e^{2\pi i u}+ \frac 12 e^{4\pi i u}.$$
    If $n=4:$

         $$g_{0}(u)=1, \quad g_{1}(u)=-\frac {11} 6 + 3 e^{2\pi i u}- \frac 32 e^{4\pi i u}+\frac 13 e^{6\pi i u},$$
         $$g_{2}(u)=1 - \frac 52 e^{2\pi i u}+ 2 e^{4\pi i u}- \frac 12 e^{6\pi i u},
         \quad g_{3}(u)=-\frac 16 + \frac 12 e^{2\pi i u}- \frac 12 e^{4\pi i u}+ \frac 16 e^{6\pi i u}.$$
    Also, recursive formulas for computing $g_k(u)$ are given in~\cite{Sk6}.

Trigonometric polynomials from the class $\Theta^A_{\beta, n},$ $\beta\in\Delta^e_n,$ can be constructed
    %using $G_{\beta}(\xi)\in\Theta_{\beta, n}$
    as follows
        $$G^{A}_{\beta}(\xi)=\frac 12 (G_{\beta}(\xi)+G_{\beta}(-\xi)), \quad G_\beta(\xi)\in\Theta_{\beta, n}.$$
        It is clear that $G^{A}_{\beta}(\xi)\in\Theta^A_{\beta, n}.$

Similarly, trigonometric polynomials from the class $\Theta^B_{\beta, n},$ $\beta\in\Delta_n,$ can be constructed
    as follows

        $$G^{B}_{\beta}(\xi)=\frac 12 (G_{\beta}(\xi)+ (-1)^{[\beta]}\overline{G_{\beta}(\xi)}),\quad G_\beta(\xi)\in\Theta_{\beta, n}.$$
        It is clear that $G^{B}_{\beta}(\xi)\in\Theta^B_{\beta, n}.$

How to construct trigonometric polynomials $G_{\beta}(\xi)$
    of semi-integer degrees
    associated certain with $\sigma\in \frac 12\zd \bigcap\left[0,\frac 12\right]^d$ such that $G_{\beta}(\xi)\in\Theta_{\beta, n}$?
    We simply find $g^s_{k}(u)$
    in the following form
    $g^s_{k}(u)=\sum\limits_{l=0}^{n-1}h^{(k)}_l e^{ 2\pi i (l+\frac 12) u}$, $u\in\r,$
    such that
    $g^s_{k}(u)\in\Theta_{k, n}$ $k=0,\dots,n-1.$
    Then, as in~(\ref{fG_bProd}), we set

        $$G_{\beta}(\xi)=\prod\limits_{r\notin o(2\sigma)} g_{\beta_r}(\xi_r)
        \prod\limits_{r\in o(2\sigma)} g^s_{\beta_r}(\xi_r), \quad\beta\in\Delta_n.$$
     Trigonometric polynomials of semi-integer degrees associated with $\sigma$
     from the class $\Theta^A_{\beta, n},$ $\beta\in\Delta^e_n$
        and from the class $\Theta^B_{\beta, n},$ $\beta\in\Delta_n,$
         can be constructed as above. %by the same formulas.

For $n\in\n$ we say that trigonometric polynomial $t$ has \emph{vanishing moments of order $n$} if
        $$D^{\beta} t(\xi)\Big|_{\xi={\bf 0}}= 0, \quad \forall  \beta\in\Delta_n.$$

For $n\in\n$ we say that trigonometric polynomial $t$ has \emph{sum rule of order $n$}
    with respect to the dilation matrix $M$ if
        $$D^\beta t({M^*}^{-1}\xi)|_{\xi=s}=0, \quad \forall s\in D(M^*)\setminus\{\nul\},
        \quad \forall  \beta\in\Delta_n.$$
    High order of sum rule for a refinable mask is very
    important in applications, since it is connected with good approximation properties of the
    corresponding wavelet system (see, e.g.~\cite{Han11lpm} and the references therein).

For $n\in\n$ we say that trigonometric polynomial $t$
has \emph{linear-phase moments of order $n$} with phase $c\in\rd$ if
        $$D^{\beta} t(\nul)=
         D^{\beta} e^{2\pi i (c,\xi)}\Big|_{\xi={\bf0}}= (2\pi i)^{[\beta]} c^\beta,
         \quad \forall  \beta\in\Delta_n.$$
    This notion was described in~\cite{Han10ComplSym} for univariate case and in~\cite{Han11lpm}
    for more general settings. The importance of linear-phase moments is
    in the fact that they are useful in the sense
    of polynomial reproduction and subdivision schemes.

It should be noted that these properties of trigonometric polynomials
    are invariant with respect to the integer shifts. If trigonometric polynomial $t$
     has vanishing moments of order $n$ (sum rule of order $n$ or linear-phase moments of order $n$ with phase $c\in\rd$)
    then trigonometric polynomial $t'(\xi):=t(\xi) e^{2\pi i (\gamma,\xi)},$ $\gamma\in\zd,$ also has vanishing moments of order $n$
    (sum rule of order $n$ or linear-phase moments of order $n$ but with phase $c+\gamma$).

If $\psi^{(\nu)}\in S'$, $\nu=1\ddd r$, then
    $\{\psi_{jk}^{(\nu)}\}$ is called a wavelet system.
    For $n\in\n$ we say that compactly supported wavelet system
    $\{\psi_{jk}^{(\nu)}\}$ has
    \emph{vanishing moments of order $n$}
    (or has {\em  $VM^{n}$ property})
    if
        $$D^{\beta}\h{\psi^{(\nu)}}({\bf0})=0, \quad \nu=1\ddd r, \quad \forall \beta\in\Delta_n.$$

A family of functions $\{f_n\}_{n\in\aleph}$ ($\aleph$ is a countable index set) in a
    Hilbert space $H$ is called a frame in $H$
    if there exist $A, B > 0$ so that, for all $f\in H,$
        $$A\|f\|^2\le\sum\limits_{n\in\aleph}|\langle f, f_n\rangle|^2\le B \|f\|^2.$$
    An important property of a frame $\{f_n\}_n$ in $H$
    is the following: every $f\in H$ can be decomposed as
    $ f=\sum_n\langle f,\w f_n\rangle f_n, $ where $\{\w f_n\}_n$ is a
    dual frame in $H$. Comprehensive characterizations of frames
    can be found in~\cite{Chris}. Wavelet frames are of great interest for many applications, especially
    in signal processing.

\subsection{Symmetry property of mask and refinable function}

A finite set of $d\times d$ integer matrices $H\subset\{E: |\det E|=1\}$ is a
    symmetry group on $\zd$ if $H$ forms a group under matrix multiplication.

Let $H$ be a symmetry group on $\zd$.
    A function $f$ is said to be \emph{$H$-symmetric} with respect to a
    center $C\in \rd$ if
        $$f(E(x-C)+C)=f(x), \quad\forall E\subset H, \quad x\in\rd.$$
    Hence, a trigonometric polynomial
        $t(\xi)=\sum\limits_{n\in\zd}h_n e^{2\pi i (n, \xi)},\, h_n \in\cn$
    is \emph{$H$-symmetric }
    %($H$-antisymmetric)}
    with respect to a center $c\in \rd$ if
        \be
        t(\xi)=e^{2\pi i (c-Ec,\xi)}t(E^*\xi), \quad
        %\left(t(\xi)=-e^{2\pi i (c-Ec,\xi)}t(E^*\xi)\right) \quad
        \forall E\subset H
        \label{TrigSym}
        \ee
    and $c-Ec\in\zd, \forall E\subset H.$
    Condition~(\ref{TrigSym}) is equivalent to

        $$h_n=h_{E(n-c)+c}, %\quad (h_n=-h_{E(n-c)+c}),
        \quad \forall n\in\zd, \quad\forall E\subset H.$$

It is known (see, e.g.~\cite{HanSym02}) that $H$-symmetry of a mask $m_0(\xi)$ not always carries over
    to its refinable function $\phi$, defined by~(\ref{Prod}). Due to this fact
    a stronger notion of symmetry was introduced in~\cite{HanSym02}.

A finite set of $d\times d$ integer matrices $H_M$ is said to be
    \emph{symmetry group with respect to the dilation matrix $M$} if
    $H_M$ is a symmetry group on $\zd$ such that
    $$MEM^{-1}\in H_M,\quad \forall E\subset H_M.$$

The following statement was shown by B.~Han.

\begin{lem}~\cite[Proposition 2.1]{Han3}
    Let $H_M$ be a symmetry group with respect to
    the dilation matrix M. Mask $m_0$ is $H_M$-symmetric with respect to
    a center $c\in\rd$ if and only if the corresponding refinable function $\phi$ given by~(\ref{Prod})
    is $H_M$-symmetric with respect to a center $C\in \rd$
%        $$ \phi (E(\cdot-C)+C)=\phi \quad \forall E\subset G_M,$$
    and $C=(M-I_d)^{-1}c.$

    \label{HanSym}
    \end{lem}

Let $H$ be a symmetry group on $\zd$ and
    trigonometric polynomial $t(\xi)$ be $H$-symmetric with a center $c\in \frac 12 \zd$.
    When $H=\{I_d, -I_d\}$, condition~(\ref{TrigSym}) is equivalent to

        \be
        t(\xi)=e^{2\pi i (2c,\xi)}t(-\xi) \quad\texttt{or}\quad h_n=h_{2c-n},\,\,\forall k\in\zd.
        \label{fSymMinus}
        \ee
    Such trigonometric polynomial $t(\xi)$ is called
    \emph{symmetric with respect to the point $c$.} Moreover, in this case $H$ is a
    symmetry group with respect to any dilation matrix and Lemma~\ref{HanSym} for $H$ is always valid.

If mask $m_0(\xi)$ is symmetric with respect to the point $c$
    then

        $$D^{e_j}m_0(\nul)=D^{e_j}\left(e^{2\pi i (2c,\xi)}m_0(-\xi)\right)\Big|_{\xi={\bf0}}=
        2\pi i 2(c)_j - D^{e_j}m_0(\nul).$$
    Therefore, $D^{e_j}m_0(\nul)=2\pi i c^{e_j},\, j=1,\dots,d,$
    i.e. the phase for linear-phase moments
    must match with the symmetry point $c$. The maximal order of linear-phase moments for symmetric mask
    with respect to the point $c$ must be an even
    integer (see~\cite[Proposition 2]{Han11lpm}).

Let $H$ be a symmetry group on $\z^2$ and
    trigonometric polynomial $t(\xi)$ be $H$-symmetric with a center $c\in \frac 12\z^2$.
    When
        $$H=\left\{\pm I_2, \pm \left(\begin{matrix}
            -1 & 0\cr
            0 & 1\cr
        \end{matrix}\right) \right\},$$
    $t(\xi)$ is called \emph{axial symmetric with respect to the center $c$.}
    When
        $$H=\left\{\pm I_2, \pm \left(\begin{matrix}
            -1 & 0\cr
            0 & 1\cr
        \end{matrix}\right), \pm \left(\begin{matrix}
            0 & 1\cr
            -1 & 0\cr
        \end{matrix}\right),
        \pm \left(\begin{matrix}
            0 & 1\cr
            1 & 0\cr
        \end{matrix}\right) \right\},$$
    $t(\xi)$ is \emph{called 4-fold symmetric with respect to the center $c$.}

In fact, for trigonometric polynomial $t(\xi)$ with complex Fourier coefficients, we can consider
    another type of point symmetry
        \be
        t(\xi)=e^{2\pi i (2c,\xi)}\overline{t(\xi)}
        \quad\texttt{or}\quad h_{k}=\overline{h_{2c-k}},\,\, \forall k\in\zd,\,\, c\in \frac 12 \zd.
        \label{fSymOver}
        \ee
    A mask with point symmetry in the sense~(\ref{fSymOver}) is often called
    a linear-phase filter in engineering.
    Linear phase filter has constant group delay, all frequency components
    have equal delay times.
    That is, there is no distortion due to the time delay of frequencies
    relative to one another. In many applications this is useful.
    If mask $m_0$ is a linear-phase filter
    then for the corresponding refinable function $\phi$ defined by~(\ref{Prod}) equality
    $\phi(2C-x)= \overline{\phi(x)}$ holds, where $C=(M-I_d)^{-1}c$ (see, e.g.,~\cite[Lemma 2]{Han10ComplSym}).

If trigonometric polynomial $t(\xi)$ has real Fourier coefficients, %i.e. $h_k\in\r$,
    then two symmetries in the sense~(\ref{fSymMinus}) and in the sense~(\ref{fSymOver})
    are the same.

A symmetry is one the most desirable property for wavelet systems. A lot of different approaches
    providing symmetry
    have been developed since the history of wavelets began. In~\cite{Daub}, I.~Daubechies
    showed that for dyadic dilation $M=2$ the Haar function $\chi_{[0,1]}$, up to an integer shift, is the only
    orthogonal refinable function with symmetry and compact support. %proved
    %Consequently, there is no compactly supported real-valued symmetric orthonormal dyadic wavelet
    %basis in dimension one which consists of continuous functions.
    Thus, to provide symmetry and compact support for orthogonal wavelet bases, we have to consider dilations $M>2$ or
    multivariate case. For dilation factor $M=3$ C.~Chui and J.~A.~Lian~\cite{ChLian} considered the construction of symmetric and antisymmetric orthonormal
    wavelets. B.~Han~\cite{HanM=4}
    constructed symmetric orthonormal scaling functions and wavelets with dilation factor $M=4.$ The method for construction of orthogonal symmetric refinable functions for
    arbitrary dilation factors $M>2$ and any order of sum rule is suggested by E.~Belogay and Y.~Wang~\cite{Bel}. For any given orthogonal symmetric refinable function algorithm for
    deriving orthogonal wavelet bases was found by A.~Petukhov~\cite{Pet} for arbitrary dilation factors $M>2$.
    The case of orthogonal symmetric refinable masks with complex coefficients,
    linear-phase moments, order of sum rule and step-by-step algorithm for construction high pass filters  was investigated by
    B. Han~\cite{Han09},~\cite{Han10ComplSym},~\cite{Han11lpm}.

For a given symmetric refinable function, C.~Chui and W.~He~\cite{ChHe} proved
    that there exists symmetric tight wavelet frame with 3 generators for dyadic dilation $M=2$.
    In some cases the number of generators can be reduced to 2~\cite{Pet03}.
    For dilation factor $M=2$ symmetric tight wavelet frames
    with 3 generators and high vanishing moments was constructed by B.~Han, Q.~Mo~\cite{HanMo}.
    For an arbitrary dilation factor $M\ge2$ systematic algorithm for constructing (anti)symmetric tight wavelet frames
    generated by a given refinable function was presented in~\cite{Pet}.

Starting from any two symmetric compactly supported refinable functions in $L_2(\r)$ with dilation
    factor $M\ge2$, I.~Daubechies and B.~Han~\cite{DH} showed that it is always possible to construct $2M$ wavelet functions with compact
    support such that they generate a pair of (anti)symmetric dual wavelet frames in $L_2(\r)$.

In multivariate case, S.~S.~Goh, Z.~Y.~Lim, Z.~Shen~\cite{GLS} suggested a general and simple method for the
    construction of point symmetric tight wavelet frames from a given tight wavelet frame and the number of generators
    for the symmetric tight wavelet frames is at most twice the number of generators for the given tight wavelet frame.
    The order of vanishing moments is preserved. S.~Karakaz’yan, M.~Skopina, M.~Tchobanou in~\cite{KST} described all
     real interpolatory masks which are symmetric with respect to the origin and generate symmetric/antisymmetric compactly supported
     biorthonormal bases or dual wavelet systems, generally speaking,
    with vanishing moments up to arbitrary
    order $n$ for matrix dilations whose determinant is odd or equal $\pm 2$.

The construction of highly symmetric wavelet systems is the least studied theme in the literature.
    One of method was suggested by Q.~Jiang in~\cite{Jiang4fold},~\cite{Jiang6fold} for the construction
    of dual wavelet frames for multiresolution surface processing with 4-fold symmetry for dyadic and $\sqrt 2$-refinement
    and 6-fold symmetry for dyadic and $\sqrt 3$-refinement.

\subsection{Paper outline}

The paper is organized as follows.  In Section 2 we introduce mixed extension principle
    and the notion of MRA-based frame-like wavelet system. Section 3 is devoted to
    the construction of point symmetric refinable masks
    and point symmetric/antisymmetric frame-like wavelet system
    with desired approximation properties. In Section 4 we construct axial symmetric
    frame-like wavelet system
    with desired approximation properties.
    In Section 5 some results about the construction of highly symmetric
    frame-like wavelet system are given.
    In Section 6 several examples are presented.

\section{Preliminary results}

%    The goal of this paper is to propose the approach for construction
%    of compactly supported MRA-based dual wavelet systems
%    $\{\psi_{jk}^{(\nu)}\}$, $\{\widetilde{\psi}_{jk}^{(\nu)}\}$, where
%    wavelet functions
%    $\psi^{(\nu)}, \widetilde{\psi}^{(\nu)}$ are functions or distributions with
%    different types of symmetry
%    that provide frame-type decomposition and desired approximation order.

A general scheme for the construction of MRA-based  wavelet frames
    was developed in~\cite{RS},
    \cite{RS2} (Mixed Extension Principle).
    To construct a pair of dual wavelet frames in $L_2(\rd)$ one
    starts with two refinable functions
    $\phi, \widetilde\phi$, $\h\phi(\nul)=1$, $\h{\w\phi}(\nul)=1,$
    (or its masks $m_0, \widetilde m_0$, respectively),
    then finds wavelet masks $m_{\nu}, \widetilde m_{\nu}$, $\nu = 1\ddd r$,
    $r\ge m-1$, which are also trigonometric
    polynomials such that matrices
        \ban
        \cal M:=\left(\begin{array}{ccc}
        \mu_{00}&\hdots & \mu_{0, m-1}\cr
        \vdots & \ddots & \vdots \cr
        \mu_{r,0}& \hdots & \mu_{r,m-1}
        \end{array}\right),\,
        \cal \widetilde M:=\left(\begin{array}{ccc}
        \widetilde{\mu}_{00}&\hdots & \widetilde{\mu}_{0, m-1}\cr
        \vdots & \ddots & \vdots \cr
        \widetilde{\mu}_{r,0}& \hdots & \widetilde{\mu}_{r,m-1}
        \end{array}\right)
        \ean
    satisfy
        \be
        {\cal M}^T\overline{\cal\widetilde M}= I_m
        \label{calM},
        \ee
    where $\mu_{\nu k},$ $\w\mu_{\nu k}$ $k=0,\dots,m-1,$ are the polyphase components of
    trigonometric polynomials $m_{\nu},$ $\w m_{\nu}$ for all $\nu=0,\dots,r,$ $r\ge m-1.$

Wavelet functions $\psi^{(\nu)}$, $\w\psi^{(\nu)},$
    $\nu = 1\ddd r$, $r\ge m-1,$ are defined by their Fourier transform
        \be \widehat{\psi^{(\nu)}}(\xi)=
        m_\nu(M^{*-1}\xi)\widehat\varphi(M^{*-1}\xi), \quad
        \widehat{\widetilde\psi^{(\nu)}}(\xi)=
        \widetilde m_\nu(M^{*-1}\xi)\widehat{\widetilde\varphi}(M^{*-1}\xi).
        \label{PsiDef}
        \ee

If wavelet functions  $\psi^{(\nu)}, \widetilde\psi^{(\nu)}$,
    $\nu=1\ddd r$, $r\ge m-1,$
    are constructed as above, then
    the set of functions $\{\psi_{jk}^{(\nu)}\}$, $\{\w\psi_{jk}^{(\nu)}\}$
    is said to be
    (MRA-based) dual wavelet systems generated by
    refinable functions $\varphi, \widetilde\varphi$
    (or their masks $m_0, \widetilde m_0$).

For wavelet functions $\w\psi^{(\nu)}$, $\nu=1\ddd r,$ $r\ge m-1,$ defined by~(\ref{PsiDef})
    $VM^{n}$ property for the corresponding
    wavelet system is equivalent to the fact that wavelet masks $\w m_\nu$, $\nu=1\ddd r$
    have vanishing moments up to order $n$.
%        \be
%        D^\beta(\w m_\nu({M^*}^{-1}x))\Big|_{x=\nul}=0,\quad \nu=1\ddd r,\quad
%        \forall \beta\in\zd_+, [\beta]\le n.
%        \label{09}
%        \ee

For a compactly supported dual wavelet systems $\{\psi_{jk}^{(\nu)}\}$, $\{\w\psi_{jk}^{(\nu)}\},$
    $\nu=1\ddd r$, $r\ge m-1,$
    generated by
    refinable functions $\varphi, \widetilde\varphi\in L_2(\rd)$
    %with $\phi, \widetilde\phi, \psi^{(\nu)}, \w\psi^{(\nu)}\in L_2(\rd)$
    necessary (see~\cite[Theorem 1]{Sk1}) and
    sufficient (see~\cite[Theorems 2.2, 2.3]{Han1}) condition
    to be a pair of dual wavelet frames in $L_2(\rd)$ is the following:
    each wavelet function $\psi^{(\nu)}, \widetilde\psi^{(\nu)}$
    should have vanishing moment at least of order $1$.
    Construction of univariate dual wavelet frames
    with vanishing moments was developed, in~\cite{DH}.
    For multidimensional case explicit method for construction
    of compactly supported dual wavelet frames with vanishing moments was suggested in~\cite{Sk1},
    but it is not easy to provide various types of symmetry for different matrix dilations.

If we reject frame requirements and aim to provide
    $VM^{n}$ property only for wavelet functions $\w\psi^{(\nu)}$, $\nu=1\ddd r$, then
    the method could be symplified. Necessary conditions with constructive proof are given in

\begin{lem} \cite[Lemma 14]{KrSk} Let  $\phi, \widetilde\phi$ be compactly supported  refinable distributions,
    $\mu_{0 k},\ \widetilde\mu_{0 k}$, $k=0\ddd m-1$,
    be the polyphase components of their masks $m_0, \w m_0$ respectively, $n\in\n$. If
        \ba
        D^\beta\mu_{0k}(\nul)=\frac {(2\pi i)^{[\beta]}} {\sqrt m}
        \sum\limits_{\alpha\in\zd_+,\,\alpha\le\beta}\lambda_\alpha
        \lll\beta\atop\alpha\rrr(- M^{-1}s_k)^{\beta-\alpha} \quad  \forall
        \beta\in\Delta_n,\,\, \forall k=0,\dots,m-1
        \label{19}
        \ea
    for some  numbers  $\lambda_\alpha \in \cn$,
    $\alpha\in\Delta_n$, $\lambda_\nul=1$, and
        \be
        D^\beta\lll 1-\sum_{k=0}^{m-1}\mu_{0k}(\xi)
        \overline{\w\mu_{0k}(\xi)}\rrr\Bigg|_{\xi=\nul}=0\ \quad  \forall
        \beta\in\Delta_n,
        \label{20}
        \ee
    then there exist (MRA-based) dual wavelet systems
    $\{\psi_{jk}^{(\nu)}\}$, $\{\w\psi_{jk}^{(\nu)}\}$, $\nu=1\ddd m,$
    such that $\{\widetilde\psi^{(\nu)}_{jk}\}$
    has $VM^{n}$ property.
   \label{lemKrSk}
\end{lem}

Later, in Theorem~\ref{theoWaveA} the extension technique
    and explicit formulas for the wavelet masks will be given. It is important to say that frames
    cannot be constructed using this method (see Remark~\ref{rNotFrame} for more details).

Note that (see, e.g.,~\cite{Sk1}) the set of numbers $\lambda_{\alpha}$ is unique, does not
    depend on $n$ and

        \ba
        \lambda_{\alpha} =\frac 1 {(2\pi i)^{[\alpha]}}
        D^{\alpha} m_0(M^{\ast-1}\xi)\Big|_{\xi={\bf0}},\quad \alpha\in\Delta_n.
        \label{Lambda}
        \ea
    Conditions~(\ref{19}),  (\ref{20}) do not depend on our choice of the set of digits.

\begin{rem} Due to Theorem~\cite[Theorem 10]{DSk} and the identity
        $$
        \sum_{k=0}^{m-1}\mu_{0k}(M^*\xi)\overline{\w\mu_{0k}}(M^*\xi)=
        \sum_{s\in D(M^{*})}m_0(\xi+{M}^{*-1}s)\overline{\w m_0(\xi+{M}^{*-1}s)},
        $$
    conditions~(\ref{19}),  (\ref{20}) in Lemma~\ref{lemKrSk} may be replaced by

              \hspace{3.4cm} (i) $m_0(\xi)$ has sum rule of order $n$,

                    \be
                    (ii)\ \ D^\beta\lll 1-m_0(\xi)\overline{\w m_0(\xi)}\rrr\Big|_{\xi=0}=0\ \
                    \forall \beta\in\Delta_n.
                    \label{20_new}
                    \ee

    \label{r2}
    \end{rem}

%    It is well known that in the case $r=m-1$, (\ref{09}) depends only on
%    the refinable mask $\widetilde m_0$, and does not depend
%    on matrix extension (which is not unique). However, it is not true if $r>m-1$,
%    in this case (\ref{09}) depends on $m_0,\widetilde m_0$ and matrix extension.
%    A criterion is given in

Thus, the method for
    the construction of MRA-based frame-like wavelet systems is quite simple. We start
    with any mask $m_0$ that
    has sum rule of order $n$.
    %condition~(\ref{19}) is satisfied.
    Next, we find dual refinable mask $\w m_0$ satisfying~(\ref{20_new}) and
    wavelet masks $m_\nu, \w m_\nu$, $\nu=1\ddd m$ using Lemma~\ref{lemKrSk}.
    The refinable functions $\phi, \w \phi$ are defined by~(\ref{Prod});
    wavelet functions $\psi^{(\nu)}, \w\psi^{(\nu)}$ are defined by~(\ref{PsiDef}).
    Thereby, we get the pair of MRA-based dual wavelet systems
    $\{\psi_{jk}^{(\nu)}\}$, $\{\w\psi_{jk}^{(\nu)}\}$
    %which are not necessary frames in $L_2(\rd)$
    %and, moreover,
    which may consist of tempered distributions, generally speaking.
    Nevertheless, we can consider frame-type decomposition
    with respect to these systems for appropriate functions $f$
    (e.g.,  $f\in S$ if $\widetilde\psi^{(\nu)}\in S'$;
    $f\in L_p$ if $\widetilde\psi^{(\nu)}\in L_q$, $\frac1p+\frac1q=1$).

Let $\{\psi_{jk}^{(\nu)}\}$, $\{\w\psi_{jk}^{(\nu)}\}$  be dual
    wavelet systems and  $A$ be a class of functions $f$  for which
    $\langle f,\widetilde\phi_{0k}\rangle$,
    $\langle f,\widetilde\psi_{jk}^{(\nu)}\rangle$ have meaning.
    We say that  $\{\psi_{jk}^{(\nu)}\}$ is {\em frame-like} if
        \be
        f=\sum\limits_{j=-\infty}^{\infty}\sum\limits_{\nu=1}^r\sum\limits_{k\in\,\zd}
        \langle f,\widetilde\psi_{jk}^{(\nu)}\rangle\psi_{jk}^{(\nu)}\quad \forall f\in A,
        \label{03}
        \ee
    and  $\{\psi_{jk}^{(\nu)}\}$  is {\em almost frame-like} if
        \be
        f=\sum\limits_{k\in\,\zd}\langle f,\widetilde\phi_{0k}\rangle\phi_{0k}+
        \sum\limits_{j=0}^{\infty}\sum\limits_{\nu=1}^r\sum\limits_{k\in\,\zd}
        \langle f,\widetilde\psi_{jk}^{(\nu)}\rangle\psi_{jk}^{(\nu)}\quad \forall f\in A,
        \label{02}
        \ee
    where the series in~(\ref{03}) and (\ref{02}) converge in some natural sense.

Here we give some results from~\cite{KrSk}.

\begin{theo} \cite[Theorem 12]{KrSk} Let $f\in S$, $\phi,\widetilde\phi\in S'$,
        $\phi, \widetilde\phi$ be compactly supported and refinable,
        $\h\phi(\nul)=\h{\widetilde\phi}(\nul)=1$.
        Then MRA-based dual wavelet systems $\{\psi_{jk}^{(\nu)}\}$, $\{\w\psi_{jk}^{(\nu)}\},$
         $\nu=1\ddd r,$
        generated by $\phi, \widetilde\phi$
        are almost frame-like, i.e.~(\ref{02}) holds
        with the series converging in  $S'$.
    \label{t6}
    \end{theo}

\begin{theo} \cite[Theorem 16]{KrSk}
    { Let $f\in S$, $\phi\in L_2(\rd)$, $\widetilde\phi\in
    S'$, $\phi, \widetilde\phi$ be compactly supported and refinable,
    $\h\phi(\nul)=\h{\widetilde\phi}(\nul)=1$, and let
    $\psi^{(\nu)}, \widetilde\psi^{(\nu)}, \nu=1\ddd r,$ be defined by~(\ref{PsiDef}).  Then
     (\ref{02}) holds with the series converging in $L_2$-norm.
    If, moreover, $\phi, \widetilde\phi$ are as in Lemma~\ref{lemKrSk}, then
        \ba
        \left\|f-\sum\limits_{k\in\,\zd}\langle
        f,\widetilde\phi_{0k}\rangle \phi_{0k}-
        \sum\limits_{i=0}^{j-1}\sum\limits_{\nu=1}^m\sum\limits_{k\in\,\zd}
        \langle f,\widetilde\psi_{ik}^{(\nu)}\rangle\psi_{ik}^{(\nu)}\right\|_2
        \le
        \frac{C\|f\|_{W_2^{n^*}}}{(|\lambda|-\varepsilon)^{jn}},
        \label{6}
        \ea
     where  $\lambda$ is a minimal (in module) eigenvalue of  $M$,
    $\varepsilon>0$, $|\lambda|-\varepsilon>1$,  $n^*\ge n$,
     C  and $n^*$ do not depend on $f$ and $j$.}
     \label{t7}
\end{theo}
In other words, almost frame-like wavelet system $\{\psi_{jk}^{(\nu)}\}$ from
    Theorem~\ref{t7}
    provide approximation order $n$.

\section{Construction of symmetric frame-like wavelets}

In this section, we will consider refinable masks
    $m_0(\xi)=\sum_{n\in\zd}h_n e^{2\pi i (n,\xi)}, $
    $h_n\in\cn$ that are symmetric with respect to a point $c\in\frac 12 \zd$
    in two different senses

        \begin{description}
          \item[(a)] $m_0(\xi)=e^{2\pi i (2c,\xi)}m_0(-\xi)$,
                    that is $h_{k}=h_{2c-k}, \forall k\in\zd$, then we say that
                    $m_0(\xi)$ is symmetric with respect to
                    the point $c\in\frac 12 \zd$ in the sense (a);
          \item[(b)] $m_0(\xi)=e^{2\pi i (2c,\xi)}\overline{m_0(\xi)}$,
                    that is $h_{k}=\overline{h_{2c-k}}, \forall k\in\zd$,
                    then we say that
                    $m_0(\xi)$ is symmetric with respect to
                    the point $c\in\frac 12 \zd$ in the sense (b).
        \end{description}

Similarly, $m_0(\xi)$ is antisymmetric with respect to
                    the point $c\in\frac 12 \zd$ in the sense (a), if  $m_0(\xi)=-e^{2\pi i (2c,\xi)}m_0(-\xi);$
            $m_0(\xi)$ is antisymmetric with respect to
                    the point $c\in\frac 12 \zd$ in the sense (b), if $m_0(\xi)=-e^{2\pi i (2c,\xi)}\overline{m_0(\xi)}.$

We will use the polyphase components $\mu_{0k},$ $k=0,\dots,m-1,$ for the construction of
    mask $m_0$ that is symmetric with respect to a point. Thus, we need to reformulate
    symmetry conditions for the mask in terms of its polyphase components.

\begin{lem}
    A mask $m_0(\xi)$ is symmetric with respect to
    a point $c\in\frac 12 \zd$ in the sense (a)
    if and only if
    its polyphase components $\mu_{0k}(\xi), $ $k=0,\dots,m-1$ satisfy

        \be
        \mu_{0k}(\xi)=e^{2\pi i (M^{-1}(2c-s_k-s_l),\xi)} \mu_{0l}(-\xi),
        \label{PolyTypeA}
        \ee
    where $l$
    is a unique integer in $\{0,\dots,m-1\}$
    %in a one way so that for the corresponding digits $s_k$ and $s_l$ we have
    such that
    $M^{-1}(2c-s_k-s_l)\in\zd.$

\label{PolyLem}
\end{lem}

{\bf Proof.} Firstly, we assume that $m_0(\xi)$ is symmetric
    with respect to the point $c\in\frac 12 \zd$ in the sense (a), i.e.
    $h_n=h_{2c-n},$ $\forall n\in\zd.$ Namely, all coefficients $h_n$
    such that indices
     $n\in\Omega_M^k$ match
    with the coefficients $h_{2c-n}$
    whose indices
     $2c-n\in\Omega_M^l$ for some $l\in\{0,\dots,m-1\}$ and

        $$e^{2\pi i (s_k,\xi)}\mu_{0k}(M^*\xi)=\sum_{n\in\zd}h_{Mn+s_k}e^{2\pi i (Mn+s_k,\xi)}=
        \sum_{n\in\zd}h_{2c-Mn-s_k}e^{2\pi i (Mn+s_k,\xi)}$$
        $$=\sum_{p\in\zd}h_{Mp+s_l}e^{2\pi i (2c-Mp-s_l,\xi)}=
        e^{2\pi i (2c-s_l,\xi)}\mu_{0l}(-M^*\xi),$$
    where $k=0,\dots,m-1$, and number $l$
    is chosen such that $2c-Mn-s_k\in \Omega_M^l$, i.e. $2c-Mn-s_k=Mp+s_l, p\in\zd$,
    or $M^{-1}(2c-s_k-s_l)\in\zd.$

Next, suppose we have~(\ref{PolyTypeA}). According with~(\ref{PR}) we obtain

        $$ m_0(-\xi)=
        \frac1{\sqrt m}\sum\limits_{k=0}^{m-1} e^{-2\pi i(s_k,\xi)}\mu_{0k}(-M^*\xi)=$$
        $$e^{-2\pi i(2c,\xi)}\frac1{\sqrt m}\sum\limits_{l=0}^{m-1} e^{2\pi i (s_l,\xi)}\mu_{0l}(M^*\xi)=
        e^{-2\pi i(2c,\xi)}m_0(\xi),
        $$
    i.e. mask $m_0$ is symmetric with respect to the point $c$ in the sense (a).
$\Diamond$

This fact for the univariate case and for multiwavelets is hidden in the proof of~\cite[Lemma 1]{HanZh}.
    So, we decided to give a direct proof of this Lemma for multivariate case.

The same statement is valid for mask $m_0$ that is symmetric
    with respect to the point $c\in\frac 12 \zd$ in the sense (b)
    with condition~(\ref{PolyTypeA}) replaced by
        \be
        \mu_{0k}(\xi)=e^{2\pi i (M^{-1}(2c-s_k-s_l),\xi)} \overline{\mu_{0l}(\xi)}.
        \label{PolyTypeB}
        \ee

The set of indices $\{0,\dots,m-1\}$ according to the equation~(\ref{PolyTypeA}) or~(\ref{PolyTypeB})
    can be divided into sets $I$ and $J$.

Index $i\in I$, if we can find index $i'\neq i$ so that $M^{-1}(2c-s_i-s_{i'})\in\zd.$
    Then index $i'\in I$ too and

        $$
        \mu_{0i}(\xi)=e^{2\pi i (M^{-1}(2c-s_i-s_{i'}),\xi)}\mu_{0i'}(-\xi),\quad M^{-1}(2c-s_i-s_{i'})\in\zd$$
    or
        $$ \quad \mu_{0i}(\xi)=e^{2\pi i (M^{-1}(2c-s_i-s_{i'}),\xi)}\overline{\mu_{0i'}(\xi)}
            \quad
        M^{-1}(2c-s_i-s_{i'})\in\zd.
        $$
    The corresponding digits $s_i, s_{i'}$ always can be choose so that

        \be
        M^{-1}(2c-s_i-s_{i'})=0.
        \label{fDigitSi}
        \ee
    Indeed, if $M^{-1}(2c-s_i-s_{i'})=n, n\in\zd,$ we simply replace digit $s_{i'}$ by $s_{i'}+Mn$
    and get the required equality.
    Then the polyphase component $\mu_{0i'}$ also change and we get %$\mu_{0i},  are said to be flipped version of each other and

        \be
        \mu_{0i}(\xi)=\mu_{0i'}(-\xi).
        \label{ItypeA}
        \ee
    or
        \be
        \mu_{0i}(\xi)=\overline{\mu_{0i'}(\xi)}.
        \label{ItypeB}
        \ee
    In what follows in this section, we assume that digits are chosen such that~(\ref{fDigitSi}) is valid.
    Also, we will use notations for the separation of the set $I$ into two as follows
    $I=I^1\cup I^2, I^1\cap I^2=\emptyset$ and for index $i\in I^1$
    the corresponding index
    $i'\in I^2$ and $i'\notin I^1.$

Index $j\in J$, if the corresponding polyphase component $\mu_{0j}$
    is symmetric with respect to the point $M^{-1}(c-s_j)$, i.e.
        \be
        \mu_{0j}(\xi)=e^{2\pi i (M^{-1}(2c-2s_j),\xi)}\mu_{0j}(-\xi),\quad
        M^{-1}(2c-2s_j)\in\zd.
        \label{IItypeA}
        \ee
    or
        \be
        \mu_{0j}(\xi)=e^{2\pi i (M^{-1}(2c-2s_j),\xi)}\overline{\mu_{0j}(\xi)},\quad
        M^{-1}(2c-2s_j)\in\zd.
        \label{IItypeB}
        \ee
    The cardinality of the set $J$ is no more then $2^d.$ Indeed,
    the digits always can be choose so that $M^{-1}s_k\in [0,1)^d, k=0,\dots,m-1.$
    The vector $2c$ can be represented as $2c=Mn+s_p,$ $n\in\zd,$ $s_p\in D(M).$
    Therefore, we should have $M^{-1}(Mn+s_p-2s_j)\in\zd$ or $2 M^{-1}s_j\in \zd+M^{-1}s_p.$
    Since $2 M^{-1}s_j\in [0,2)^d$, there is only $2^d$ possible variants for  $2 M^{-1}s_j$.

\subsection{Construction of symmetric initial mask $m_0$}

Using the above consideration, Lemma~\ref{lemKrSk} and Remark~\ref{r2} we suggest a simple
    algorithm for the construction of mask $m_0$ that is
    symmetric with respect to
    the point $c\in\frac 12 \zd$ in the sense (a) or in the sense (b) and
    has arbitrary order $n$ of sum rule (or equivalently condition~(\ref{19}) is fulfilled with some
    numbers $\lambda_{\alpha}\in \cn,$ $\alpha\in\Delta_n$).
    Also we indicate how to provide
    linear-phase moments for $m_0$. Firstly, we need the following simple

\begin{lem}
        Let $t'(\xi)$ be a trigonometric polynomial, $n\in\n,$
            $t(\xi):=e^{2\pi i (a,\xi)}t'(\xi),$ $a\in\rd.$
        Then
            $D^{\beta}t'(\xi)\Bigg|_{\xi=\nul}=(2\pi i)^{[\beta]}\kappa'_{\beta}, $
             $\forall \beta\in\Delta_n$
        if and only if

            $$D^{\beta}t(\xi)\Bigg|_{\xi=\nul}=(2\pi i)^{[\beta]}\sum\limits_{\alpha\in\zd_+,\,\alpha\le\beta}
            \kappa'_{\alpha} \lll\beta\atop\alpha\rrr a^{\beta-\alpha}, \quad  \forall \beta\in\Delta_n.$$

    \label{LambdaLem}
\end{lem}

{\bf Proof.} Let
        $D^{\beta}t'(\xi)\Bigg|_{\xi=\nul}=(2\pi i)^{[\beta]}\kappa'_{\beta}$, $\forall \beta\in\Delta_n.$
    Then the statement is obvious.

    Conversely, let us denote
        $\kappa_\beta:=\sum\limits_{\alpha\in\zd_+,\,\alpha\le\beta} \kappa'_{\alpha}
        \lll\beta\atop\alpha\rrr a^{\beta-\alpha},$ $\forall \beta\in\Delta_n$
    and
        $D^{\beta}t(\xi)\Bigg|_{\xi=\nul}=(2\pi i)^{[\beta]}\kappa_{\beta}.$
    Then

        $$D^{\beta}t'(\xi)\Bigg|_{\xi=\nul}=(2\pi i)^{[\beta]}\sum\limits_{\alpha\in\zd_+,\,\alpha\le\beta}
        (-a)^{\beta-\alpha} \lll\beta\atop\alpha\rrr \kappa_{\alpha}=$$$$
        (2\pi i)^{[\beta]}\sum\limits_{\alpha\in\zd_+,\,\alpha\le\beta} (-a)^{\beta-\alpha} \lll\beta\atop\alpha\rrr
        \sum\limits_{\gamma\in\zd_+,\, \gamma\le\alpha} \kappa'_{\gamma} \lll\alpha\atop\gamma\rrr a^{\alpha-\gamma},
        \quad \forall \beta\in\Delta_n.$$
    After the change of variables in the sums and taking into account that

        $$ \lll\alpha\atop{\gamma}\rrr \lll\beta\atop\alpha\rrr=
        \frac{\alpha!}{\gamma!(\alpha-\gamma)!} \frac{\beta!}{\alpha!(\beta-\alpha)!} =$$
        $$\frac {\beta! (\beta-\gamma)!} {\gamma!(\alpha-\gamma)! (\beta-\gamma-(\alpha-\gamma))!(\beta-\gamma)!}=
        \lll\beta\atop\gamma\rrr \lll\beta-\gamma\atop\alpha-\gamma\rrr$$
     we have

        $$D^{\beta}t'(\xi)\Bigg|_{\xi=\nul}=
        (2\pi i)^{[\beta]}\sum\limits_{\gamma\in\zd_+,\, \gamma\le\beta}\kappa'_{\gamma} \sum\limits_{\alpha\in\zd_+\atop\gamma\le\alpha\le\beta}
        (-a)^{\beta-\alpha} \lll\beta\atop\alpha\rrr \lll\alpha\atop\gamma\rrr a^{\alpha-\gamma}=$$
        $$(2\pi i)^{[\beta]}\sum\limits_{\gamma\in\zd_+,\,\gamma\le\beta}\kappa'_{\gamma} \lll\beta\atop\gamma\rrr
        \sum\limits_{\alpha\in\zd_+\atop\gamma\le\alpha\le\beta} (-a)^{\beta-\alpha}
        \lll\beta-\gamma\atop\alpha-\gamma\rrr a^{\alpha-\gamma}=$$
        $$(2\pi i)^{[\beta]}\sum\limits_{\gamma\in\zd_+,\,\gamma\le\beta}\kappa'_{\gamma}
        \lll\beta\atop\gamma\rrr (a-a)^{\beta-\gamma} = (2\pi i)^{[\beta]}\kappa'_{\beta},$$
        where $\beta\in\Delta_n.$ $\Diamond$

The construction of mask $m_0(\xi)$ will be carried out using its polyphase components
    $\mu_{0k}(\xi),$ $k=0,\dots,m-1.$ So,
    we reformulate condition~(\ref{19}) for mask $m_0$ that is symmetric in the sense (a) or in the sense (b).
\begin{lem}
    Let $n\in\n,$ $m_0$ be a mask satisfying condition~(\ref{19}) with some
    numbers $\lambda_{\alpha}\in \cn,$ $\alpha\in\Delta_n$.
    If $m_0$ is symmetric with
    respect to the point $c\in\frac 12\zd$ in the sense (a) then
    \begin{enumerate}
      \item  numbers $\lambda_{\alpha}\in \cn$ are given by
      \be
      \lambda_{\alpha}=
        \sum\limits_{\gamma\in\zd_+,\,\gamma\le\alpha} \lambda'_{\gamma}
        \lll\alpha\atop\gamma\rrr (M^{-1}c)^{\alpha-\gamma}, \quad\forall \alpha\in\Delta_n
        \label{lam'}
        \ee
        where $\lambda'_{\alpha}\in\cn$, $\alpha\in\Delta_n$, $\lambda_{\nul}=1,$
      $\lambda'_{\alpha}=0$, if $[\alpha]$ is odd;

      \item  condition~(\ref{19}) is equivalent to
        \be
        D^\beta\mu_{0k}(\nul)=\frac {(2\pi i)^{[\beta]}} {\sqrt m} \sum\limits_{\gamma\in\zd_+\atop \gamma\le\beta}
        \lambda'_{\gamma} \lll\beta\atop\gamma\rrr (M^{-1}c-M^{-1}s_k)^{\beta-\gamma},
        \quad \forall \beta\in\Delta_n,\, \forall k=0,\dots,m-1.
        \label{19mod}
        \ee
    \end{enumerate}
        If $m_0$ is symmetric with
    respect to the point $c\in\frac 12\zd$ in the sense (b) then
    \begin{enumerate}
      \item   numbers $\lambda_{\alpha}\in \cn$ are given by~(\ref{lam'})
      where $\lambda'_{\alpha}\in\cn$, $\alpha\in\Delta_n$, $\lambda_{\nul}=1,$
      $Re \lambda'_{\alpha}=0,$ if $[\alpha]$ is odd,
      $Im\lambda'_{\alpha}=0,$ if $[\alpha]$ is even;
      \item  condition~(\ref{19}) is equivalent to~(\ref{19mod}).
    \end{enumerate}

    \label{LambdaLemA}
\end{lem}
{\bf Proof.} Let $m'_0(\xi)=e^{-2\pi i (c,\xi)}m_0(\xi).$
     Therefore, $m'_0(\xi)$ is an even trigonometric polynomial of
    semi-integer degrees associated with $\sigma=c-\lfloor c \rfloor$. Define numbers $\lambda'_{\alpha},$ $\alpha\in\Delta_n,$ by
        $$\lambda'_{\alpha}:=\frac 1 {(2\pi i)^{[\alpha]}}
        D^{\alpha} m'_0(M^{\ast-1}\xi)\Big|_{\xi={\bf0}}.$$
    It is clear that $\lambda'_{\nul}=1,$  $\lambda'_{\alpha}=0$, if $[\alpha]$ is odd.
            For mask $m_0(\xi)=e^{2\pi i (c,\xi)} m'_0(\xi)$ due to Lemma~\ref{LambdaLem} and~(\ref{Lambda})

        $$\lambda_{\alpha}=\frac 1 {(2\pi i)^{[\alpha]}}D^{\alpha}m_0(M^{*-1}\xi)\Bigg|_{\xi=\nul}=
        \sum\limits_{\gamma\in\zd_+,\,\gamma\le\alpha} \lambda'_{\gamma}
        \lll\alpha\atop\gamma\rrr (M^{-1}c)^{\alpha-\gamma}, \quad\forall\alpha\in\Delta_n.$$
    After combining~(\ref{lam'}) with~(\ref{19}) and using the same reasoning %considerations
    as in Lemma's~\ref{LambdaLem} proof we %get the required conditions~(\ref{19mod})
    establish the equivalence of conditions~(\ref{19mod}) and~(\ref{19}).
    %on the polyphase components $\mu_{0k}, $ $ k=0\ddd m-1.$

    For point symmetry in the sense (b), the proof is similar.
    $\Diamond$

\begin{theo}
    Let $M$ be an arbitrary matrix dilation, $c\in\frac 12\zd,$ $n\in\n$,
    $\lambda'_{\alpha}\in\cn$, $\alpha\in\Delta_n$, $\lambda'_{\alpha}=0$, if $[\alpha]$ is odd.
    Then there exists mask $m_0$
    which is symmetric with
    respect to the point $c$ in the sense (a) and
    satisfies condition~(\ref{19})
    with numbers $\lambda_\alpha$ defined by~(\ref{lam'}).
    The mask $m_0$ can be represented by %looks as follows

        \ba
        m_0(\xi)=\frac 1 m \sum\limits_{j\in J} \sum\limits_{\beta\in\Delta_n}
        \lambda'_{\beta} G^{A}_{\beta j}(M^*\xi)e^{2\pi i (c,\xi)}+
        \hspace{7cm}
        \nonumber
        \\
        \frac 1 m \sum\limits_{i\in I^1}\sum\limits_{\beta\in\Delta_n}
        \sum\limits_{\gamma\in\zd_+\atop\gamma\le\beta}\lambda'_{\gamma} \lll\beta\atop\gamma\rrr (M^{-1}(c-s_i))^{\beta-\gamma}
        \left(G_{\beta i}(M^*\xi) e^{2\pi i (s_i,\xi)}+
        G_{\beta i}(-M^*\xi) e^{2 \pi i (2c-s_i,\xi)}\right),
        \label{m_0A}
        \ea
    where
    $G_{\beta i}(\xi)\in\Theta_{\beta,n},$ $\beta\in\Delta_n,$ $i\in I^1,$ are trigonometric polynomials;
    $G^{A}_{\beta j}(\xi)\in\Theta^A_{\beta, n},$ $\beta\in\Delta^e_n,$ $j\in J,$ are
    trigonometric polynomials of semi-integer degrees associated with $\sigma=M^{-1}(c-s_j)-\lfloor M^{-1}(c-s_j)\rfloor.$

    \label{theoMaskA}
\end{theo}

{\bf Proof.} %The construction of mask $m_0(\xi)$ will be carried out using the polyphase components
%    $\mu_{0k}(\xi),$ $k=0,\dots,m-1.$ Firstly,
%    we reformulate condition~(\ref{19}) for symmetric in the sense (a) mask $m_0$.
%    Let us consider an even mask $m'_0(\xi)$
%    %i.e. $m'_0(\xi)=m'_0(-\xi)$
%    and denote
%        $$\lambda'_{\alpha}:=\frac 1 {(2\pi i)^{[\alpha]}}
%        D^{\alpha}m'_0(M^{*-1}\xi)\Bigg|_{\xi=\nul},
%        \quad \alpha\in\Delta_n.$$
%    Obviously, $\lambda_{\nul}=1,$ $\lambda'_{\alpha}=0$, if $[\alpha]$ is odd.
%
%    For the trigonometric polynomial $m_0(\xi)=e^{2\pi i (c,\xi)} m'_0(\xi)$ due to Lemma~\ref{LambdaLem}
%
%        $$\lambda_{\alpha}=D^{\alpha}m_0(M^{*-1}\xi)\Bigg|_{\xi=\nul}=
%        \sum\limits_{\gamma\in\zd_+,\,\gamma\le\alpha} \lambda'_{\gamma}
%        \lll\alpha\atop\gamma\rrr (M^{-1}c)^{\alpha-\gamma}, \quad\alpha\in\Delta_n.$$
%
%    After combining $\lambda_{\alpha}$ with~(\ref{19}) and using the same reasoning %considerations
%    as in Lemma's~\ref{LambdaLem} proof we get the following conditions
%    on the polyphase components $\mu_{0k}, $ $ k=0\ddd m-1$
%    %$\gamma\in\zd_+$,  $[\gamma]\le n$, $\lambda_\nul=1$, such that
%        \be
%        D^\beta\mu_{0k}(\nul)=\frac {(2\pi i)^{[\beta]}} {\sqrt m} \sum\limits_{\gamma\in\zd_+,\, \gamma\le\beta}
%        \lambda'_{\gamma} \lll\beta\atop\gamma\rrr (M^{-1}c-M^{-1}s_k)^{\beta-\gamma},
%        \quad \forall \beta\in\Delta_n.
%        \label{19mod}
%        \ee
Let us construct the polyphase components  $\mu_{0k}(\xi),$ $k=0,\dots,m-1$ such that
they satisfy conditions~(\ref{PolyTypeA}) and~(\ref{19mod}). Therefore,
we will get the required mask $m_0$.

For the polyphase components $\mu_{0i},$ $\mu_{0i'}$, $i\in I^1,$ $2c-s_i-s_{i'}=0,$
    we set

        \ba
        \mu_{0i}(\xi)= \sum\limits_{\beta\in\Delta_n}
        \left[ \frac 1 {(2\pi i)^{[\beta]}} D^\beta\mu_{0k}(\nul)\right]
        G_{\beta i}(\xi)= \hspace{6cm}
        \nonumber
        \\
        \hspace{3cm}\frac {1} {\sqrt m}\sum\limits_{\beta\in\Delta_n}
        \left[\sum\limits_{\gamma\in\zd_+,\, \gamma\le\beta}\lambda'_{\gamma} \lll\beta\atop\gamma\rrr (M^{-1}c-M^{-1}s_i)^{\beta-\gamma} \right]
        G_{\beta i}(\xi)
        \label{mu_0iA}
        \ea
    and

        \be
        \mu_{0i'}(\xi)=\mu_{0i}(-\xi),
        \label{mu_0i'A}
        \ee
    where $G_{\beta i}(\xi)\in\Theta_{\beta,n},$  $\beta\in\Delta_n,$ are trigonometric polynomials.
    The condition~(\ref{19mod}) for $\mu_{0i}$ is obviously valid.
    Since $ c - s_i = -(c - s_{i'}) , $ and
    $\lambda'_{\alpha}=0$, when $[\alpha]$ is odd, we have

        $$D^{\beta}\mu_{0i'}({\bf 0}) =
        \frac {(2\pi i)^{[\beta]}} {\sqrt m} \sum\limits_{\alpha\in\zd_+,\, \alpha\le\beta\atop [\alpha]\, \texttt{is \,even}}
        \lambda'_{\alpha} \lll\beta\atop\alpha\rrr (M^{-1}c-M^{-1}s_{i'})^{\beta-\alpha}=\hspace{2.2cm}$$
        $$ \hspace{1.5cm}(-1)^{[\beta]}\frac {(2\pi i)^{[\beta]}} {\sqrt m}
        \sum\limits_{\alpha\in\zd_+,\, \alpha\le\beta\atop [\alpha]\, \texttt{is \,even}}
        \lambda'_{\alpha} \lll\beta\atop\alpha\rrr (M^{-1}c-M^{-1}s_i)^{\beta-\alpha}=(-1)^{[\beta]}D^{\beta}\mu_{0i}({\bf 0}).$$
    Therefore~(\ref{19mod}) for $\mu_{0i'}$ is also satisfied.
    %because $D^{\alpha}G_{\beta i}(-\xi)\Bigg|_{\xi=\nul}= (-2\pi i)^{[\beta]}\delta_{\alpha\beta}.$

For the polyphase components $\mu_{0j}$, $j\in J,$
    we have to satisfy conditions~(\ref{IItypeA}) and~(\ref{19mod}).
    Let us
     %consider $\mu'_{0j}(\xi)=e^{-2\pi i (M^{-1}c-M^{-1}s_j, \xi)}\mu_{0j}(\xi),$
    %and require that
    construct an even trigonometric polynomial $\mu'_{0j}$ of semi-integer degrees associated with
     $\sigma=M^{-1}(c-s_j)-\lfloor M^{-1}(c-s_j) \rfloor$
    %, i.e. $\mu'_{0j}(\xi)=\mu'_{0j}(-\xi)$
    such that
    %due to
    %Lemma~\ref{LambdaLem} must hold
    $D^\beta\mu'_{0j}(\nul)=\frac {(2\pi i)^{[\beta]}} {\sqrt m}\lambda'_{\beta},$
    $\forall \beta\in\Delta_n.$
To provide it we set

    $$\mu'_{0j}(\xi)=\frac 1 {\sqrt m}\sum\limits_{\beta\in\Delta_n} \lambda'_{\beta} G^{A}_{\beta j}(\xi),$$
    where $G^{A}_{\beta j}\in\Theta^A_{\beta, n},$ $\beta\in\Delta^e_n,$ are
    trigonometric polynomials of semi-integer degrees associated with $\sigma=M^{-1}(c-s_j)-\lfloor M^{-1}(c-s_j) \rfloor.$

    %$G^{A}_{\beta j}(\xi)$ has semi-integer degrees by variables $\xi_r$, where $r\in o(2M^{-1}(c-s_j)).$
        Therefore,

        \ba
        \mu_{0j}(\xi)=e^{2\pi i (M^{-1}c-M^{-1}s_j, \xi)}\frac 1 {\sqrt m} \sum\limits_{\beta\in\Delta_n} \lambda'_{\beta} G^{A}_{\beta j}(\xi).
        \label{mu_0jA}
        \ea
    Condition~(\ref{19mod}) is satisfied due to Lemma~\ref{LambdaLem}. Condition~(\ref{IItypeA})
    is obviously valid.

Hence, we construct the polyphase components $\mu_{0k}$, $k=0,\dots,m-1$ that
    satisfy symmetry conditions~(\ref{PolyTypeA}) in Lemma~\ref{PolyLem}
    %symmetric with respect to the point $c$ in the sense (a) by Lemma~\ref{PolyLem}
    and condition~(\ref{19mod}) with $\lambda'_{\alpha}$.
    %Parameters $\lambda'_\gamma=1$, if $\gamma=\nul$, $\lambda'_\gamma=0$, if $\gamma$ is odd and $\lambda'_\gamma\in\cn$, for else.
    To get the mask $m_0$ it remains to use formula~(\ref{PR}) together with~(\ref{mu_0iA}),~(\ref{mu_0i'A}),~(\ref{mu_0jA})
    and take into account that $s_{i'}=2c-s_i.$ $\Diamond$

\begin{theo}
    Let $M$ be an arbitrary matrix dilation, $c\in\frac 12\zd,$ $n\in\n$,
    $\lambda'_{\alpha}\in\cn$, $\alpha\in\Delta_n$, $Re \lambda'_{\alpha}=0,$ if $[\alpha]$ is odd,
    $Im\lambda'_{\alpha}=0,$ if $[\alpha]$ is even.
    Then there exists mask $m_0$ which is
    symmetric with
    respect to the point $c$ in the sense (b) and satisfies condition~(\ref{19})
     with numbers $\lambda_\alpha$ defined by~(\ref{lam'}).
    The mask $m_0$ can be represented by

        \ba
        m_0(\xi)=\frac 1 m \sum\limits_{j\in J} \sum\limits_{\beta\in\Delta_n}
        \lambda'_{\beta} G^{B}_{\beta j}(M^*\xi)e^{2\pi i (c,\xi)}+
        \hspace{6cm}
        \nonumber
        \\
        \frac 1 m \sum\limits_{i\in I^1}\left[\sum\limits_{\beta\in\Delta_n}
        \left(\sum\limits_{\gamma\in\zd_+,\,\gamma\le\beta}\lambda'_{\gamma} \lll\beta\atop\gamma\rrr (M^{-1}c-M^{-1}s_i)^{\beta-\gamma} \right)
        G_{\beta i}(M^*\xi) e^{2\pi i (s_i,\xi)}
        +
        \right.
        \hspace{1cm}
        \nonumber
        \\
        \left.\sum\limits_{\beta\in\Delta_n}
         \left(\sum\limits_{\gamma\in\zd_+,\,\gamma\le\beta}\overline{\lambda'_{\gamma}} \lll\beta\atop\gamma\rrr (M^{-1}c-M^{-1}s_i)^{\beta-\gamma} \right)
        \overline{G_{\beta i}(M^*\xi)} e^{2 \pi i (2c-s_i,\xi)}\right],
        \label{m_0B}
        \ea
    where $G_{\beta i}(\xi)\in\Theta_{\beta,n},$ $\beta\in\Delta_n,$ $i\in I^1,$ are trigonometric polynomials;
    $G^{B}_{\beta j}(\xi)\in\Theta^B_{\beta, n},$ $\beta\in\Delta_n,$ $j\in J,$ are
    trigonometric polynomials of semi-integer degrees associated with $\sigma=M^{-1}(c-s_j)-\lfloor M^{-1}(c-s_j) \rfloor.$

    \label{theoMaskB}
\end{theo}

{\bf Proof.}
The construction of mask $m_0(\xi)$ is very similar to Theorem's~\ref{theoMaskA} proof.
    Let us construct the polyphase components  $\mu_{0k}(\xi),$ $k=0,\dots,m-1$ such that
    they satisfy conditions~(\ref{PolyTypeB}) and~(\ref{19mod}). Therefore,
    we will get the required mask $m_0$.
%Firstly,
%    we reformulate condition~(\ref{19}) for symmetric in the sense (b) mask $m_0$.
%    Let us consider a real mask $m'_0(\xi)$
%    %, i.e. $m'_0(\xi)=\overline{m'_0(\xi)}$
%    and
%    denote
%        $$\lambda'_{\alpha}:=\frac 1 {(2\pi i)^{[\alpha]}}D^{\alpha}m'_0(M^{*-1}\xi)\Bigg|_{\xi=\nul},
%        \quad\alpha\in\Delta_n$$
%       % \quad\lambda'_{\alpha}=\lambda'^r_{\alpha}+i\lambda'^i_{\alpha} ,$$
%     Obviously, $\lambda'_\nul=1$,
%      $Re \lambda'_{\alpha}=0,$ if $[\alpha]$ is odd, $Im\lambda'_{\alpha}=0,$ if $[\alpha]$ is even.
%
%    For the trigonometric polynomial $m_0(\xi)=e^{2\pi i (c,\xi)} m'_0(\xi)$ due to Lemma~\ref{LambdaLem}
%    the conditions on the polyphase components $\mu_{0k}, $ $ k=0\ddd m-1$ remain the same as in~(\ref{19mod}).

For the polyphase components $\mu_{0i},$ $\mu_{0i'}$, $i\in I^1,$ $2c-s_i-s_{i'}=0,$
    we set
        \ba
        \mu_{0i}(\xi)=\sum\limits_{\beta\in\Delta_n}\left[\frac 1 {(2\pi i )^{[\beta]}}
        D^\beta\mu_{0i}(\nul) \right] G_{\beta i}(\xi)= \hspace{6cm}
        \nonumber
        \\
        \frac 1 {\sqrt m}\sum\limits_{\beta\in\Delta_n}
        \left[\sum\limits_{\gamma\in\zd_+,\,\gamma\le\beta}
        \lambda'_{\gamma} \lll\beta\atop\gamma\rrr (M^{-1}c-M^{-1}s_k)^{\beta-\gamma} \right] G_{\beta i}(\xi)
        \label{mu_0iB}
        \ea
    and
        \be
        \mu_{0i'}(\xi)=\overline{\mu_{0i}(\xi)},
        \label{mu_0i'B}
        \ee
    where $G_{\beta i}(\xi)\in\Theta_{\beta,n},$ $\beta\in\Delta_n,$ are trigonometric polynomials.
    The condition~(\ref{19mod}) for $\mu_{0i}$ is obviously valid.
    Since $ c - s_i = -(c- s_{i'}), $ we have

        $$D^{\beta}\mu_{0i'}({\bf 0}) =
        \frac {(2\pi i )^{[\beta]}} {\sqrt m} \sum\limits_{\alpha\in\zd_+,\,\alpha\le\beta}
        \lambda'_{\alpha} \lll\beta\atop\alpha\rrr (M^{-1}c-M^{-1}s_{i'})^{\beta-\alpha}=$$

        $$\frac {(2\pi i )^{[\beta]}} {\sqrt m} \sum\limits_{\alpha\in\zd_+,\,\alpha\le\beta \atop [\alpha]\, \texttt{is even}}
        Re\lambda'_{\alpha} \lll\beta\atop\alpha\rrr (M^{-1}c-M^{-1}s_{i'})^{\beta-\alpha}+$$
        $$\frac {(2\pi i )^{[\beta]}} {\sqrt m}\sum\limits_{\alpha\in\zd_+,\,\alpha\le\beta\atop[\alpha]\, \texttt{is odd}}
        i Im\lambda'_{\alpha} \lll\beta\atop\alpha\rrr (M^{-1}c-M^{-1}s_{i'})^{\beta-\alpha}=$$

        $$ (-1)^{[\beta]}\frac {(2\pi i )^{[\beta]}} {\sqrt m} \sum\limits_{\alpha\in\zd_+\alpha\le\beta}
        \overline{\lambda'_{\alpha}} \lll\beta\atop\alpha\rrr (M^{-1}c-M^{-1}s_i)^{\beta-\alpha}=
        (-1)^{[\beta]}\overline{D^{\beta}\mu_{0i}({\bf 0})}.$$
    Therefore~(\ref{19mod}) for $\mu_{0i'}$ is also satisfied.
    %, because $D^{\alpha} \overline{G_{\beta i}(\xi)}\Bigg|_{\xi=\nul}= (-2\pi i)^{[\beta]}\delta_{\alpha\beta}.$

For the polyphase components $\mu_{0j}$, $j\in J$ we have to satisfy condition~(\ref{IItypeB}) and~(\ref{19mod}).
    Let us
    construct a real trigonometric polynomial $\mu'_{0j}$ of semi-integer degrees associated with
     $\sigma=M^{-1}(c-s_j)-\lfloor M^{-1}(c-s_j) \rfloor$
    %consider $\mu'_{0j}(\xi)=e^{-2\pi i (M^{-1}c-M^{-1}s_j, \xi)}\mu_{0j}(\xi)$
    %and require that $\mu'_{0j}(\xi)=\overline{\mu'_{0j}(\xi)}$ and
    %due to
    %Lemma~\ref{LambdaLem} must  hold
    such that
    $D^\beta\mu'_{0j}(\nul)=\frac {(2\pi i)^{[\beta]}}{\lambda'_{\beta}},$ $\forall \beta\in\Delta_n.$
To provide it we set

    $$\mu'_{0j}=\sum\limits_{\beta\in\Delta_n} \lambda'_{\beta} G^{B}_{\beta j}(\xi),$$
    where $G^{B}_{\beta j}\in\Theta^B_{\beta, n},$ $\beta\in\Delta_n,$ are
     trigonometric polynomials of semi-integer degrees associated with $\sigma=M^{-1}(c-s_j)-\lfloor M^{-1}(c-s_j) \rfloor.$
        Therefore,

        \ba
        \mu_{0j}(\xi)=e^{2\pi i (M^{-1}c-M^{-1}s_j, \xi)}
        \frac 1 {\sqrt m}\sum\limits_{\beta\in\Delta_n} \lambda'_{\beta} G^{B}_{\beta j}(\xi).
        \label{mu_0jB}
        \ea
    Condition~(\ref{19mod}) is satisfied due to Lemma~\ref{LambdaLem}. Condition~(\ref{IItypeB})
    is obviously valid.

Hence, we construct the polyphase components $\mu_{0k}$, $k=0,\dots,m-1$ that
    satisfy symmetry conditions~(\ref{PolyTypeB})
    %symmetric with respect to the point $c$ in the sense (a) by Lemma~\ref{PolyLem}
    and condition~(\ref{19mod}) with $\lambda'_{\alpha}$.
    %Parameters $\lambda'_\gamma=1$, if $\gamma=\nul$, $\lambda'_\gamma=0$, if $\gamma$ is odd and $\lambda'_\gamma\in\cn$, for else.
    To get the mask $m_0$ it remains to use formula~(\ref{PR}) together with~(\ref{mu_0iB}),~(\ref{mu_0i'B}),~(\ref{mu_0jB})
    and take into account that $s_{i'}=2c-s_i.$
$\Diamond$

\begin{rem}
    To provide linear-phase moments of order $n$, $n\in\n,$ for the mask $m_0$
    we simply have to take $\lambda'_{\beta}=\delta_{\beta\nul},$ $\beta\in\Delta_n$ in
    Theorem~\ref{theoMaskA}.
    Therefore, due to Lemma~\ref{LambdaLemA}

        $$D^{\beta}m_0(M^{*-1}\xi)\Bigg|_{\xi=\nul}= (2\pi i)^{[\beta]}(M^{-1}c)^{\beta} \quad \forall \beta\in\Delta_n,$$
    that means that $m_0$ has
     linear-phase moments of order $n$.

According to~(\ref{m_0A}), mask $m_0$ that
    is symmetric
     with respect to the point $c\in\frac 12 \zd,$ has sum rule of order $n$ and
     linear-phase moments of order $n,$
     %if $n$ is even, or $n+1$, if $n$ is odd,
     looks as follows

   $$
    m_0(\xi)=\frac 1 {\sqrt m} \sum\limits_{j\in J} G^{A}_{\nul j}(M^*\xi)e^{2\pi i (c,\xi)}+
    \hspace{7cm}
    $$$$
    \frac 1 {\sqrt m} \sum\limits_{i\in I^1}\sum\limits_{\beta\in\Delta_n}
         (M^{-1}c-M^{-1}s_i)^{\beta}\left(G_{\beta i}(M^*\xi) e^{2\pi i (s_i,\xi)}+
        G_{\beta i}(-M^*\xi) e^{2 \pi i (2c-s_i,\xi)}\right),
   $$
     where $G_{\beta i}(\xi)\in\Theta_{\beta,n},$ $\beta\in\Delta_n,$ $i\in I^1$ are trigonometric polynomials;
       $G^{A}_{\nul j}\in\Theta^A_{\nul, n},$ $j\in J,$ are
       trigonometric polynomials of semi-integer degrees associated with $\sigma=M^{-1}(c-s_j)-\lfloor M^{-1}(c-s_j) \rfloor.$

    \label{rLpm}
\end{rem}

\begin{rem}
    The construction of mask $m_0$ with real coefficients which is
    symmetric with
    respect to the point $c$ and satisfies condition~(\ref{19}) can be done as in
    Theorem~\ref{theoMaskA} by~(\ref{m_0A}) but $\lambda'_{\alpha}\in\r$
    %$\alpha\in\zd_+,$ $[\alpha]< n$, $\lambda'_{\alpha}=0$, if $[\alpha]$ is odd,
    and  $G_{\beta i}(\xi),$ $G^{A}_{\beta j}(\xi)$ for $\beta\in\Delta_n,$ $j\in J,$ $i\in I^1$
    are the same trigonometric polynomial as in Theorem~\ref{theoMaskA} but with real coefficients.
    \label{rRealCoef}
\end{rem}

\begin{rem} Let $m_0$ be a mask which is symmetric with
    respect to the point $c$ and satisfies condition~(\ref{19}). Such masks are constructed in
    Theorem~\ref{theoMaskA} or  Theorem~\ref{theoMaskB}.
    How to get mask $m_0$ that has the minimal number of its non-zero coefficients
    among the other masks with the same properties?
    It is enough to choose functions
    $G_{\beta}(\xi)\in\Theta_{\beta, n},$  $\beta\in\Delta_n,$
    $G^{A}_{\beta}(\xi)\in\Theta^A_{\beta, n},$ $\beta\in\Delta^e_n,$ or
    $G^{B}_{\beta}(\xi)\in\Theta^B_{\beta, n},$ $\beta\in\Delta_n,$ $\xi\in \rd$
    so that they have the minimal possible number of their non-zero coefficients.
    This numbers are equal to the quantities of the imposed requirements on the
    coefficients of these functions.
    For instance, functions $G_{\beta}(\xi)\in\Theta_{\beta, n}$ have to satisfy
    $D^{\alpha}G_{\beta}({\bf 0})=
    (2\pi i)^{[\beta]}\delta_{\alpha\beta},$ $\forall \alpha\in\Delta_n.$
    We will seek these functions in the following form $G_{\beta}(\xi)=\sum\limits_{k\in\zd}h^{(\beta)}_k e^{ 2\pi i \left(k+\sigma,\xi\right)},$
    where $\sigma\in \frac 12\zd \bigcap\left[0,\frac 12\right]^d,$
    such that $\left|\texttt{csupp}\,\left(\sum\limits_{\beta\in\Delta_n}G_{\beta}\right)\right|$ is equal to the cardinality of the set $\Delta_n.$
    According to the requirements $D^{\alpha}G_{\beta}({\bf 0})=
    (2\pi i)^{[\beta]}\delta_{\alpha\beta},$ $\alpha\in\Delta_n,$ we get the linear system of equations, and if
    coefficients support is chosen in appropriate way, then the matrix of the linear system is invertible. For example,
    the matrix is always invertible if
    $\texttt{csupp}\,\left(\sum\limits_{\beta\in\Delta_n}G_{\beta}\right)=\Delta_n$ (for more information see, e.g.~\cite{Gasca}).
    The functions $G^{B}_{\beta}(\xi)$ with the minimal number of their non-zero coefficients
    can be chosen as

        $$G^{B}_{\beta}(\xi)=\frac 12 (G_{\beta}(\xi)+ (-1)^{[\beta]}\overline{G_{\beta}(\xi)}),\quad\forall\beta\in\Delta_n.$$

    To get functions  $G^{A}_{\beta}(\xi),$ $\beta\in\Delta^e_n$ with the minimal number of their non-zero coefficients,
    we fistly look for $G_{\beta}(\xi)$ such that $D^{\alpha}G_{\beta}({\bf 0})=
    (2\pi i)^{[\beta]}\delta_{\alpha\beta},$  $\forall \alpha\in\Delta_n^e.$
    These functions can be found %in the following form $G^{A}_{\beta}(\xi)=\sum\limits_{k\in\zd}h^{(\beta)}_k e^{ 2\pi i (k,\xi)}$
    as above from the linear system of equations
    such that $\left|\texttt{csupp}\,\left(\sum\limits_{\beta\in\Delta^e_n}G_{\beta}\right)\right|$ is equal to the cardinality of the set $\Delta^e_n.$
    Then, as previously stated, we set
    $$G^{A}_{\beta}(\xi)=\frac 12 (G_{\beta}(\xi)+G_{\beta}(-\xi)),\quad\forall\beta\in\Delta_n^e.$$

    \label{rMinimSup}
\end{rem}

\subsection{Class of symmetric initial masks}

To describe the whole class of masks which are symmetric with respect to a point in the sense (a) or in the sense (b) and
    has arbitrary order of sum rule, %satisfy condition~(\ref{19}) with numbers $\lambda_\alpha$ as in~(\ref{lam'}),
    we need some auxiliary statements.

\begin{lem}
    Let $n \in \n$, $\sigma\in \frac 12\zd \bigcap\left[0,\frac 12\right]^d.$ A general form for all even (odd) trigonometric
    polynomials $T(\xi)$ of semi-integer degrees associated with $\sigma$ %, i.e. $T(\xi)=T(-\xi)$ ($T(\xi)=-T(-\xi)$)
    such that
    $D^{\beta} T (\nul) = 0, \forall\beta\in\Delta_n$, is given by
    \be
    T(\xi)=\sum\limits_{\alpha\in\zd_+ \atop [\alpha]=n}
    \left(A_{\alpha}(\xi)Re \Pi_{\alpha}(\xi)+i B_{\alpha}(\xi) Im \Pi_{\alpha}(\xi)\right),
    \label{EOpolA}
    \ee
    %$$ i\left( B_{\alpha}(\xi) Re \Pi_{\alpha}(\xi)+ \w A_{\alpha}(\xi) Im \Pi_{\alpha}(\xi)\right),$$
    where $A_{\alpha}(\xi)$, $B_{\alpha}(\xi)$ are arbitrary respectively even (odd) and odd (even)
    trigonometric polynomials of semi-integer degrees associated with $\sigma$.
    \label{GenSymA}
\end{lem}

{\bf Proof.} Let $T(\xi)$ defined by~(\ref{EOpolA}). For $\alpha\in\zd_+, [\alpha]=n$ it follows from
    $$D^{\beta} \Pi_{\alpha}(\nul)=D^{\beta}\left( \prod_{j=1}^d \left(1-e^{2\pi i \xi_j}\right)^{\alpha_j}\right)\Big|_{\xi=\nul}=0,$$
    that

    $$D^{\beta}( Re\Pi_{\alpha}(\xi))\Big|_{\xi=\nul}=0, \quad
    D^{\beta}( Im\Pi_{\alpha}(\xi))\Big|_{\xi=\nul}=0\quad \forall\beta\in\Delta_n.$$
    Then due to Leibniz formula, we obtain that $D^{\beta} T (\nul) = 0$ for all
    $\beta\in\Delta_n$. Trigonometric polynomials $\Pi_{\alpha}(\xi)$
    have real Fourier coefficients, then (see, e.g.~\cite[Lemma 4]{KST}) $Re\Pi_{\alpha}(\xi)$ is an even function and $Im\Pi_{\alpha}(\xi)$ is
    an odd function. Therefore, it is clear that $T(\xi)$ is an even (odd) trigonometric polynomial of semi-integer degrees associated with $\sigma$.

Now let $T(\xi)$ be an even trigonometric polynomial of semi-integer degrees associated with $\sigma$, such that
    $D^{\beta} T (\nul) = 0, \forall\beta\in\Delta_n$.
    By Taylor formula, there exist trigonometric polynomials $T_{\alpha}(\xi)$
    of semi-integer degrees associated with $\sigma,$ where
    $\alpha\in\zd_+,$ $[\alpha]=n,$ so  that

        $$T(\xi)=\sum\limits_{\alpha\in\zd_+ \atop [\alpha]=n} \Pi_{\alpha}(\xi) T_{\alpha}(\xi)=\hspace{9cm}$$
        $$\sum\limits_{\alpha\in\zd_+ \atop [\alpha]=n}
        \left(A_{\alpha}(\xi)Re \Pi_{\alpha}(\xi)+ B_{\alpha}(\xi)Re \Pi_{\alpha}(\xi)
        +i B_{\alpha}(\xi) Im \Pi_{\alpha}(\xi)+iA_{\alpha}(\xi)Im \Pi_{\alpha}(\xi)\right),$$
where $A_{\alpha}(\xi)= \frac {T_{\alpha}(\xi)+T_{\alpha}(-\xi)} 2,$
    $B_{\alpha}(\xi)= \frac {T_{\alpha}(\xi)-T_{\alpha}(-\xi)} 2,$ $\alpha\in\zd_+, [\alpha]=n.$
    It is clear that $A_{\alpha}(\xi)$ are even and $B_{\alpha}(\xi)$ are odd trigonometric polynomials of semi-integer degrees associated with $\sigma$.
    Rewrite
    the previous equation as

        $$T(\xi)-
        \sum\limits_{\alpha\in\zd_+ \atop [\alpha]=n}
        \left(A_{\alpha}(\xi)Re \Pi_{\alpha}(\xi)+ i B_{\alpha}(\xi) Im \Pi_{\alpha}(\xi)\right)
        =\sum\limits_{\alpha\in\zd_+ \atop [\alpha]=n}
        \left(B_{\alpha}(\xi)Re \Pi_{\alpha}(\xi)+iA_{\alpha}(\xi)Im \Pi_{\alpha}(\xi)\right).$$

The left hand side is an even trigonometric polynomial of semi-integer degrees associated with $\sigma$ and
    the right hand side is an odd trigonometric polynomial
    of semi-integer degrees associated with $\sigma$.
    Hence both the sides are identically equal to zero. Therefore,~(\ref{EOpolA}) holds.

For odd $T(\xi)$, the proof is similar.
$\Diamond$

\begin{lem}
    Let $n \in \n$, $\sigma\in \frac 12\zd \bigcap\left[0,\frac 12\right]^d.$ A general form for all real (purely imaginary) trigonometric
polynomials $T(\xi)$ of semi-integer degrees associated with $\sigma$ such that
    $D^{\beta} T (\nul) = 0, \forall\beta\in\Delta_n$, is given by
    \be
    T(\xi)=\sum\limits_{\alpha\in\zd_+ \atop [\alpha]=n}
    \left(\w A_{\alpha}(\xi)Re \Pi_{\alpha}(\xi)+\w B_{\alpha}(\xi) Im \Pi_{\alpha}(\xi)\right),
    \label{EOpolB}
    \ee
    %$$ i\left( B_{\alpha}(\xi) Re \Pi_{\alpha}(\xi)+ \w A_{\alpha}(\xi) Im \Pi_{\alpha}(\xi)\right),$$
    where $\w A_{\alpha}(\xi)$, $\w B_{\alpha}(\xi)$ are arbitrary respectively real (purely imaginary)
trigonometric polynomials of semi-integer degrees associated with $\sigma$.
    \label{GenSymB}
\end{lem}

{\bf Proof.} Let $T(\xi)$ defined by~(\ref{EOpolB}).
As in Lemma's~\ref{GenSymA} proof, due to Leibniz formula, we have $D^{\beta} T (\nul) = 0$ for all
    $\beta\in\Delta_n$.
    Also, it is clear that $T(\xi)$ is a real trigonometric polynomial of semi-integer degrees associated with $\sigma$.

Now let $T(\xi)$ be a real trigonometric polynomial of semi-integer degrees associated with $\sigma$, such that
    $D^{\beta} T (\nul) = 0,  \forall\beta\in\Delta_n$.
    By Taylor formula, there exist trigonometric polynomials $T_{\alpha}(\xi)$
    of semi-integer degrees associated with $\sigma,$ where
    $\alpha\in\zd_+,$ $[\alpha]=n$, so  that

        $$T(\xi)=\sum\limits_{\alpha\in\zd_+ \atop [\alpha]=n} \Pi_{\alpha}(\xi) T_{\alpha}(\xi)=\hspace{9cm}$$
           $$\sum\limits_{\alpha\in\zd_+ \atop [\alpha]=n} \left(Re T_{\alpha}(\xi) Re \Pi_{\alpha}(\xi)- Im T_{\alpha}(\xi)Im \Pi_{\alpha}(\xi)
        +i Im T_{\alpha}(\xi)Re \Pi_{\alpha}(\xi) +i Re T_{\alpha}(\xi) Im \Pi_{\alpha}(\xi)\right).$$

    Since $T(\xi)=\overline{T(\xi)}$, then imaginary part is
    equal to zero. It remains to set $\w A_{\alpha}(\xi)=Re T_{\alpha}(\xi),$
    $\w B_{\alpha}(\xi)=- Im T_{\alpha}(\xi).$

    For purely imaginary $T(\xi)$, the proof is similar.
 $\Diamond$

\begin{theo}
    Let $M$ be an arbitrary matrix dilation, $c\in\frac 12\zd$  and $n\in\n.$
    A general form for all masks $m_0,$ which are
    symmetric with
    respect to the point $c$ in the sense (a) and have sum rule of order $n$,
    %satisfying condition~(\ref{19mod}) with $\lambda'_{\alpha}$
    is given by
    \ba
    m_0(\xi)=m_0^*(\xi) + e^{2\pi i (c,\xi)} \sum\limits_{j\in J} \sum\limits_{\alpha\in\zd_+ \atop [\alpha]=n}
    (A_{\alpha j}(M^*\xi)Re \Pi_{\alpha}(M^*\xi)+i B_{\alpha j}(M^*\xi) Im \Pi_{\alpha}(M^*\xi)) +
    \nonumber
    \\
     \sum\limits_{i\in I^1 }
    \sum\limits_{\alpha\in\zd_+ \atop [\alpha]=n} \left(
    T_{\alpha i}(M^*\xi)\Pi_{\alpha}(M^*\xi) e^{2\pi i (s_i,\xi)}+
    T_{\alpha i}(-M^*\xi)\Pi_{\alpha}(-M^*\xi) e^{2\pi i (2c-s_i,\xi)}\right),
    \label{GenFormA}
    \ea
    where for $\alpha\in\zd_+,$ $[\alpha]=n,$ $j\in J$, $i\in I^1,$ $A_{\alpha j}(\xi)$, $B_{\alpha j}(\xi)$ are arbitrary respectively even and odd
    trigonometric polynomials  of semi-integer degrees associated with $\sigma=M^{-1}(c-s_j)-\lfloor M^{-1}(c-s_j)\rfloor$,
    $T_{\alpha i}(\xi)$ are arbitrary trigonometric polynomials, $m_0^*(\xi)$ is constructed
    by Theorem~\ref{theoMaskA} %using formula~(\ref{m_0A})
    for some parameters
    $\lambda'_{\alpha}\in\cn$, $\alpha\in\Delta_n$, $\lambda'_{\alpha}=0$, if $[\alpha]$ is odd.
    \label{GenFormLemA}
\end{theo}

{\bf Proof.} Let $m_0$ be a mask
    that is symmetric with
    respect to the point $c$ in the sense (a) and
     has sum rule of order $n$. Define
    complex numbers $\lambda'_{\alpha}$ for all $\alpha\in\Delta_n$ by

        $$\lambda'_{\alpha}:=\frac 1 {(2\pi i)^{[\alpha]}}
        D^{\alpha}\left(e^{-2\pi i (c,M^{*-1}\xi)}m_0(M^{*-1}\xi)\right)\Bigg|_{\xi=\nul}.$$
    It is clear that $\lambda'_{\nul}=1,$ $\lambda'_{\alpha}=0$, if $[\alpha]$ is odd and condition~(\ref{19mod}) holds
    with $\lambda'_{\alpha},$ $\alpha\in\Delta_n.$
        Due to Theorem~\ref{theoMaskA} there exists mask $m_0^*$
        that is symmetric with
    respect to the point $c$ in the sense (a) and satisfies condition~(\ref{19})
    with numbers $\lambda_\alpha$ defined by~(\ref{lam'}) or equivalently satisfying
     condition~(\ref{19mod}) with $\lambda'_\alpha$.
    %that is constructed by~(\ref{m_0A}).

The polyphase components of $m_0-m_0^*$ are $\mu_{0k}-\mu_{0k}^*$, $k=0,\dots,m-1,$
    where $\mu_{0k}$ and $\mu_{0k}^*$ are the polyphase components of $m_0$ and $m_0^*$
    respectively.
    Then, due to the fact that condition~(\ref{19mod}) holds for $m_0$ and $m_0^*$ with
    the same $\lambda'_{\alpha}$,  we have $D^{\beta}(\mu_{0k}(\xi)-\mu_{0k}^*(\xi))\Big|_{\xi=\nul}=0,$
    $k=0,\dots,m-1,$ for all $\beta\in\Delta_n.$

For $j\in J$, $\mu_{0j}(\xi)-\mu_{0j}^*(\xi)$
    is symmetric with respect to the point $M^{-1}(c-s_j)$ in the sense (a),
    so together with Lemma~\ref{GenSymA} we obtain

    \ba
    \mu_{0j}(\xi)-\mu_{0j}^*(\xi)=
    e^{2\pi i \left(M^{-1}(c-s_j),\xi\right)}\sum\limits_{\alpha\in\zd_+ \atop [\alpha]=n}
    (A_{\alpha j}(\xi)Re \Pi_{\alpha}(\xi)+i B_{\alpha j}(\xi) Im \Pi_{\alpha}(\xi)),
    \label{Work2}
    \ea
    where for $\alpha\in\zd_+,$ $[\alpha]=n,$ $A_{\alpha j}(\xi)$ and $B_{\alpha j}(\xi)$ are respectively even and odd
    trigonometric polynomials of semi-integer degrees associated with $\sigma=M^{-1}(c-s_j)-\lfloor M^{-1}(c-s_j)\rfloor$.

    For $i\in I$, $\mu_{0i}(\xi)-\mu_{0i}^*(\xi)$ can be represented as

        \be
        \mu_{0i}(\xi)-\mu_{0i}^*(\xi)=  \sum\limits_{\alpha\in\zd_+ \atop [\alpha]=n}
        T_{\alpha i}(\xi)\Pi_{\alpha}(\xi),
        \label{Work1}
        \ee
    where for $\alpha\in\zd_+,$ $[\alpha]=n,$ $T_{\alpha i}(\xi)$ are trigonometric polynomials and due to~(\ref{ItypeA}) for $i'\in I,$
    such that $2c-s_i-s_{i'}=0,$ we have

        $$\mu_{0i'}(\xi)-\mu_{0i'}^*(\xi)=\mu_{0i}(-\xi)-\mu_{0i}^*(-\xi)=
        \sum\limits_{\alpha\in\zd_+ \atop [\alpha]=n}
        T_{\alpha i}(-\xi)\Pi_{\alpha}(-\xi).$$
    It remains to combine $\mu_{0k}(\xi)-\mu_{0k}^*(\xi),$ $k=0,\dots,m-1,$ by~(\ref{PR}) and take into account that $s_{i'}=2c-s_i$.
    %and set $T_{\alpha i}=2T'_{\alpha i}.$

Next, we assume that trigonometric polynomial $m_0$ defined by~(\ref{GenFormA}).
    Trigonometric polynomial $m_0^*$ is symmetric
    with respect to the point $c$ in the sense (a) and satisfy condition~(\ref{19}) with numbers $\lambda_\alpha\in\cn$.
    It is clear that the polyphase representatives of $m_0-m_0^*$ are
given by~(\ref{Work1}) and~(\ref{Work2}). Due to Lemma~\ref{GenSymA},
$D^{\beta}(\mu_{0k}(\xi)-\mu_{0k}^*(\xi))\Big|_{\xi=\nul}=0,$
    $k=0,\dots,m-1,$ for all $\beta\in\Delta_n.$
    Since

        $$D^\beta\mu_{0k}^*(\nul)=\frac {(2\pi i)^{[\beta]}} {\sqrt m}
        \sum\limits_{\nul\le\alpha\le\beta}\lambda_\alpha
        \lll\beta\atop\alpha\rrr(- M^{-1}s_k)^{\beta-\alpha} \quad  \forall
        \beta\in\Delta_n,$$
       then (\ref{19}) is valid for the polyphase components of $m_0$ with numbers $\lambda_\alpha.$ Therefore, $m_0$ has sum rule of order $n.$
   It is easy to see that $m_0$ is symmetric with respect to the point $c$ in the sense (a).
$\Diamond$
%
%It is clear that expression
%        $T_{\alpha i}(M^*\xi)\Pi_{\alpha}(M^*\xi) e^{2\pi i (s_i,\xi)}+
%        T_{\alpha i}(-M^*\xi)\Pi_{\alpha}(-M^*\xi) e^{2\pi i (2c-s_i,\xi)}$
%    is symmetric with respect to the point $c$ in the sense (a).

The same theorem is valid for the case of symmetric masks with
    respect to a point in the sense (b).

\begin{theo}
     Let $M$ be an arbitrary matrix dilation, $c\in\frac 12\zd$  and $n\in\n.$
    A general form for all masks $m_0,$ which are symmetric with
    respect to the point $c$ in the sense (b) and have sum rule of order $n$, is given by
    $$
    m_0(\xi)=m_0^*(\xi) + e^{2\pi i (c,\xi)} \sum\limits_{j\in J}
    \sum\limits_{\alpha\in\zd_+ \atop [\alpha]=n} \left(
    \w A_{\alpha j}(M^*\xi)Re \Pi_{\alpha}(M^*\xi)+\w B_{\alpha j}(M^*\xi) Im \Pi_{\alpha}(M^*\xi) \right)+
    $$
    %\nonumber
    %\\
    $$
     \sum\limits_{i\in I^1 }
    \sum\limits_{\alpha\in\zd_+ \atop [\alpha]=n} \left(
    T_{\alpha i}(M^*\xi)\Pi_{\alpha}(M^*\xi) e^{2\pi i (s_i,\xi)}+
    \overline{T_{\alpha i}(M^*\xi)\Pi_{\alpha}(M^*\xi)} e^{2\pi i (2c-s_i,\xi)}\right),
    %\label{GenFormB}
    $$
    where for $\alpha\in\zd_+,$ $[\alpha]=n,$ $j\in J$, $i\in I^1,$ $\w A_{\alpha j}(\xi)$, $\w B_{\alpha j}(\xi)$ are arbitrary real
trigonometric polynomials of semi-integer degrees associated with $\sigma=M^{-1}(c-s_j)-\lfloor M^{-1}(c-s_j)\rfloor$,
  $T_{\alpha i}(\xi)$ are arbitrary trigonometric polynomials, $m_0^*$ is constructed
by Theorem~\ref{theoMaskB} using formula~(\ref{m_0B}) for some numbers $\lambda'_{\alpha}\in\cn$, $\alpha\in\Delta_n$, $Re \lambda'_{\alpha}=0,$ if $[\alpha]$ is odd,
    $Im\lambda'_{\alpha}=0,$ if $[\alpha]$ is even.
    \label{GenFormLemB}
\end{theo}

\subsection{Construction of symmetric frame-like wavelets}

We know how to get the initial symmetric refinable mask $m_0$ that has arbitrary order of sum rule.
    On the way for constructing symmetric/antisymmetric frame-like wavelet system
    such that wavelet functions $\w\psi^{(\nu)},$ $\nu=1,\dots,m,$ have $VM^n$ property we need to find appropriate dual mask $\w m_0$.
    The method is given by the following theorem.

\begin{theo}
Let $M$ be an arbitrary matrix dilation, $c\in\frac 12\zd,$ $n\in\n$, mask $m_0$ be symmetric
    with respect to the point $c$ in the sense (a) or in the sense (b)
    and has sum rule of order $n$.
    %satisfying condition~(\ref{19}).
    Then there exists mask $\w m_0$ which is symmetric
    with respect to the point $c$ in the sense (a) or in the sense (b) and satisfies condition~(\ref{20_new}).
    \label{DualTheo}
\end{theo}

{\bf Proof.} Firstly, we transform condition~(\ref{20_new}) into
        $$D^\beta\lll 1-e^{-2\pi i (c, \xi)} m_0(\xi)\overline{e^{-2\pi i (c, \xi)}\w m_0(\xi)}\rrr\Big|_{\xi=0}=0,
        \quad \forall \beta\in\Delta_n.$$
    From this equation we consequentially find
    $\rho_{\beta}=\frac 1 {(2\pi i)^{[\beta]}}
        D^\beta \lll e^{-2\pi i ( c , \xi)}\w m_0(\xi)\rrr\Big|_{\xi=0}.$
        For the case of symmetry in the sense (a) $\rho_{\beta}=0$, if $[\beta]$ is odd.
    Then we set %for this case of symmetry
        \be
        e^{-2\pi i (c, \xi)}\w m_0(\xi)=
        \sum\limits_{\beta\in\Delta^e_n} \rho_{\beta} G^{A}_{\beta}(\xi),
        \label{fDualMaskA}
        \ee
    where trigonometric polynomials of semi-integer degrees associated with
    $\sigma=c-\lfloor c \rfloor$
     $G^{A}_{\beta}(\xi)\in\Theta^A_{\beta, n}$ for $\beta\in\Delta^e_n.$
    For the case of symmetry in the sense (b) we set

        $$e^{-2\pi i (c, \xi)}\w m_0(\xi)=
        \sum\limits_{\beta\in\Delta_n} \rho_{\beta} G^{B}_{\beta}(\xi),$$
        where trigonometric polynomials of semi-integer degrees associated with
    $\sigma=c-\lfloor c \rfloor$
    $G^{B}_{\beta}(\xi)\in\Theta^B_{\beta, n}$ for $\beta\in\Delta_n.$

Thus, $\w m_0(\xi) $ is symmetric with respect to the point $c$ in the sense (a) or in the sense (b)
    and satisfy condition~(\ref{20_new}).
$\Diamond$
%
%\begin{rem}
%Full and axial symmetry??
%\label{rSym}
%\end{rem}

\begin{theo} Let $M$ be an arbitrary matrix dilation, $c\in\frac 12\zd,$ $n\in\n.$ Masks $m_0$ and $\w m_0$
    are symmetric with
    respect to the point $c$ in the sense (a), mask $m_0$
    has sum rule of order $n$,
    %satisfy condition~(\ref{19}),
    mask $\w m_0$ satisfies
    condition~(\ref{20_new}).
    Then there exist wavelet masks  $m_{\nu}$ and $\w m_{\nu}$, $\nu=1,\dots,m,$
    which are symmetric/antisymmetric  with
    respect to some points in the sense (a) and
    %provide approximation order $n$ for the corresponding frame-like wavelet system.
    wavelet masks $\w m_{\nu},$ $\nu=1,\dots,m,$ have vanishing moments of order $n$.
\label{theoWaveA}
\end{theo}

{\bf Proof.} Matrix extension providing vanishing moments of order $n$ for $\w m_{\nu},$ $\nu=1,\dots,m,$
    realises as follows (see~\cite[Lemma 14]{KrSk})
    \ba \cal N:=\left(\begin{matrix}
        \mu_{00}&\mu_{01}&\hdots&\mu_{0,m-1}&\mu_{0,m}\cr
        1 & 0 & \hdots &0& -\overline{\widetilde\mu}_{00} \cr
        0 & 1 & \hdots &0& -\overline{\widetilde\mu}_{01} \cr
        \vdots & \vdots & \ddots &\vdots& \vdots \cr
        0 & 0 & \hdots &0& -\overline{\widetilde\mu}_{0, m-2} \cr
        0 & 0 & \hdots &1& -\overline{\widetilde\mu}_{0, m-1} \cr
        \end{matrix}\right),\,
        \hspace{1.5cm}
        \label{NNMatr0}
    \\
        \cal \widetilde N:=\left(\begin{matrix}
        \widetilde\mu_{00}&\widetilde\mu_{01}&\hdots&\widetilde\mu_{0m-1}&1\cr
        1-\w\mu_{00}\overline\mu_{00} & -\w\mu_{01}\overline\mu_{00} & \hdots &-\w\mu_{0,m-1}\overline\mu_{00}&-\overline\mu_{00} \cr
        -\w\mu_{00}\overline\mu_{01} & 1-\w\mu_{01}\overline\mu_{01} & \hdots &-\w\mu_{0,m-1}\overline\mu_{01}&-\overline\mu_{01} \cr
        \vdots & \vdots & \ddots &\vdots& \vdots \cr
        -\w\mu_{00}\overline\mu_{0,m-2}&-\w\mu_{01}\overline\mu_{0,m-2}&\hdots&-\w\mu_{0,m-1}\overline\mu_{0,m-2}&-\overline\mu_{0,m-2}\cr
        -\w\mu_{00}\overline\mu_{0,m-1}&-\w\mu_{01}\overline\mu_{0,m-1}&\hdots&1-\w\mu_{0,m-1}\overline\mu_{0,m-1}&-\overline\mu_{0,m-1}\cr
        \end{matrix}\right),
        \label{NNMatr}
        \ea
    where $\mu_{0k}(\xi),$ $\w\mu_{0k}(\xi),$ $k=0,\dots,m-1$ are the polyphase components of masks
    $m_0$ and $\w m_0$ respectively, $\mu_{0,m}=1-\sum_{k=0}^{m-1}\mu_{0k}\overline{\w\mu_{0k}}.$
    It is easy to check that ${\cal N}{\cal \widetilde N}^*=I_{m+1}$ holds.
    Therefore, the equality~(\ref{calM}) is valid.

 Let us denote matrices elements by

        $${\cal N}=\{\mu_{kl}\}_{k,l=0}^m,\quad
        {\cal \w N}=\{\w\mu_{kl}\}_{k,l=0}^m,$$
        the first rows of the matrices by $P,$ $\w P$ and the other rows of the matrices by  $Q_{\nu},$ $\w Q_{\nu},$ respectively, $\nu=0,\dots,m-1.$

Due to the symmetry of refinable masks $m_0$ and $\w m_0$ by Lemma~\ref{PolyLem} we have

    $$\mu_{0i}(\xi)=\mu_{0i'}(-\xi),\quad  \w \mu_{0i}(\xi)=\w\mu_{0i'}(-\xi), \quad
    2c-s_i-s_{i'}=0,\,\,i, i'\in I,$$

    $$\mu_{0j}(\xi)=e^{2\pi i (M^{-1}(2c-2s_j),\xi)}\mu_{0j}(-\xi),\quad \w\mu_{0j}(\xi)=e^{2\pi i (M^{-1}(2c-2s_j),\xi)}\w\mu_{0j}(-\xi),\quad   j\in J.$$

To provide symmetry for wavelet masks
    we make some transformations and obtain symmetric/an\-ti\-symmetric wavelet masks from the above rows.

Let $k^*\in J$.  The row $\w Q_{k^*}$ looks as
        $$\w Q_{k^*}=(\dots,-\widetilde\mu_{0i}\overline\mu_{0k^*},\hdots,1-\widetilde\mu_{0k^*}\overline\mu_{0k^*},
        \dots,-\widetilde\mu_{0j}\overline\mu_{0k^*},\dots),$$
    where $i\in I,$ $j\in J,$ $j\neq k^*.$ Now, we show that first $m$
    elements of the row $\w Q_{k^*}$ (that are the polyphase components of wavelet mask $\w m_{k^*}$) satisfy condition~(\ref{PolyTypeA})

        $$\w\mu_{k^* i}(\xi)=-\widetilde\mu_{0i}(\xi)\overline{\mu_{0k^*}(\xi)}=
        -e^{-2\pi i (M^{-1}(2c-2s_{k^*}),\xi)}\widetilde\mu_{0i'}(-\xi)\overline{\mu_{0k^*}(-\xi)}
        =$$
        $$=e^{2\pi i (M^{-1}(2s_{k^*}-s_i-s_{i'}),\xi)}\w\mu_{k^* i'}(-\xi),\quad \forall i, i'\in I,\quad  2c-s_i-s_{i'}=0;$$

        $$\w\mu_{k^* j}(\xi)=-\widetilde\mu_{0j}(\xi)\overline{\mu_{0k^*}(\xi)}=
        -e^{2\pi i (M^{-1}(2c-2s_j),\xi)}e^{-2\pi i (M^{-1}(2c-2s_{k^*}),\xi)}
            \widetilde\mu_{0j}(-\xi)\overline{\mu_{0k^*}(-\xi)}
        =$$
        $$=e^{2\pi i (M^{-1}(2s_{k^*}-2s_j),\xi)}\w\mu_{k^* j}(-\xi),\quad \forall  j \in J,\quad j\neq k^*;$$

        $$\w\mu_{k^* k^*}(\xi)=1-\widetilde\mu_{0k^*}(\xi)\overline{\mu_{0k^*}(\xi)}=
        1-\widetilde\mu_{0k^*}(-\xi)\overline{\mu_{0k^*}(-\xi)}
        =\w\mu_{k^* k^*}(-\xi).$$
        %$$e^{2\pi i (M^{-1}(2s_{k^*}-2s_{k^*}),\xi)}\w\mu_{k^* k^*}(-\xi)=$$\w\mu_{k^* k^*}(-\xi).$$
Thus, by Lemma~\ref{PolyLem} wavelet mask $\w m_{k^*}$  is symmetric with respect
    to the point $s_{k^*}$ in the sense (a) for all $k^*\in J.$

If $k^*\in I^1$ we find $k'\in I^2$, such that
    $2c-s_{k^*}-s_{k'}=0.$
    %$\mu_{0k'}$, $\mu_{0k^*}$ we have
    %are flipped version of each other,
    Therefore, $\mu_{0k^*}(\xi)=\mu_{0k'}(-\xi)$.
    Let us consider the rows

        $$\w Q_{k^*}=(\dots,-\widetilde\mu_{0i}\overline\mu_{0k^*},\dots,1-\widetilde\mu_{0k^*}\overline\mu_{0k^*},\dots,
        -\widetilde\mu_{0k'}\overline\mu_{0k^*},
        \dots,-\widetilde\mu_{0j}\overline\mu_{0k^*},\dots),$$
        $$\w Q_{k'}=(\dots,-\widetilde\mu_{0i}\overline\mu_{0k'},\dots,-\widetilde\mu_{0k^*}\overline\mu_{0k'},\dots,
        1-\widetilde\mu_{0k'}\overline\mu_{0k'},
        \dots,-\widetilde\mu_{0j}\overline\mu_{0k'},\dots),$$
    where $i\in I,$ $j\in J,$ $i\neq k^*,$  $i\neq k',$
    and transform them into

        %\ba
        $$\w Q'_{k^*}:=\frac 1 2 (\w Q_{k^*}+\w Q_{k'}),\quad
        \w Q'_{k'}:=\frac 1 2 (\w Q_{k^*}-\w Q_{k'}).$$
        %\label{BioNew}
        %\ea
    Define
        $A(\xi)=\frac 12\left(\overline{\mu_{0k^*}(\xi)}+\overline{\mu_{0k'}(\xi)}\right),$
        $B(\xi)=\frac 12\left(\overline{\mu_{0k^*}(\xi)}-\overline{\mu_{0k'}(\xi)}\right).$
        It is clear that
        $A(\xi)$ is even, i.e. $A(-\xi)=A(\xi);$
     and $B(\xi)$ is odd, i.e. $B(-\xi)=-B(\xi).$
    New rows take the form

        $$\w Q'_{k^*}=
        (\dots,-\widetilde\mu_{0i}A,\dots,\frac 12-\widetilde\mu_{0k^*}A,\dots,\frac 12-\widetilde\mu_{0k'}A,\dots,-\widetilde\mu_{0j}A,\dots),$$
         $$\w Q'_{k'}=
         (\dots,-\widetilde\mu_{0i}B,\dots,\frac 12-\widetilde\mu_{0k^*}B,\dots,-\frac 12 -\widetilde\mu_{0k'}B,\dots,-\widetilde\mu_{0j}B,\dots).$$
         Redenote elements such that
            $$\w Q'_{k^*}=\{\w\mu_{k^*l}\}_{l=0,\dots,m},\quad
        \w Q'_{k'}=\{\w\mu_{k'l}\}_{l=0,\dots,m}.$$
    Now, we show that first $m$
    elements of the row $\w Q'_{k^*}$ (that are the polyphase components of wavelet mask $\w m_{k^*}$) satisfy condition~(\ref{PolyTypeA})
        $$\w\mu_{k^* i}(\xi)=-\widetilde\mu_{0i}(\xi)A(\xi)=
        -\widetilde\mu_{0i'}(-\xi)\,A(-\xi)=
        \w\mu_{k^* i'}(-\xi), \quad \forall i, i'\in I,\, i\neq k^*, i\neq k';$$

        $$\w\mu_{k^* j}(\xi)=-\widetilde\mu_{0j}(\xi)A(\xi)=
        -e^{2\pi i (M^{-1}(2c-2s_j),\xi)}\widetilde\mu_{0j}(-\xi)\, A(-\xi)=$$
        $$e^{2\pi i (M^{-1}(2c-2s_j),\xi)}\w\mu_{k^* j}(-\xi), \quad  \forall j \in J;$$

        $$\w\mu_{k^* k^*}(\xi)=\frac 12-\widetilde\mu_{0k^*}(\xi)A(\xi)=
        \frac 12-\widetilde\mu_{0k'}(-\xi)\,A(-\xi)=
        \w\mu_{k^* k'}(-\xi).$$
    Hence, by Lemma~\ref{PolyLem} wavelet mask $\w m_{k^*}$ is symmetric with respect to the point $c$ in the sense (a).
    The same equations are valid for the elements of the row $\w Q'_{k'}$

        $$\w\mu_{k' i}(\xi)=-\w\mu_{k' i'}(-\xi), \quad \forall i, i'\in I,\, i\neq k^*, i\neq k';$$

        $$\w\mu_{k' j}(\xi)=-e^{2\pi i (M^{-1}(2c-2s_j),\xi)}\w\mu_{k' j}(-\xi), \quad  \forall j \in J;$$

        $$\w\mu_{k' k^*}(\xi)=-\w\mu_{k' k'}(-\xi).$$
    Therefore, wavelet mask $\w m_{k'}$ is antisymmetric with respect to the point $c$ in the sense (a).

For the rows $Q_{k^*},$ $Q_{k'}$ we make analogous transformations
        $$Q'_{k^*}:=Q_{k^*}+Q_{k'},\quad
        Q'_{k'}:=Q_{k^*}-Q_{k'}.$$
    New rows
        $$Q'_{k^*}=(\dots,1,\dots,1,\dots),$$
        $$Q'_{k'}=(\dots,1,\dots,-1,\dots)$$
    form wavelet masks $m_{k^*}$ and $m_{k'}$ that are accordingly symmetric and
    antisymmetric with respect to the point $c$.

It is clear that the equality ${\cal N}\w{\cal N}^*=I_{m+1}$ remains correct with new rows.

$\Diamond$

\begin{rem} It is important to note that dual wavelet frames cannot be constructed using such extension technique due to
    the fact that in matrix~(\ref{NNMatr}) we take $\w\mu_{0m} \equiv 1.$
    Therefore, wavelet masks $m_{\nu},$ $\nu=1,\dots,m$ do not have any order of vanishing moments
    according with~\cite[Theorem 8]{Sk1}.
\label{rNotFrame}
\end{rem}

The same theorem is valid for the case of the symmetry in the sense (b)

\begin{theo} Let $M$ be an arbitrary matrix dilation, $c\in\frac 12\zd,$ $n\in\n.$ Masks $m_0$ and $\w m_0$
    are symmetric with
    respect to the point $c$ in the sense (b), mask $m_0$
    has sum rule of order $n$,
    %satisfy condition~(\ref{19}),
    mask $\w m_0$ satisfies
    condition~(\ref{20_new}).
    Then there exist wavelet masks  $m_{\nu}$ and $\w m_{\nu}$, $\nu=1,\dots,m,$
    which are symmetric/antisymmetric with
    respect to some points in the sense (b) and
    %provide approximation order $n$ for the corresponding frame-like wavelet system.
    masks $\w m_{\nu},$ $\nu=1,\dots,m,$  have vanishing moments of order $n$.
\label{theoWave}
\end{theo}

\section{Bivariate axial symmetric frame-like wavelets.}

Let us denote $G^{axis}:=\left\{\pm I_2, \pm Y\right\},$
    where $Y=\left(\begin{matrix}
            -1 & 0\cr
            0 & 1\cr
        \end{matrix}\right)$. The set $G^{axis}$ is called axial symmetry group on $\z^2$.
    Recall that trigonometric polynomial $t(\xi)$ is axial symmetric
    with respect to a center $c=(c_1,c_2)\in \frac 12 \z^2$, if

       $$t(\xi)=e^{2\pi i (c-Ec,\xi)}t(E^*\xi),
       \quad\forall E\in G^{axis},\quad \xi\in\r^2$$
    or
       $$ h_n=h_{E(n-c)+c}, \quad \forall n\in\z^2,\quad \forall E\in G^{axis}.$$
       Similarly, trigonometric polynomial $t(\xi)$ is axial antisymmetric
       with respect to the center $c$, if

       $$t(\xi)=\pm e^{2\pi i (c-Ec,\xi)}t(E^*\xi),
       \quad\forall E\in G^{axis}.$$
    To provide axial symmetry for $t(\xi)$ it is enough to satisfy %~(\ref{TrigSym})
        only two equalities
              \be
              t(\xi)=e^{2\pi i ( 2c, \xi)}{t(-\xi)},
              \quad t(\xi)=e^{2\pi i 2 c_1 \xi_1 } t(Y^* \xi),\quad (c_1,c_2)\in \frac 12 \z^2.
              \label{AxisMask}
              \ee

For which dilation matrices $G^{axis}$-symmetry of mask $m_0$ carries over to its refinable function?
     The answer is given in

\begin{lem} (see~\cite[Proposition 4.1]{Han3})
     Let symmetry group $G^{axis}$ be a symmetry group with respect to the dilation matrix $M$. Then
     \be M=\left(\begin{matrix}
            \gamma_1 & 0\cr
            0 & \gamma_2\cr
      \end{matrix}\right) \quad  \texttt{or} \quad M=\left(\begin{matrix}
            0 & \gamma_1\cr
            \gamma_2 & 0\cr
      \end{matrix}\right), \quad \gamma_1, \gamma_2\in\z.
      \label{AxisM}\ee
\label{AxisLem}
\end{lem}

In this section, we assume that dilation matrices are as in~(\ref{AxisM}), $m=|\gamma_1\gamma_2|.$
     Also, we set $d_1=|\gamma_1|$, $d_2=|\gamma_2|$.

Firstly, we consider the case when $c\in\z^2.$
    Due to the fact that orders of sum rule, vanishing moments, linear-phase moments and axial symmetry of trigonometric polynomial
    are invariant with respect to an integer shift it is enough to consider only the case when $c=(0,0).$
    For dilation matrix $M$
    the digits always can be choose in the rectangular set

     $$D(M)\subset\left[\frac {-d_1+1} 2,\frac {d_1} 2\right]\times
     \left[\frac {-d_2+1} 2,\frac {d_2} 2\right].$$

     Let us define the set of digits $D(M)$ as

     $$D(M)=\left\{\left(v_1, v_2\right) \Big| \quad
        v_l=-\left\lfloor\frac {d_l-1} 2\right\rfloor,\dots, \left\lfloor\frac {d_l} 2\right\rfloor, l=1,2\right\}.$$
        %v_2=\left\lfloor\frac {-d_2+1} 2\right\rfloor,\dots, \left\lfloor\frac {d_2} 2\right\rfloor
     Enumerate all digits from $0$ to $m-1$ and divide them into several sets with respect to the parity of integers $d_1$ and $d_2$.

     I. If $d_1, d_2$ are odd, then
       digits can be divided as follows
            \begin{enumerate}
              \item $J_I:=\{j_1  \Big| s_{j_1}=(0, 0)\};$
              \item $R:=\{r \Big| s_r=(v_1, 0),  \quad v_1=1,\dots, \left\lfloor\frac {d_1-1} 2\right\rfloor\},$

                    $R':=\{r' \Big| s_{r'}=(-v_1, 0), \quad v_1=1,\dots, \left\lfloor\frac {d_1-1} 2\right\rfloor\},$

                        $Q:=\{q  \Big| s_q=(0, v_2), \quad v_2=1,\dots, \left\lfloor\frac {d_2-1} 2\right\rfloor\},$

                        $Q':=\{q'  \Big|s_{q'}=(0, -v_2), \quad v_2=1,\dots, \left\lfloor\frac {d_2-1} 2\right\rfloor\};$
              \item %оставшиеся могут быть сгруппированы в четверки вида
              $
              K:=\{ k \Big| s_k=(v_1, v_2),\quad v_1=1,\dots, \left\lfloor\frac {d_1-1} 2\right\rfloor, v_2=1,\dots, \left\lfloor\frac {d_2-1} 2\right\rfloor\},$

              $K':=\{ k' \Big| s_{k'}=(v_1, -v_2),\quad v_1=1,\dots, \left\lfloor\frac {d_1-1} 2\right\rfloor, v_2=1,\dots, \left\lfloor\frac {d_2-1} 2\right\rfloor\},$

              $K'':=\{ k'' \Big| s_{k''}=(-v_1, v_2),\quad v_1=1,\dots, \left\lfloor\frac {d_1-1} 2\right\rfloor, v_2=1,\dots, \left\lfloor\frac {d_2-1} 2\right\rfloor\},$

              $K''':=\{ k''' \Big| s_{k'''}=(-v_1, -v_2),\quad v_1=1,\dots, \left\lfloor\frac {d_1-1} 2\right\rfloor, v_2=1,\dots, \left\lfloor\frac {d_2-1} 2\right\rfloor\};
              $
            \end{enumerate}

        II. If $d_1$ is even, $d_2$ is odd, then
        digits can be divided as follows
            \begin{enumerate}
              \item $J_{II}:=\{ j_1 , j_2 \Big| s_{j_1}=(0, 0), s_{j_2}=(\frac {d_1} 2, 0)\};$
              \item $R, R'$ and $Q, Q'$ are the same as in I;
              \item $T:=\{t \Big| s_t=(\frac {d_1} 2, v_2), v_2=1,\dots,
              \left\lfloor\frac {d_2-1} 2\right\rfloor\};$

              $T':=\{t' \Big| s_{t'}=(\frac {d_1} 2, -v_2), v_2=1,\dots,
              \left\lfloor\frac {d_2-1} 2\right\rfloor\};$
              \item $K, K', K'', K'''$ are the same as in I.
            \end{enumerate}

        III. If $d_1$ is odd, $d_2$ is even, then
        digits can be divided as follows
            \begin{enumerate}
              \item $J_{III}:=\{ j_1 , j_3 \Big| s_{j_1}=(0, 0), s_{j_3}=(0, \frac {d_2} 2)\};$
              \item $R, R'$ and $Q, Q'$ are the same as in I;
              \item $U:=\{u \Big| s_u=(v_1, \frac {d_2} 2),
              v_1=1,\dots, \left\lfloor\frac {d_1-1} 2\right\rfloor\};$

              $U':=\{u' \Big| s_{u'}=(-v_1, \frac {d_2} 2),
              v_1=1,\dots, \left\lfloor\frac {d_1-1} 2\right\rfloor\};$
              \item $K, K', K'', K'''$ are the same as in I.
            \end{enumerate}

        IV. If $d_1, d_2$ are even, then
        digits can be divided as follows
            \begin{enumerate}
              \item $J_{IV}:=\{j_1 , j_2, j_3 , j_4 \Big| s_{j_1}=(0, 0), s_{j_2}=(\frac {d_1} 2,0),
              s_{j_3}=(0, \frac {d_2} 2), s_{j_4}=(\frac {d_1} 2, \frac {d_2} 2)\};$
              \item $R, R'$ and $Q, Q'$ are the same as in I;
              \item $T, T'$ and $U, U'$ are the same as in II and III;
              \item $K, K', K'', K'''$ are the same as in I.
            \end{enumerate}
     Let us reformulate axial symmetry conditions~(\ref{AxisMask}) for mask $m_0$ in terms of its polyphase components. To cover all
     cases we give them for the IV type of dilation matrices.

For the I type of dilation matrices mask $m_0$ can be represented as

            %\ba
            $$m_0(\xi)\sqrt m=\mu_{0j_1}(M^*\xi) e^{2\pi i  (s_{j_1},\xi)}+
            \sum\limits_{r\in R}\mu_{0r}(M^*\xi)e^{2\pi i (s_r,\xi)}+
            \sum\limits_{r'\in R'}\mu_{0r'}(M^*\xi)e^{2\pi i (s_{r'},\xi)}+$$
          %   \nonumber
           % \\
            $$\hspace{0.7cm}\sum\limits_{q\in Q}\mu_{0q}(M^*\xi)e^{2\pi i (s_q,\xi)}+\
            \sum\limits_{q'\in Q'}\mu_{0q'}(M^*\xi)e^{2\pi i (s_{q'},\xi)}+%\hspace{0.7cm}
            $$
          %   \nonumber
          %  \\
            $$\sum\limits_{k\in K}\mu_{0k}(M^*\xi)e^{2\pi i (s_k,\xi)}+
            \sum\limits_{k'\in K'}\mu_{0k'}(M^*\xi)e^{2\pi i (s_{k'},\xi)}+\hspace{0.7cm}$$
           %  \nonumber
           % \\
             $$\sum\limits_{k''\in K''}\mu_{0k''}(M^*\xi)e^{2\pi i (s_{k''},\xi)}+
             \sum\limits_{k'''\in K'''}\mu_{0k'''}(M^*\xi)e^{2\pi i (s_{k'''},\xi)}:=m_0^I(\xi),$$
            % \nonumber
            %\label{AxisMaskPolyI}
            %\ea
      %where for any summand $\mu_{0l}e^{2\pi i \left(a\atop b \right)\xi}$ the set $\Omega_M^l$ include digit $s_l=(a, b).$

For the II type of dilation matrices mask $m_0$ can be represented as
            $m_0(\xi)\sqrt m=m_0^I(\xi)+m_0^{II}(\xi),$
        where

        % \ba
           $$ m_0^{II}(\xi)=\sum\limits_{t\in T}
            \mu_{0t}(M^*\xi)e^{2\pi i (s_t,\xi)}+
            \sum\limits_{t'\in T'}
            \mu_{0t'}(M^*\xi)e^{2\pi i (s_{t'}\xi)}+
            \mu_{0j_2}(M^*\xi)e^{2\pi i (s_{j_2},\xi)}.$$
          %  \nonumber
            %\label{AxisMaskPolyII}
           % \ea

For the III type of dilation matrices mask $m_0$ can be represented as $m_0(\xi)\sqrt m=m_0^I(\xi)+m_0^{III}(\xi),$ where

         % \ba
          $$  m_0^{III}(\xi)=\sum\limits_{u\in U}
            \mu_{0u}(M^*\xi)e^{2\pi i (s_u,\xi)}+
            \sum\limits_{u'\in U'}
            \mu_{0u'}(M^*\xi)e^{2\pi i (s_{u'},\xi)}+
            \mu_{0j_3}(M^*\xi)e^{2\pi i  (s_{j_3},\xi)}.$$
          %  \nonumber
            %\label{AxisMaskPolyIII}
           % \ea

For the IV type of dilation matrices mask $m_0$ can be represented as

        $$ m_0(\xi)\sqrt m=m_0^I(\xi)+m_0^{II}(\xi)+m_0^{III}(\xi)+m_0^{IV}(\xi),$$
         where

          %\ba
           $$m_0^{IV}(\xi):=
          \mu_{0j_4}(M^*\xi)e^{2\pi i (s_{j_4},\xi)}.$$
            %\nonumber
            %\label{AxisMaskPolyIV}
            %\ea

From conditions~(\ref{AxisMask}) we can conclude the following symmetry conditions on
    the polyphase components according with the chose digits
       and taking into account that $c=(0,0)$
            \be
            s_{j_1}=(0, 0):\quad \mu_{0j_1}(\xi)=\mu_{0j_1}(-\xi), \quad\mu_{0j_1}(M^*\xi)=\mu_{0j_1}(M^*Y^*\xi),\hspace{3.6cm}
            \label{MaskPolyJ}
            \ee
            %\nonumber
            $$s_{j_2}=\left(\frac {d_1} 2,0\right):
            \quad \mu_{0j_2}(M^*\xi)=\mu_{0j_2}(-M^*\xi)e^{-2\pi i d_1\xi_1},\quad \mu_{0j_2}(M^*\xi)=\mu_{0j_2}(M^*Y^*\xi)e^{-2\pi i d_1\xi_1},$$
            %\nonumber
            %\\
            $$s_{j_3}=\left(0, \frac {d_2} 2\right):
            \quad \mu_{0j_3}(M^*\xi)=\mu_{0j_3}(-M^*\xi)e^{-2\pi i  d_2\xi_2},\quad \mu_{0j_3}(M^*\xi)=\mu_{0j_3}(M^*Y^*\xi),\hspace{1.3cm}$$
            %\nonumber
            %\\

              $$s_{j_4}=\left(\frac {d_1} 2, \frac {d_2} 2\right):
              \quad \mu_{0j_4}(\xi)=\mu_{0j_4}(-\xi)e^{-2\pi i (\xi_1+\xi_2)},\quad\mu_{0j_4}(M^*\xi)=\mu_{0j_4}(M^*Y^*\xi)e^{-2\pi i d_1\xi_1};\hspace{0.3cm}$$

            $$s_k=(v_1, v_2),\,  s_{k'}=(v_1, -v_2),\,  s_{k''}=(-v_1, v_2),\,
            s_{k'''}=(-v_1, -v_2):\hspace{4cm}$$
            $$\mu_{0k}(\xi)=\mu_{0k'''}(-\xi), \quad \mu_{0k'}(\xi)=\mu_{0k''}(-\xi),$$
            \be
            \mu_{0k}(M^*\xi)=\mu_{0k''}(M^*Y^*\xi), \quad\mu_{0k'}(M^*\xi)=\mu_{0k'''}(M^*Y^*\xi),
            \label{MaskPolyK}
            \ee

           % $$ $$
                   \be
            s_r=(v_1, 0),\,  s_{r'}=(-v_1, 0):
            \quad \mu_{0r}(\xi)=\mu_{0r'}(-\xi),\quad \mu_{0r}(M^*\xi)=\mu_{0r'}(M^*Y^*\xi),\hspace{1.5cm}
            \label{MaskPolyR}
            \ee

            \ba
            s_q=(0, v_2),\,s_{q'}=(0, -v_2):
            \quad \mu_{0q}(\xi)=\mu_{0q'}(-\xi), \quad \mu_{0q}(M^*\xi)=\mu_{0q}(M^*Y^*\xi),\hspace{1.5cm}
            \nonumber
            \\
            \quad \mu_{0q'}(M^*\xi)=\mu_{0q'}(M^*Y^*\xi),
            \label{MaskPolyQ}
            \ea

            $$ s_t=\left(\frac {d_1} 2, v_2\right),\,s_{t'}=\left(\frac {d_1} 2, -v_2\right):
            \quad \mu_{0t}(M^*\xi)=\mu_{0t'}(-M^*\xi)e^{-2\pi i d_1\xi_1},\hspace{3.2cm}$$
            %\nonumber
            %\\
            \be
            \mu_{0t}(M^*\xi)=\mu_{0t}(M^*Y^*\xi)e^{-2\pi i d_1\xi_1},\quad\mu_{0t'}(M^*\xi)=\mu_{0t'}(M^*Y^*\xi)e^{-2\pi i d_1\xi_1},
            \label{MaskPolyT}
            \ee

            \ba
             s_u=\left(v_1, \frac {d_2} 2\right),\, s_{u'}=\left(-v_1, \frac {d_2} 2\right):
             \quad\mu_{0u}(M^*\xi)=\mu_{0u'}(-M^*\xi)e^{-2\pi i d_2\xi_2}, \hspace{2.3cm}
            \nonumber
            \\
              \quad\mu_{0u}(M^*\xi)=\mu_{0u'}(M^*Y^*\xi),
            \label{MaskPolyU}
            \ea
    where $  v_1=1,\dots, \left\lfloor\frac {d_1-1} 2\right\rfloor,
    v_2=1,\dots, \left\lfloor\frac {d_2-1} 2\right\rfloor.$

 \begin{lem}
    Let $M$ be the matrix dilation.
    Then mask $m_0$ is  axial symmetric with respect to the origin if and only if
    its polyphase components $\mu_{0l}(\xi), $ $l=0,\dots,m-1,$ satisfy conditions~(\ref{MaskPolyJ})-(\ref{MaskPolyU})
    according to the type of dilation matrix.
    \label{lemAxisMask}
\end{lem}

{\bf Proof.} Let mask $m_0$ be axial symmetric with respect to the origin. Then with the above
    considerations conditions~(\ref{MaskPolyJ})-(\ref{MaskPolyU}) hold. Conversely,
    combine mask $m_0$ by~(\ref{PR}). It is clear that axial symmetry conditions~(\ref{AxisMask}) are valid.
    $\Diamond$

    Next, we will construct axial symmetric mask $m_0$ with arbitrary order of sum rule.
    The construction will be carried out only for $M=\left(\begin{matrix}
            \gamma_1 & 0\cr
            0 & \gamma_2\cr
      \end{matrix}\right)$,  $\gamma_1,\gamma_2\in\z.$ The construction for $M=\left(\begin{matrix}
            0 & \gamma_1\cr
            \gamma_2 & 0\cr
      \end{matrix}\right)$ is similar.

\begin{theo}
    Let $M=\left(\begin{matrix}
            \gamma_1 & 0\cr
            0 & \gamma_2\cr
      \end{matrix}\right)$, $\gamma_1,\gamma_2\in\z,$ $n\in\n.$
    Then there exists mask $m_0$ that is axial symmetric with
    respect to the origin and has sum rule of order $n.$
    %satisfying condition~(\ref{19}).
    \label{theoAxisM}
\end{theo}

{\bf Proof.} The construction of mask $m_0(\xi)$ will be carried out using the polyphase components
    $\mu_{0l}(\xi),$ $l=0,\dots,m-1.$  The polyphase components should satisfy condition~(\ref{19})
    with some number $\lambda_\gamma,$  $\gamma\in\Delta_n.$
    Also, the above expressions for symmetry
    for the corresponding polyphase components~(\ref{MaskPolyJ})-(\ref{MaskPolyU})
    must be fulfilled according to the type of
    matrix dilation $M$. To cover all cases we consider IV type of dilation matrices.
    Firstly, we set
        $\lambda_{\gamma}=\delta_{\gamma\nul}, \gamma\in\Delta_n.$
        Therefore, by~(\ref{19}) we have the following conditions
    on the polyphase components
    %$\gamma\in\zd_+$,  $[\gamma]\le n$, $\lambda_\nul=1$, such that
        \be
        D^\beta\mu_{0l}(\nul)=\frac {(2\pi i)^{[\beta]}} {\sqrt m}  (-M^{-1}s_l)^{\beta},
        \quad \forall \beta\in\Delta_n, \quad l=0\ddd m-1,
        \label{19Axis}
        \ee
    %To construct such trigonometric polynomial $m_0(\xi)$,
%    firstly we consider $m'_0(\xi)= e^{-2\pi i (c,\xi)} m_0(\xi)$ which is even and axial symmetric, i.e. $m'_0(\xi)=m'_0(-\xi),$ $m'_0(\xi)=m'_0(Y^*\xi)$,
%$m'_0(\xi)=m'_0(E_2\xi)$ and
%    denote $\lambda'_{\alpha}:=\frac 1 {(2\pi i)^{[\alpha]}}D^{\alpha}m'_0(M^{*-1}\xi)\Bigg|_{\xi=\nul}$, $\alpha\in\zd_+,$ $[\alpha]< n$.
%    Obviously, $\lambda'_{\alpha}=0$, if $o(\alpha)$ is not an empty set.
%    Then due to Lemma~\ref{LambdaLem}
%    $$\lambda_{\alpha}=D^{\alpha}m'_0(M^{*-1}\xi)\Bigg|_{\xi=\nul}=
%    \sum\limits_{0\le\gamma\le\alpha} \lambda'_{\gamma} \lll\alpha\atop\gamma\rrr (M^{-1}c)^{\alpha-\gamma}, \alpha\in\zd_+, [\alpha]< n.$$
%
%    We have the following conditions
%    on the polyphase components $\mu_{0k}, $ $ k=0\ddd d_1d_2-1$ as in~(\ref{19mod})
%    %$\gamma\in\zd_+$,  $[\gamma]\le n$, $\lambda_\nul=1$, such that
%        $$
%        D^\beta\mu_{0k}(\nul)=\frac {(2\pi i)^{[\beta]}} {\sqrt m} \sum\limits_{\nul\le\gamma\le\beta}
%        \lambda'_{\gamma} \lll\beta\atop\gamma\rrr (M^{-1}c-M^{-1}s_k)^{\beta-\gamma},$$
%
%    where $\beta\in\zd_+, [\beta]< n. $

For polyphase components $\mu_{0k}$, $k\in K$ with digits $s_k=(v_1, v_2),$
    $v_l=1,\dots, \left\lfloor\frac {d_l-1} 2\right\rfloor,$ $l=1,2,$ we set

    %$$\mu_{0k}(\xi)=\frac {1} {\sqrt m} \sum\limits_{{\beta\in\zd_+}\atop{0\le[\beta]< n}}%\left[\sum\limits_{\nul\le\gamma\le\beta}\lambda'_{\gamma}
    %\lll\beta\atop\gamma\rrr (M^{-1}c-M^{-1}s_k)^{\beta-\gamma} \right]
%        \left[ \frac 1 {(2\pi i)^{[\beta]}} D^\beta\mu_{0k}(\nul)\right]
%        G_{\beta}(\xi)=$$
%        $$\frac {1} {\sqrt m}\sum\limits_{{\beta\in\zd_+}\atop{0\le[\beta]< n}}
%        \left[\sum\limits_{\nul\le\gamma\le\beta \atop o(\gamma)=\emptyset}\lambda'_{\gamma} \lll\beta\atop\gamma\rrr (M^{-1}c-M^{-1}s_k)^{\beta-\gamma} \right]
%        G_{\beta}(\xi)$$

            $$\mu_{0k}(\xi)=\frac {1} {\sqrt m}\sum\limits_{\beta\in\Delta_n}  (-M^{-1}s_k)^{\beta}G_{\beta}(\xi),\quad$$
    where $G_{\beta}(\xi)\in\Theta_{\beta,n},$  $\beta\in\Delta_n,$ are trigonometric polynomials.
    The polyphase components with digits $s_{k'}=(v_1, -v_2),$ $s_{k''}=(-v_1, v_2),$ $s_{k'''}=(-v_1, -v_2)$
    by symmetry conditions~(\ref{MaskPolyK})
    are as follows

           $$\mu_{0k'}(\xi)=\mu_{0k}(-Y^*\xi), \quad \mu_{0k''}(\xi)=\mu_{0k}(Y^*\xi), \quad \mu_{0k'''}(\xi)=\mu_{0k}(-\xi).$$
    Then condition~(\ref{19Axis}) for $\mu_{0k}$ is obviously valid. It follows from

        $$  - s_k =  s_{k'''},\,\,-(s_k)_2= (s_{k'})_2,\,\,-(s_k)_1=(s_{k''})_1$$
 %   $$D^{\beta}\mu_{0k}({\bf 0})=D^{\beta}\mu_{0k'''}({\bf 0}) (-1)^{[\beta]}, \quad D^{\beta}\mu_{0k'}({\bf 0})=D^{\beta}\mu_{0k''}({\bf 0}) (-1)^{[\beta]},$$
%    $$ D^{\beta}\mu_{0k}({\bf 0})=D^{\beta}\mu_{0k'}({\bf 0}) (-1)^{\beta_1}, \quad D^{\beta}\mu_{0k''}({\bf 0})=D^{\beta}\mu_{0k'''}({\bf 0}) (-1)^{\beta_1},$$
%    $$D^{\beta}\mu_{0k}({\bf 0})=D^{\beta}\mu_{0k''}({\bf 0}) (-1)^{\beta_2}, \quad D^{\beta}\mu_{0k'}({\bf 0})=D^{\beta}\mu_{0k'''}({\bf 0}) (-1)^{\beta_2}.$$
%    $$ D^{\beta}G_{\beta}({\bf 0})=(-1)^{\beta_1} D^{\beta}G_{\beta}(Y^*\xi)\Bigg|_{\xi=\nul} , \quad
%    D^{\beta}G_{\beta}({\bf 0})=(-1)^{\beta_2} D^{\beta}G_{\beta}(-Y^*\xi)\Bigg|_{\xi=\nul} $$
%    $$D^{\beta}G_{\beta}({\bf 0})= (-1)^{[\beta]} D^{\beta}G_{\beta}(-\xi)\Bigg|_{\xi=\nul}, $$
that condition~(\ref{19Axis}) holds for  $\mu_{0k'}, \mu_{0k''},  \mu_{0k'''}$.

For polyphase components $\mu_{0r}$, $r\in R,$ with digits $s_r=(v_1, 0),$
    $v_1=1,\dots, \left\lfloor\frac {d_1-1} 2\right\rfloor,$ condition~(\ref{19Axis}) looks as follows

        $$D^\beta\mu_{0r}(\nul)=\frac {(2\pi i)^{[\beta]}} {\sqrt m}
        \left(-\frac {v_1} {d_1}, 0\right)^{\beta}.$$  %(M^{-1}({\nu_1\atop 0}))^{\beta} \delta_{\beta_2 0}
    Then we set

        $$\mu_{0r}(\xi)=\frac 1 {\sqrt m} \sum\limits_{0\le\beta_1<n}
        \left(-\frac {v_1} {d_1}\right)^{\beta_1}g_{\beta_1}(\xi_1),$$
    where $g_{\beta_1}(\xi_1)\in\Theta_{\beta_1,n},$ $0\le\beta_1<n,$ are trigonometric polynomials.
    The polyphase components with digits $s_{r'}=(-v_1, 0)$ by symmetry conditions~(\ref{MaskPolyR})
    are as follows

        $$\quad \mu_{0r'}(\xi)=\mu_{0r}(-\xi).$$
     Then condition~(\ref{19Axis}) is obviously valid for $\mu_{0r}$ and $\mu_{0r'}$, $r\in R,$ $r'\in R'$.

For polyphase components $\mu_{0q}$, $q\in Q,$ with digits $s_q=(0, v_2),$
    $v_2=1,\dots, \left\lfloor\frac {d_2-1} 2\right\rfloor, $ condition~(\ref{19Axis}) looks as follows

        $$D^\beta\mu_{0q}(\nul)=\frac {(2\pi i)^{[\beta]}} {\sqrt m}
        \left(0, -\frac {v_2} {d_2}\right)^{\beta}.$$  %(M^{-1}({\nu_1\atop 0}))^{\beta} \delta_{\beta_2 0}
Then we set

     $$\mu_{0q}(\xi)=\frac 1 {\sqrt m}\sum\limits_{0\le\beta_2<n}
     \left(-\frac {v_2} {d_2}\right)^{\beta_2}g_{\beta_2}(\xi_2),$$
     where $g_{\beta_2}(\xi_2)\in\Theta_{\beta_2,n},$ $0\le\beta_2<n,$ are trigonometric polynomials.
     The polyphase components with digits $s_{q'}=(0, -v_2)$ by symmetry conditions~(\ref{MaskPolyQ})
     are as follows

        $$\mu_{0q'}(\xi)=\mu_{0q}(-\xi).$$
     Then condition~(\ref{19Axis}) is obviously valid for $\mu_{0q}$ and $\mu_{0q'}$, $q\in Q,$ $q'\in Q'$.

 For polyphase components $\mu_{0t}$, $t\in T$ with digits $s_t=(\frac {d_1} 2, v_2),$
     $v_2=1,\dots,\left\lfloor\frac {d_2-1} 2\right\rfloor$ condition~(\ref{19Axis}) looks as follows
         $$D^\beta\mu_{0t}(\nul)=\frac {(2\pi i)^{[\beta]}} {\sqrt m}
        \left(-\frac 12, -\frac {\nu_2} {d_2}\right)^{\beta}.$$
Then we set
        $$\mu_{0t}(\xi)=e^{-2\pi i \frac {\xi_1} 2 } g^A_0(\xi_1)\sum\limits_{0\le\beta_2< n}
        \left(-\frac {\nu_2} {d_2}\right)^{\beta_2}g_{\beta_2}(\xi_2)$$
    where $g_{\beta_2}(\xi_2)\in\Theta_{\beta_2,n},$ $0\le\beta_2<n,$ are trigonometric polynomials
    and $g^A_0(\xi_1)\in\Theta^A_{0,n}$ is trigonometric polynomial of semi-integer degrees associated with $\sigma=\frac 12.$
    The polyphase components with digits $s_{t'}=(\frac {d_1} 2, -v_2)$ by symmetry conditions~(\ref{MaskPolyT})  are as follows
        $$\mu_{0t'}(\xi) = \mu_{0t}(-\xi)e^{-2\pi i \xi_1}.$$
    The condition~(\ref{19Axis}) is obviously valid  for $\mu_{0t}$ and $\mu_{0t'}$, $t\in T,$ $t'\in T'$.

For polyphase components $\mu_{0u}$, $u\in U$  with digits  $s_u=(v_1, \frac {d_2} 2),$
     $v_1=1,\dots, \left\lfloor\frac {d_1-1} 2\right\rfloor$ condition~(\ref{19Axis}) looks as follows
         $$D^\beta\mu_{0u}(\nul)=\frac {(2\pi i)^{[\beta]}} {\sqrt m}
        \left(-\frac {\nu_1} {d_1}, -\frac 12 \right)^{\beta}.$$
    Then we set
        $$\mu_{0u}(\xi)=\sum\limits_{{\bf 0}\le\beta_1\le n}
        \left(-\frac {\nu_1} {d_1}\right)^{\beta_1}g_{\beta_1}(\xi_1)\, e^{-2\pi i \frac {\xi_2} 2 } g^A_0(\xi_2),$$
    where $g_{\beta_1}(\xi_1)\in\Theta_{\beta_1,n}$ $0\le\beta_1<n,$  are trigonometric polynomials
    and $g^A_0(\xi_2)\in\Theta^A_{0,n}$ is trigonometric polynomial of semi-integer degrees associated with $\sigma=\frac 12.$
    The polyphase components with digits $s_{u'}=(-v_1, \frac {d_2} 2)$ by symmetry conditions~(\ref{MaskPolyU})  are as follows
        $$\mu_{0u'}(\xi) = \mu_{0u}(-\xi)e^{-2\pi i \xi_2}.$$
    Then condition~(\ref{19Axis}) is obviously valid for $\mu_{0u}$ and $\mu_{0u'}$, $u\in U,$ $u'\in U'$.

    For polyphase components $\mu_{0{j_p}}$, $p=1,\dots,4$ with digits $ s_{j_1}=(0, 0),$ $s_{j_2}=\left(\frac {d_1} 2,0\right),$
    $s_{j_3}=\left(0, \frac {d_2} 2\right),$ $s_{j_4}=\left(\frac {d_1} 2, \frac {d_2} 2\right)$
     condition~(\ref{19Axis}) looks as follows
        $$D^\beta\mu_{0j_1}(\nul)=\frac {(2\pi i)^{[\beta]}} {\sqrt m}  \delta_{\nul\beta},\quad
        D^\beta\mu_{0j_2}(\nul)=\frac {(2\pi i)^{[\beta]}} {\sqrt m}  \left(-\frac 12, 0\right)^{\beta},$$
        $$D^\beta\mu_{0j_3}(\nul)=\frac {(2\pi i)^{[\beta]}} {\sqrt m}  \left(0, -\frac 12\right)^{\beta}, \quad
       D^\beta\mu_{0j_4}(\nul)=\frac {(2\pi i)^{[\beta]}} {\sqrt m}  \left(-\frac 12, -\frac 12\right)^{\beta}.$$
        Then we set
        $$\mu_{0j_1}(\xi)=\frac 1  {\sqrt m},\quad
        \mu_{0j_2}(\xi)=\frac 1 {\sqrt m} e^{-2\pi i \frac {\xi_1} 2}g^A_0(\xi_1),$$
        $$\mu_{0j_3}(\xi)=\frac 1 {\sqrt m} e^{-2\pi i \frac {\xi_2} 2}g^A_0(\xi_2),\quad
        \mu_{0j_4}(\xi)=\frac 1 {\sqrt m} e^{-2\pi i \frac {\xi_1+\xi_2} 2}g^A_{\nul}(\xi).$$
        where $g^A_0(\xi_i)\in\Theta^A_{0,n},$ $i=1,2,$ are trigonometric polynomials of semi-integer degrees associated with $\sigma=\frac 12$
         and $g^A_{\nul}(\xi)\in\Theta^A_{\nul,n}$ is trigonometric polynomials of semi-integer degrees
         associated with $\sigma=(\frac 12,\frac 12).$  Then condition~(\ref{19Axis}) is obviously valid for $\mu_{0j_p},$ $p=1,2,3,4$.
$\Diamond$

\begin{theo}
Let matrix dilation $M$ be as in~(\ref{AxisM}), $c\in\frac 12\zd,$ mask $m_0$ be axial symmetric
    with respect to the center $c$
    and has sum rule of order $n$.
    %satisfying condition~(\ref{19}).
    Then there exist mask $\w m_0$ that is axial symmetric
    with respect to the center $c$ and satisfies condition~(\ref{20_new}).
    \label{DualTheoAxis}
\end{theo}

The proof remains almost the same as %for point symmetry case
    Theorem's~\ref{DualTheo} proof but
    $\rho_{\beta}=0$, if $o(\beta)\neq\emptyset$, i.e. if $\beta_1$ or $\beta_2$ is odd. Also,
    instead of functions $G^{A}_{\beta}(\xi),$ $\beta\in\Delta^e_n,$ in~(\ref{fDualMaskA}) we should take
    $G^{AxSym}_{\beta}(\xi)$ for $o(\beta)=\emptyset$ that are axial symmetric trigonometric polynomial of
    semi-integer degree associated with $\sigma=c-\lfloor c \rfloor$ such that $G^{AxSym}_{\beta}(\xi)\in\Theta_{\beta,n}.$
    They can be easily constructed as follows
        $$G^{AxSym}_{\beta}(\xi)=\frac 14 (G_{\beta}(\xi)+G_{\beta}(Y^*\xi)+G_{\beta}(-Y^*\xi)+G_{\beta}(-\xi)),
        \quad \forall \beta\in\Delta_n,\quad o(\beta)=\emptyset,$$
    where $G_{\beta}(\xi)\in\Theta_{\beta,n},$ $\beta\in\Delta_n,$ $o(\beta)=\emptyset,$ are trigonometric polynomials of
    semi-integer degree associated with $\sigma=c-\lfloor c \rfloor.$

\begin{theo}Let matrix dilation $M$ be as in~(\ref{AxisM}), % $c\in\frac 12\zd,$
    $n\in\n.$ Masks $m_0$ and $\w m_0$
    are axial symmetric with
    respect to the origin, %center $c,$
    mask $m_0$
    has sum rule of order $n,$
    %satisfy condition~(\ref{19}),
    mask $\w m_0$ satisfies
    condition~(\ref{20_new}).
    Then there exist wavelet masks $m_{\nu}$ and $\w m_{\nu}$, $\nu=1,\dots,m,$
    which are axial symmetric/antisymmetric  with
    respect to some centers and
    %provide approximation order $n$ for the corresponding frame-like wavelet system.
    masks $\w m_{\nu},$ $\nu=1,\dots,m,$  have vanishing moments of order $n$.
\label{theoWaveAxis}
\end{theo}

{\bf Proof.} Matrix extension providing vanishing moments of order $n$ for $\w m_{\nu},$ $\nu=1,\dots,m,$
    realises as in~(\ref{NNMatr0}) and~(\ref{NNMatr}).
    We show that after some transformations axial symmetric/antisymmetric wavelet masks
    can be received.
    Denote matrices elements and matrices rows as in Theorem's~\ref{theoWaveA} proof.
    %by $P,$ $\w P$ the first row of ${\cal  M},$ $\w{\cal  M}$, and by $Q_{0},\dots,Q_{m-1},$ $\w Q_{0},\dots,\w Q_{m-1}$ the second
    %and the others respectively.
    To cover all cases we consider IV type of dilation matrices.
    Due to the symmetry of masks $m_0$ and $\w m_0$ conditions~(\ref{MaskPolyJ})-(\ref{MaskPolyU})
    are satisfied for the corresponding polyphase components $\mu_{0,s}(\xi)$ and $\w\mu_{0,s}(\xi),$ $s=0,\dots,m-1.$

Let us consider the rows $\w Q_{j_p},$ $p=1, 2, 3, 4.$
The first $m$ elements of the rows $\w Q_{j_p}$ look as follows

        $$\w \mu_{j_pk}(\xi)=-\w \mu_{0k}(\xi)\overline{\mu_{0 j_p}(\xi)},\quad
        \w \mu_{j_pk'}(\xi)=-\w \mu_{0k'}(\xi)\overline{\mu_{0 j_p}(\xi)},\quad k\in K,\,k'\in K';$$
        $$\w \mu_{j_pk''}(\xi)=-\w \mu_{0k''}(\xi)\overline{\mu_{0 j_p}(\xi)},\quad
        \w \mu_{j_pk'''}(\xi)=-\w \mu_{0k'''}(\xi)\overline{\mu_{0 j_p}(\xi)},\quad k''\in K'',\,k'''\in K''';$$
        $$\w \mu_{j_pq}(\xi)=-\w \mu_{0q}(\xi)\overline{\mu_{0 j_p}(\xi)},\quad
        \w \mu_{j_pq'}(\xi)=-\w \mu_{0q'}(\xi)\overline{\mu_{0 j_p}(\xi)},\quad q\in Q,\,q'\in Q';$$
        $$\w \mu_{j_pr}(\xi)=-\w \mu_{0r}(\xi)\overline{\mu_{0 j_p}(\xi)},\quad
        \w \mu_{j_pr'}(\xi)=-\w \mu_{0r'}(\xi)\overline{\mu_{0 j_p}(\xi)},\quad r\in R,\,r'\in R';$$
        $$\w \mu_{j_pt}(\xi)=-\w \mu_{0t}(\xi)\overline{\mu_{0 j_p}(\xi)},\quad
        \w \mu_{j_pt'}(\xi)=-\w \mu_{0t'}(\xi)\overline{\mu_{0 j_p}(\xi)},\quad t\in T,\,t'\in T';$$
        $$\w \mu_{j_pu}(\xi)=-\w \mu_{0u}(\xi)\overline{\mu_{0 j_p}(\xi)},\quad
        \w \mu_{j_pu'}(\xi)=-\w \mu_{0u'}(\xi)\overline{\mu_{0 j_p}(\xi)},\quad u\in U,\,u'\in U';$$
        $$\w \mu_{j_pj_r}(\xi)=-\w \mu_{0j_r}(\xi)\overline{\mu_{0 j_p}(\xi)}, r=1,2,3,4, r\neq p;$$
        $$\w \mu_{j_pj_p}(\xi)=1-\w \mu_{0j_p}(\xi)\overline{\mu_{0 j_p}(\xi)}.$$
        %It is clear that for $\mu_{j_ps},$ $s\in\{Q,Q',R,R',T,T',U,U'\}$ holds
         %$$\mu_{j_1s}(\xi)=\mu_{j_1s'}(-\xi)$$.
         Due to the fact that

            \be
            \mu_{0j_1}(\xi)=\mu_{0j_1}(-\xi),
            \quad\mu_{0j_1}(M^*\xi)=\mu_{0j_1}(M^*Y^*\xi),
            \label{fAxisJ_1}
            \ee
         it is clear that conditions~(\ref{MaskPolyJ})-(\ref{MaskPolyU}) are valid for $\w\mu_{j_1s}(\xi),$ $s=0,\dots,m-1.$
         These elements are the polyphase components of wavelet mask $\w m_{j_1}(\xi)$.
         Therefore, by Lemma~\ref{lemAxisMask} $\w m_{j_1}(\xi)$ is axial symmetric with respect to the origin.
         Using the same considerations and
                  \be
                  \mu_{0j_2}(M^*\xi)=\mu_{0j_2}(-M^*\xi)e^{-2\pi i d_1\xi_1},\quad \mu_{0j_2}(M^*\xi)=\mu_{0j_2}(M^*Y^*\xi)e^{-2\pi i d_1\xi_1}
                  \label{fAxisJ_2}
                  \ee
            %\nonumber
            %\\
            $$\mu_{0j_3}(M^*\xi)=\mu_{0j_3}(-M^*\xi)e^{-2\pi i  d_2\xi_2},\quad \mu_{0j_3}(M^*\xi)=\mu_{0j_3}(M^*Y^*\xi)$$
            %\nonumber
            %\\
            $$
            \mu_{0j_4}(M^*\xi)=\mu_{0j_4}(-M^*\xi)e^{-2\pi i (d,\xi)},\quad\mu_{0j_4}(M^*\xi)=\mu_{0j_4}(M^*Y^*\xi)e^{-2\pi i d_1\xi_1}$$
          we can conclude, that
         $\w m_{j_2}$ is axial symmetric with respect to the center $\left(\frac {d_1} 2,0\right)$;
         $\w m_{j_3}$ is axial symmetric with respect to the center $\left(0,\frac {d_2} 2\right)$;
         $\w m_{j_4}$ is axial symmetric with respect to the center $\frac d 2$.

Let us consider the rows $\w Q_k, \w Q_{k'}, \w Q_{k''}, \w Q_{k'''},$
for $k\in K$, $k'\in K'$, $k''\in K''$, $k'''\in K'''$ such that

            $$s_k=(v_1, v_2),\,  s_{k'}=(v_1, -v_2),\,  s_{k''}=(-v_1, v_2),\,
            s_{k'''}=(-v_1, -v_2),$$
                where $  v_1=1,\dots, \left\lfloor\frac {d_1-1} 2\right\rfloor,
    v_2=1,\dots, \left\lfloor\frac {d_2-1} 2\right\rfloor.$
    Transform them into

        $$\w Q'_{k}:=\frac 1 4 (\w Q_k+\w Q_{k'}+\w Q_{k''}+\w Q_{k'''}),\quad
        \w Q'_{k'}:=\frac 1 4 (\w Q_k-\w Q_{k'}-\w Q_{k''}+\w Q_{k'''}),$$
         $$\w Q'_{k''}:=\frac 1 4 (\w Q_k+\w Q_{k'}-\w Q_{k''}-\w Q_{k'''}),\quad
         \w Q'_{k'''}:=\frac 1 4 (\w Q_k-\w Q_{k'}+\w Q_{k''}-\w Q_{k'''}).$$
    Denote

          $$A(\xi)=\overline{\mu_{0k}(\xi)}+\overline{\mu_{0k'}(\xi)}+\overline{\mu_{0k''}(\xi)}+\overline{\mu_{0k'''}(\xi)},\quad
          B(\xi)=\overline{\mu_{0k}(\xi)}-\overline{\mu_{0k'}(\xi)}-\overline{\mu_{0k''}(\xi)}+\overline{\mu_{0k'''}(\xi)},$$

          $$C(\xi)=\overline{\mu_{0k}(\xi)}-\overline{\mu_{0k'}(\xi)}+\overline{\mu_{0k''}(\xi)}-\overline{\mu_{0k'''}(\xi)},\quad
          D(\xi)=\overline{\mu_{0k}(\xi)}+\overline{\mu_{0k'}(\xi)}-\overline{\mu_{0k''}(\xi)}-\overline{\mu_{0k'''}(\xi)}.$$
    Therefore, by~(\ref{MaskPolyK})

            $$A(\xi)=A(-\xi),\quad A(M^*\xi)=A(M^*Y^*\xi),\quad B(\xi)=B(-\xi),\quad B(M^*\xi)=-B(M^*Y^*\xi),$$
            $$C(\xi)=-C(-\xi),\quad C(M^*\xi)=C(M^*Y^*\xi),\quad D(\xi)=-D(-\xi),\quad D(M^*\xi)=-D(M^*Y^*\xi).$$
          Consider the first $m$ elements
          %polyphase compomets
          of the row $\w Q'_{k}.$
          Due to the fact that $A(\xi)$ has the same symmetry properties as $\mu_{0j_1}$ in~(\ref{fAxisJ_1}),
          we have almost the same situation as for the row $\w Q_{j_1}$
          except the elements
          %we can receive axial symmetric
          %wavelet mask $m_k$.
  %        $$\w Q'_{k}=
%        (\dots,-\widetilde\mu_{0{j_1}}A, -\widetilde\mu_{0q}A,\dots, -\widetilde\mu_{0q'}A,\dots,
%        -\widetilde\mu_{0r}A,\dots, -\widetilde\mu_{0r'}A,\dots,
%        ,\dots,
%        -\widetilde\mu_{0u}A,\dots, -\widetilde\mu_{0u'}A,\dots,
%        -\widetilde\mu_{0t}A,\dots,$$ $$ -\widetilde\mu_{0t'}A,\dots,
%        \frac 14-\frac 14\widetilde\mu_{0k}A,\dots,
%        \frac 14-\frac 14\widetilde\mu_{0k'}A,\dots,
%        \frac 14-\frac 14\widetilde\mu_{0k''}A,\dots,
%        \frac 14-\frac 14\widetilde\mu_{0k'''}A,\dots,).$$

        %It is clear that for the elements
        $$\w \mu_{kk}(\xi)=\frac 14-\frac 14\widetilde\mu_{0k}(\xi)A(\xi),\quad
        \w \mu_{kk'}(\xi)=\frac 14-\frac 14\widetilde\mu_{0k'}(\xi)A(\xi),\quad$$
        $$\w \mu_{kk''}(\xi)=\frac 14-\frac 14\widetilde\mu_{0k''}(\xi)A(\xi),\quad
        \w \mu_{kk'''}(\xi)=\frac 14-\frac 14\widetilde\mu_{0k'''}(\xi)A(\xi), \quad\w \mu_{kj_1}=-\w \mu_{0j_1}(\xi)A(\xi).$$
        It is clear that those elements also satisfy the conditions~(\ref{MaskPolyJ}) and~(\ref{MaskPolyK}).
        Then wavelet mask  $\w m_{k}$ is axial symmetric with respect to the origin.
%        The same is valid for
%        $\mu_{kq},
%        \mu_{kq'},$
%        $\mu_{kr},
%        \mu_{kr'},$
%        $\mu_{kt},
%        \mu_{kt'},$
%        $\mu_{ku},
%        \mu_{ku'},$
%        $\mu_{kj_p}, p=1,2,3,4.$
        With the same reasoning from the other rows $\w Q'_{k'},$  $\w Q'_{k''},$ $\w Q'_{k'''},$
        we receive wavelet masks  $\w m_{k'}$, $\w m_{k''}$, $\w m_{k'''}$ that are axial antisymmetric as follows
            $$\w m_{k'}(\xi)={\w m_{k'}(-\xi)},
              \quad \w m_{k'}(\xi)=-\w m_{k'}(Y^* \xi).$$
              $$\w m_{k''}(\xi)=-{\w m_{k''}(-\xi)},
              \quad \w m_{k''}(\xi)=\w m_{k''}(Y^* \xi).$$
              $$\w m_{k'''}(\xi)=-{\w m_{k'''}(-\xi)},
              \quad \w m_{k'''}(\xi)=-\w m_{k'''}(Y^* \xi).$$

Let us consider the rows $\w Q_r, \w Q_{r'}$ for $r\in R$, $r'\in R'$
    such that $s_r=(v_1, 0),\,  s_{r'}=(-v_1, 0),$ where $v_1=1,\dots, \left\lfloor\frac {d_1-1} 2\right\rfloor.$
    Transform them into

        $$\w Q'_{r}:=\frac 1 2 (\w Q_r+\w Q_{r'}),\quad \w Q'_{r'}:=\frac 1 2 (\w Q_r-\w Q_{r'}).$$
        Denote,

            $$A(\xi)=\frac 12\left(\overline{\mu_{0r}(\xi)}+\overline{\mu_{0r'}(\xi)}\right),\quad
            B(\xi)=\frac12 \left(\overline{\mu_{0r}(\xi)}-\overline{\mu_{0r'}(\xi)}\right).$$
    Therefore,  by~(\ref{MaskPolyR})

            $$A(\xi)=A(-\xi),\quad A(M^*\xi)=A(M^*Y^*\xi),\quad
               B(\xi)=-B(-\xi),\quad B(M^*\xi)=-B(M^*Y^*\xi).$$
    Due to the fact that $A(\xi)$ has the same symmetry properties as $\mu_{0j_1}$ in~(\ref{fAxisJ_1}),
          we have almost the same situation as for the row $\w Q_{j_1}$
          except the elements
          %Here we have almost the same situation as for $m_{j_1},$
                %except the elements
                $$\w \mu_{r j_1}=-\w \mu_{0j_1}(\xi)A(\xi),\quad
                %\w \mu_{r' j_1}=-\w \mu_{0j_1}(\xi)B(\xi),$$
                \w \mu_{r r}=\frac 12-\w \mu_{0r}(\xi)A(\xi),\quad
                \w \mu_{r r'}=\frac 12-\w \mu_{0r'}(\xi)A(\xi),$$
                %$$\w \mu_{r' r}=\frac 12-\w \mu_{0r}(\xi)B(\xi),\quad
                %\w \mu_{r' r'}=-\frac 12-\w \mu_{0r'}(\xi)B(\xi),$$

                It is clear that those elements also satisfy conditions~(\ref{MaskPolyR}) and~(\ref{MaskPolyJ}).
            Then wavelet mask $\w m_{r}$ is axial symmetric with respect to the origin.
                With the same reasoning from the other row $\w Q'_{r'}$
                 we receive wavelet mask $\w m_{r'}$ that is axial antisymmetric as follows
                $$\w m_{r'}(\xi)=-{\w m_{r'}(-\xi)},
              \quad \w m_{r'}(\xi)=-\w m_{r'}(Y^* \xi).$$

Let us consider the rows $\w Q_q, \w Q_{q'}$ for $q\in Q$, $q'\in Q'$
    such that  $s_q=(0, v_2),\,s_{q'}=(0, -v_2),$ where $ v_2=1,\dots, \left\lfloor\frac {d_2-1} 2\right\rfloor.$
    Transform them into

            $$\w Q'_{q}:=\frac 1 2 (\w Q_q+\w Q_{q'}),\quad
            \w Q'_{q'}:=\frac 1 2 (\w Q_q-\w Q_{q'}).$$
    Due to the same reasoning as for the rows $\w Q_r, \w Q_{r'}$ we conclude
    that wavelet mask $\w m_q$ is axial symmetric with respect to the origin,
    wavelet mask $\w m_{q'}$ is axial antisymmetric as follows
            $$\w m_{q'}(\xi)=-{\w m_{q'}(-\xi)},
            \quad \w m_{q'}(\xi)=\w m_{q'}(Y^* \xi).$$

Let us consider the rows  $\w Q_t, \w Q_{t'}$ for $t\in T$, $t'\in T'$ such that
        $s_t=\left(\frac {d_1} 2, v_2\right),\,s_{t'}=\left(\frac {d_1} 2, -v_2\right),$
        where $ v_2=1,\dots, \left\lfloor\frac {d_2-1} 2\right\rfloor.$
    Transform them into

         $$\w Q'_{t}:=\frac 1 2 (\w Q_t+\w Q_{t'}),\quad
         \w Q'_{t'}:=\frac 1 2 (\w Q_t-\w Q_{t'}).$$
    Denote,

        $$A(\xi)=\frac 12\left(\mu_{0t}(\xi)+\mu_{0t'}(\xi)\right),\quad
        B(\xi)=\frac 12\left(\mu_{0t}(\xi)-\mu_{0t'}(\xi)\right).$$
    Therefore,  by~(\ref{MaskPolyT})
        $$A(M^*\xi)=e^{-2\pi i d_1\xi_1}A(-M^*\xi),\quad A(M^*\xi)=e^{-2\pi i d_1\xi_1}A(M^*Y^*\xi),$$
     $$B(M^*\xi)=-e^{-2\pi i d_1\xi_1}B(-M^*\xi),\quad B(M^*\xi)=e^{-2\pi i d_1\xi_1}B(M^*Y^*\xi).$$
     Due to the fact that $A(\xi)$ has the same symmetry properties as $\mu_{0j_2}$ in~(\ref{fAxisJ_2}),
          we have almost the same situation as for the row $\w Q_{j_2}$
          except the elements
          %Here we have almost the same situation as for $m_{j_2},$
                %except the elements

                $$\w \mu_{t j_2}(\xi)=-\w \mu_{0j_2}(\xi)\overline{A(\xi)}, \quad
                %\mu_{t' j_2}=-\w \mu_{0j_2}(\xi)\overline{B(\xi)},$$
                \w \mu_{t t}(\xi)=\frac 12-\w \mu_{0t}(\xi)\overline{A(\xi)}, \quad
                \w \mu_{t t'}(\xi)=\frac 12-\w \mu_{0t'}(\xi)\overline{A(\xi)}.$$
               % $$\mu_{t' t}=\frac 12-\w \mu_{0t}(\xi)\overline{B(\xi)},
                %\mu_{t' t'}=-\frac 12-\w \mu_{0t'}(\xi)\overline{B(\xi)}.$$

                It is clear that those elements also satisfy
                conditions~(\ref{MaskPolyJ}) and~(\ref{MaskPolyT}).
    Then wavelet mask $\w m_t$ is axial symmetric with respect to the the center $\left(\frac {d_1} 2,0\right)$,
               With the same reasoning from the other row $\w Q'_{t'}$ we receive
                wavelet mask $\w m_{t'}$ that is  axial antisymmetric as follows
                $$\w m_{t'}(\xi)=-e^{2\pi i ((d_1,0), \xi)}{\w m_{t'}(-\xi)},
              \quad \w m_{t'}(\xi)=e^{2\pi i ((d_1,0), \xi)}\w m_{t'}(Y^* \xi).$$

Let us consider the rows  $\w Q_u, \w Q_{u'}$ for $u\in U$, $u\in U$ such that
    $ s_u=\left(v_1, \frac {d_2} 2\right),\, s_{u'}=\left(-v_1, \frac {d_2} 2\right),$
     where $ v_1=1,\dots, \left\lfloor\frac {d_1-1} 2\right\rfloor.$
    Transform them into
        $$\w Q'_{u}:=\frac 1 2 (\w Q_u+\w Q_{u'}),
        \quad\w Q'_{u'}:=\frac 1 2 (\w Q_u-\w Q_{u'}).$$
        Due to the same reasoning as for the rows $\w Q_t, \w Q_{t'}$ we conclude
    that wavelet mask $\w m_u$ is axial symmetric with respect to the center $(0,\frac {d_2} 2)$,
    wavelet mask $\w m_{u'}$ is axial antisymmetric as follows
                $$\w m_{u'}(\xi)=-e^{2\pi i ( (0,d_2), \xi)}{\w m_{u'}(-\xi)},
              \quad \w m_{u'}(\xi)=-\w m_{u'}(Y^* \xi).$$

For the rows $Q_{l}$ of matrix ${\cal N}$ we make analogous transformations
    as for $\w Q_{l},$ $l=0,\dots,m-1.$
    For instance,
            $$ Q'_{u}:= Q_u+ Q_{u'},
        \quad Q'_{u'}:= Q_u- Q_{u'},\quad \forall u\in U, u'\in U'.$$
    New rows $Q'_{l}$
    form wavelet masks $m_{l}$ that are axial symmetric/antisymmetric
    with respect to some points.
    It is clear that the equality ${\cal N}\w{\cal N}^*=I_{m+1}$ remains correct with new rows.$\Diamond$

              Now, we consider the case when symmetry center $c$ is semi-integer, i.e. $c_i\in\frac 12 \z^2, i=1,2.$
    Due to the fact that orders of  sum rule, vanishing moments, linear-phase moments and axial symmetry of trigonometric polynomial
    are invariant with respect to the integer shift it is enough to consider only the case when $c=\left(\frac 12,\frac 12\right).$
     Let us define the set of digits $D(M)$ as

     $$D(M)=\left\{\left(v_1, v_2\right) \Big| \quad
        v_l=-\left\lfloor\frac {d_l-1} 2\right\rfloor,\dots, \left\lfloor\frac {d_l} 2\right\rfloor, l=1,2\right\}.$$
        %v_2=\left\lfloor\frac {-d_2+1} 2\right\rfloor,\dots, \left\lfloor\frac {d_2} 2\right\rfloor
     Enumerate all digits from $0$ to $m-1$ and divide them into several sets with respect to the parity of integers $d_1$ and $d_2$.

          I. If $d_1, d_2$ are even, then
        the set of digits can be divided as follows

              $
              K:=\{ k \Big| s_k=(v_1, v_2),\quad v_l=1,\dots, \left\lfloor\frac {d_l} 2\right\rfloor, l=1,2\},$

              $K':=\{ k' \Big| s_{k'}=(v_1, 1-v_2),\quad v_l=1,\dots, \left\lfloor\frac {d_l} 2\right\rfloor, l=1,2\},$

              $K'':=\{ k'' \Big| s_{k''}=(1-v_1, v_2),\quad v_l=1,\dots, \left\lfloor\frac {d_l} 2\right\rfloor, l=1,2\},$

              $K''':=\{ k''' \Big| s_{k'''}=(1-v_1, 1-v_2),\quad v_l=1,\dots, \left\lfloor\frac {d_l} 2\right\rfloor, l=1,2\};
              $

             II. If $d_1$ is odd, $d_2$ is even, then
        the set of digits can be divided as follows
            \begin{enumerate}
              \item $T:=\{t \Big| s_t=(-\frac {d_1-1} 2, v_2), v_2=1,\dots,
              \left\lfloor\frac {d_2} 2\right\rfloor\};$

              $T':=\{t' \Big| s_{t'}=(-\frac {d_1-1} 2, 1-v_2), v_2=1,\dots,
              \left\lfloor\frac {d_2} 2\right\rfloor\};$
              \item $K, K', K'', K'''$ are the same as in I.
            \end{enumerate}

        III. If $d_1$ is even, $d_2$ is odd, then
       the set of digits can be divided as follows
            \begin{enumerate}
              \item $U:=\{u \Big| s_u=(v_1, -\frac {d_2-1} 2),
              v_1=1,\dots, \left\lfloor\frac {d_1} 2\right\rfloor\};$

              $U':=\{u' \Big| s_{u'}=(1-v_1, -\frac {d_2-1} 2),
              v_1=1,\dots, \left\lfloor\frac {d_1} 2\right\rfloor\};$
              \item $K, K', K'', K'''$ are the same as in I.
            \end{enumerate}

                    IV. If $d_1, d_2$ are odd, then
       the set of digits can be divided as follows
            \begin{enumerate}
              \item $J:=\{j_1 \Big| s_{j_1}=(-\frac {d_1-1} 2, -\frac {d_2-1} 2)\};$
              \item $T, T'$ and $U, U'$ are the same as in II and III;
              \item $K, K', K'', K'''$ are the same as in I.
            \end{enumerate}
     Let us reformulate the axial symmetry conditions~(\ref{AxisMask}) for mask $m_0$ in terms of its polyphase components. To cover all
     cases we give them for the IV type of dilation matrices.

     For the I type of dilation matrices mask $m_0$ can be represented as

            %\ba
            $$m_0(\xi)\sqrt m=
            \sum\limits_{k\in K}\mu_{0k}(M^*\xi)e^{2\pi i (s_k,\xi)}+
             \sum\limits_{k'\in K'}\mu_{0k'}(M^*\xi)e^{2\pi i (s_{k'},\xi)}$$
             %\nonumber
             %\\
             $$\sum\limits_{k''\in K''}\mu_{0k''}(M^*\xi)e^{2\pi i (s_{k''},\xi)}+
             \sum\limits_{k'''\in K'''}\mu_{0k'''}(M^*\xi)e^{2\pi i (s_{k'''},\xi)}:=m_0^I(\xi).$$
            %\label{AxisMaskPolyIS}
            %\ea
      %where for any summand $\mu_{0l}e^{2\pi i \left(a\atop b \right)\xi}$ the set $\Omega_M^l$ include digit $s_l=(a, b).$

For the II type of dilation matrices mask $m_0$ can be represented as
            $m_0(\xi)\sqrt m=m_0^I(\xi)+m_0^{II}(\xi),$
        where

         %\ba
            $$m_0^{II}(\xi)=\sum\limits_{t\in T}
            \mu_{0t}(M^*\xi)e^{2\pi i (s_t,\xi)}+
            \sum\limits_{t'\in T'}
            \mu_{0t'}(M^*\xi)e^{2\pi i (s_{t'}\xi)}.$$
            %\label{AxisMaskPolyIIS}
            %\ea

For the III type of dilation matrices mask $m_0$ can be represented as $m_0(\xi)\sqrt m=m_0^I(\xi)+m_0^{III}(\xi),$ where

         % \ba
            $$m_0^{III}(\xi)=\sum\limits_{u\in U}
            \mu_{0u}(M^*\xi)e^{2\pi i (s_u,\xi)}+
            \sum\limits_{u'\in U'}
            \mu_{0u'}(M^*\xi)e^{2\pi i (s_{u'},\xi)}.$$
            %\label{AxisMaskPolyIIIS}
            %\ea

For the IV type of dilation matrices mask $m_0$ can be represented as

        $$m_0(\xi)\sqrt m=m_0^I(\xi)+m_0^{II}(\xi)+m_0^{III}(\xi)+m_0^{IV}(\xi),$$
         where

          %\ba
          $$  m_0^{IV}(\xi):=
            \mu_{0j_1}(M^*\xi)e^{2\pi i (s_{j_1},\xi)}.$$
            %\label{AxisMaskPolyIVS}
            %\ea

            From conditions~(\ref{AxisMask}) we can conclude the following conditions on
    the polyphase components according to the chosen digits and taking into account
    that  $c=\left(\frac 12,\frac 12\right)$

            $$s_{j_1}=\left(-\frac {d_1-1} 2, -\frac {d_2-1} 2\right):
            \quad \mu_{0j_1}(\xi)=\mu_{0j_1}(-\xi)e^{-2\pi i (2c,\xi)},
            %\ba
            \quad\mu_{0j_1}(M^*\xi)=\mu_{0j_1}(M^*Y^*\xi)e^{-2\pi i d_1\xi_1},$$
            %\label{MaskPolyJS}
            %\ea
            %\nonumber

            $$s_k=(v_1, v_2),\,  s_{k'}=(v_1, 1-v_2),\,  s_{k''}=(1-v_1, v_2),\,
            s_{k'''}=(1-v_1, 1-v_2):\hspace{3.2cm}$$
            $$\mu_{0k}(\xi)=\mu_{0k'''}(-\xi), \quad \mu_{0k'}(\xi)=\mu_{0k''}(-\xi),$$
            %\be
            $$\mu_{0k}(M^*\xi)=\mu_{0k''}(M^*Y^*\xi), \quad\mu_{0k'}(M^*\xi)=\mu_{0k'''}(M^*Y^*\xi),$$
            %\label{MaskPolyKS}
            %\ee

            $$ s_t=\left(-\frac {d_1-1} 2, v_2\right),\,s_{t'}=\left(-\frac {d_1-1} 2, 1-v_2\right):
            \quad \mu_{0t}(M^*\xi)=\mu_{0t'}(-M^*\xi)e^{-2\pi i d_1\xi_1},\hspace{1.6cm}$$
            %\nonumber
            %\\
            %\be
           $$\mu_{0t}(M^*\xi)=\mu_{0t}(M^*Y^*\xi)e^{-2\pi i d_1\xi_1},\quad\mu_{0t'}(M^*\xi)=\mu_{0t'}(M^*Y^*\xi)e^{-2\pi i d_1\xi_1},$$
            %\label{MaskPolyTS}
            %\ee

            %\ba
            $$ s_u=\left(v_1,  -\frac {d_2-1} 2\right),\, s_{u'}=\left(1-v_1,  -\frac {d_2-1} 2\right):\hspace{6.5cm}$$
            %\nonumber
           % \\
             $$\quad\mu_{0u}(M^*\xi)=\mu_{0u'}(-M^*\xi)e^{-2\pi i d_2\xi_2}, \quad\mu_{0u}(M^*\xi)=\mu_{0u'}(M^*Y^*\xi),$$
            %\label{MaskPolyUS}
            %\ea

            where $  v_1=1,\dots, \left\lfloor\frac {d_1} 2\right\rfloor,
            v_2=1,\dots, \left\lfloor\frac {d_2} 2\right\rfloor,
            $

Next, we will construct axial symmetric mask $m_0$ with arbitrary order of sum rule.

\begin{theo}
    Let $M=\left(\begin{matrix}
            \gamma_1 & 0\cr
            0 & \gamma_2\cr
      \end{matrix}\right)$, $\gamma_1,\gamma_2\in\z,$ $n\in\n.$
    Then there exists mask $m_0$
    that is axial symmetric with
    respect to the center $c=(1/2,1/2)$ and has sum rule of order $n.$
    %satisfying condition~(\ref{19}).
    \label{theoAxisMS}
\end{theo}

\begin{theo}Let matrix dilation $M$ be as in~(\ref{AxisM}), % $c\in\frac 12\zd,$
    $n\in\n.$ Masks $m_0$ and $\w m_0$
    are axial symmetric with
    respect to the center $c$, %center $c,$
    mask $m_0$ has sum rule of order $n,$
    %satisfy condition~(\ref{19}),
    mask $\w m_0$ satisfies
    condition~(\ref{20_new}).
    Then there exist wavelet masks  $m_{\nu}$ and $\w m_{\nu}$, $\nu=1,\dots,m,$
    which are axial symmetric/antisymmetric  with
    respect to some centers and
    %provide approximation order $n$ for the corresponding frame-like wavelet system.
    masks $\w m_{\nu},$ $\nu=1,\dots,m,$  have vanishing moments of order $n$.
\label{theoWaveAxisS}
\end{theo}

Proofs of these theorems almost repeat the proofs of Theorem~\ref{theoAxisM} and Theorem~\ref{theoWaveAxis}.

\section{Construction of highly symmetric frame-like wavelets}

This section is devoted to the construction of highly symmetric frame-like wavelet systems
    for some dilation matrices in $\rd.$ Let $H$ be a symmetry group on $\zd.$
    In what follows in this section, we consider dilation matrices $M$ such that symmetry group
    $H$ is a symmetry group with respect to these dilation matrices $M$, i.e.
    $$MEM^{-1}\in H,\quad \forall E\subset H$$
    and also for all digits $s_j\in D(M)$ %$2M^{-1}$ is an integer matrix.
        we have
            \be
            E s_j = Mr_j^{E} + s_j,\quad r_j^{E}\in\zd,\quad \forall E\in H, \quad j=0,\dots,m-1,
            \label{fR_j}
            \ee
            i.e. $Es_j\in \Omega_M^j,$ $j=0,\dots,m-1.$ Let $c$ be in $\zd.$
            With the imposed requirements %on dilation matrices $M$
            the polyphase components of $H$-symmetric mask remains $H$-symmetric.

        %    for a mask $m_0$ which is
        %    $H$-symmetric with respect to the center $c$
       % each of its polyphase components of $H$-symmetric mask $\mu_{0j}(\xi),$ $j=0,\dots,m-1,$ remains $H$-symmetric.

%    For matrix $M$, such that $2M^{-1}$ is an integer matrix, the set of digits $D(M)$
%        consists only of digits $s_j$, $j\in J.$ Indeed, due to Lemma~\ref{PolyLem} for each
%        $k=0,\dots,m-1$ we can find index $l$ such that $M^{-1}(s_k-s_l)\in\zd$ in a one way.
%        It is clear that $2M^{-1}s_k\in\zd$ for all $k=0,\dots,m-1$. Therefore, $I=\emptyset,$
%        $J=\{0,\dots,m-1\}.$

\begin{lem} %$c\in\zd$,
    Let $H$ be a symmetry group on $\zd,$
    $M$ be the dilation matrix, $c\in\zd.$
    %$H$ is a symmetry group with respect to that dilation matrix $M$ and $2M^{-1}$ is an integer matrix.
    A mask $m_0(\xi)$ is $H$-symmetric with respect to the center $c$ if and only if
    its polyphase components $\mu_{0j}(\xi),$ $j=0,\dots,m-1,$ satisfy

                \be
                \mu_{0j}(\xi)=\mu_{0j}(M^*E^*M^{*-1}\xi)e^{2\pi i(r_j^{E},\xi)}e^{2\pi i(M^{-1}(c-Ec),\xi)}, \quad \forall E\in H,
                \label{fPolySymH}
                \ee
       where $r_j^{E}$ are given by~(\ref{fR_j}).
\label{lemHSym}
\end{lem}

 {\bf Proof.} Firstly, we assume that mask $m_0$ is $H$-symmetric with respect to the center $c$, i.e.
             $$e^{2\pi i (-c,\xi)}m_0(\xi)=e^{2\pi i (-Ec,\xi)}m_0(E^*\xi), \quad\forall E\in H.$$
         Define $m'_0(\xi)=e^{2\pi i (-c,\xi)}m_0(\xi).$ It is clear that $m'_0(\xi)$ is $H$-symmetric with respect to the
         origin. The polyphase components of $m'_0(\xi)$ are $\mu'_{0j}(\xi).$ Then, by~(\ref{fR_j})
                $$ m'_0(E^* \xi)= \frac1{\sqrt m}\sum\limits_{j=0}^{m-1} e^{2\pi i(E s_j,\xi)}\mu'_{0j}(M^*E^*\xi) =
                \frac1{\sqrt m}\sum\limits_{j=0}^{m-1} e^{2\pi i(Mr_j^{E} + s_j,\xi)}\mu'_{0j}(M^*E^*\xi)=
                $$
                $$
                \frac1{\sqrt m}\sum\limits_{j=0}^{m-1} e^{2\pi i(s_j,\xi)}\mu'_{0j}(M^*\xi)
                =m'_0(\xi), \quad \forall E\in H.$$
             Hence, we obtain %that for each $j\in \{0,\dots,m-1\}$
                $$\mu'_{0j}(M^*\xi)=\mu'_{0j}(M^*E^*\xi)e^{2\pi i(Mr_j^{E},\xi)},\quad \forall E\in H, \quad \forall j\in \{0,\dots,m-1\}.$$
                %where $r_j^{E}$ from~(\ref{fR_j}), i.e. $r_j^{E}=M^{-1}(E s_j-s_j).$
                Finally,

                $$\mu'_{0j}(\xi)=\mu'_{0j}(M^*E^*M^{*-1}\xi)e^{2\pi i(r_j^{E},\xi)}, \quad \forall E\in H, \quad \forall j\in \{0,\dots,m-1\}.$$
         Due to $m_0(\xi)=e^{2\pi i (c,\xi)}m'_0(\xi),$ the polyphase components of $m_0$ and $m'_0$ are related as follows

            $$\mu_{0j}(M^*\xi)=e^{2\pi i (c,\xi)}\mu'_{0j}(M^*\xi),\quad
         \mu_{0j}(M^*E^*\xi)=e^{2\pi i (Ec,\xi)}\mu'_{0j}(M^*E^*\xi), \quad \forall j\in \{0,\dots,m-1\}.$$
         Thereby,
          $$\mu_{0j}(M^*\xi)=\mu_{0j}(M^*E^*\xi)e^{2\pi i(Mr_j^{E},\xi)}e^{2\pi i(c-Ec,\xi)},\quad \forall E\in H,  \quad \forall j\in \{0,\dots,m-1\}.$$
         If $c\in\Omega_M^j$ for some  $j\in \{0,\dots,m-1\},$ then $Ec\in\Omega_M^j$. It is clear that $M^{-1}(c-Ec)\in\zd.$
         Therefore,~(\ref{fPolySymH}) is valid.

      The proof of the conversely statement is obvious. $\Diamond$

       \begin{theo} Let $M$ be the matrix dilation, $n\in\n,$ $c\in\zd,$
       $H$ be a symmetry group on $\zd.$ Masks $m_0$ and $\w m_0$
    are $H$-symmetric with
    respect to the center $c$, mask $m_0$
    has sum rule of order $n,$
    %satisfy condition~(\ref{19}),
    mask $\w m_0$ satisfies
    condition~(\ref{20_new}). Then there exist wavelet masks  $m_{\nu}$ and $\w m_{\nu}$, $\nu=1,\dots,m,$
    which are $H$-symmetric with
    respect to some points and
    %provide approximation order $n$ for the corresponding frame-like wavelet system.
    wavelet masks $\w m_{\nu},$ $\nu=1,\dots,m,$  have vanishing moments of order $n$.
\label{theoWaveHSym}
\end{theo}

{\bf Proof.} Matrix extension providing vanishing moments of order $n$ for $\w m_{\nu},$ $\nu=1,\dots,m,$
    realises as in~(\ref{NNMatr0}) and~(\ref{NNMatr}).
        We denote matrices elements as in Theorem's~\ref{theoWaveA} proof

        $${\cal N}=\{\mu_{kn}\}_{k,n=0}^m,\quad
        {\cal N}=\{\w\mu_{kn}\}_{k,n=0}^m,$$
        the first rows of matrices ${\cal N},$ ${\cal \w N}$ by $P,$ $\w P$ and the other rows by
        $Q_{\nu},$ $\w Q_{\nu},$ $\nu=0,\dots,m-1.$

        Due to $H$-symmetry of masks $m_0(\xi)$ and $\w m_0(\xi),$ by Lemma~\ref{lemHSym}~(\ref{fPolySymH}) is valid
        for its polyphase components $\mu_{0j}(\xi)$ and $\widetilde\mu_{0j}(\xi),$ $j=0,\dots,m-1,$ respectively. %we have~(\ref{fPolySymH}) for all $j=0,\dots,m-1.$

        Let us consider the row $\w Q_{j},$ $j\in \{0,\dots,m-1\},$

         $$\w Q_{j}=(\dots,-\widetilde\mu_{0l}\overline\mu_{0j},\hdots,1-\widetilde\mu_{0j}\overline\mu_{0j},
        \dots),\quad l\in\{0,\dots,m-1\},\,\, l\neq j.$$
        It is clear that

            $$\w\mu_{j l}(M^*E^*M^{*-1}\xi)=-\widetilde\mu_{0l}(M^*E^*M^{*-1}\xi)\overline{\mu_{0j}(M^*E^*M^{*-1}\xi)}=
        -e^{-2\pi i (r_l^E,\xi)}e^{2\pi i (r_j^E,\xi)}
            \widetilde\mu_{0l}(\xi)\overline{\mu_{0j}(\xi)}=$$
        $$=e^{2\pi i (r_j^E-r_l^E,\xi)}\w\mu_{j l}(\xi),\quad \forall E\in H,\,\, \forall  l=0,\dots,m-1,\,\, l\neq j,$$
        where $r_j^{E},$ $r_l^{E}$ defined by~(\ref{fR_j});

         $$\w\mu_{j j}(M^*E^*M^{*-1}\xi)=1-\widetilde\mu_{0j}(M^*E^*M^{*-1}\xi)\overline{\mu_{0j}(M^*E^*M^{*-1}\xi)}=
        1-\widetilde\mu_{0j}(\xi)\overline{\mu_{0j}(\xi)}=\w\mu_{j j}(\xi)\quad \forall E\in H.$$
    By Lemma~\ref{lemHSym} wavelet mask $\w m_{j}(\xi)$ is $H$-symmetric with respect
    to the point $s_{j},$ $j\in \{0,\dots,m-1\}.$ Wavelet masks $m_j,$ $j=0,\dots,m-1$
    are obviously $H$-symmetric.
    $\Diamond$

Using Lemma~\ref{lemKrSk} and Lemma~\ref{lemHSym} for appropriate symmetry groups $H$
    and matrix dilations $M$ we suggest a simple
    algorithm for the construction of mask $m_0$ that is
    $H$-symmetric with respect to
     the center $c$ and
    has arbitrary order of sum rule. The symmetry center $c\in\zd,$
    but due to the integer shifts invariance of the main properties (like orders of sum rule, linear-phase moments
    and $H$-symmetry) it is enough to consider only the case when $c=\nul.$

Let $H$ be a symmetry group on $\zd,$
    $M$ be the dilation matrix. In what follows, we assume that
    for each $E\in H$, $s_j\in D(M)$ and any trigonometric polynomial $t(\xi)$
    of semi-integer degrees associated with $\sigma= M^{-1}s_j-\lfloor M^{-1}s_j\rfloor$
    trigonometric polynomial $t(E^*\xi)$ remains trigonometric polynomial
    of semi-integer degrees associated with $\sigma= M^{-1}s_j-\lfloor M^{-1}s_j\rfloor.$ %This additional condition we will call condition $S$.
\begin{theo} %$c\in\zd$,
    Let $n\in\n$, $H$ be a symmetry group on $\zd,$
    $M$ be the dilation matrix. %The condition $S$ is valid.
    %$H$ is a symmetry group with respect to that dilation matrix $M$ and $2M^{-1}$ is an integer matrix.
    Then, there exist masks $m_0$ and $\w m_0$ that are $H$-symmetric with respect to the origin, mask $m_0$ has sum rule of order $n,$
    mask $\w m_0$ satisfies condition~(\ref{20_new}).

\label{theoHSym}
\end{theo}

{\bf Proof.} To construct such trigonometric polynomial $m_0(\xi)$ we
set
$\lambda_{\gamma}=\delta_{\gamma\nul}, \gamma\in\Delta_n.$ Therefore, by~(\ref{19}) we have the following conditions
    on the polyphase components $\mu_{0j}(\xi)$ of mask $m_0$ %$l=0\ddd m-1,$
    %$\gamma\in\zd_+$,  $[\gamma]\le n$, $\lambda_\nul=1$, such that
        \be
        D^\beta\mu_{0j}(\nul)=\frac {(2\pi i)^{[\beta]}} {\sqrt m}  (-M^{-1}s_j)^{\beta},
        \quad \forall j\in \{0,\dots,m-1\},\quad \forall \beta\in\Delta_n.
        \label{19H}
        \ee
Also, to provide  $H$-symmetry with respect to the origin
    for mask $m_0$ the polyphase components $\mu_{0j}(\xi),$ $j\in \{0,\dots,m-1\}$
    should satisfy~(\ref{fPolySymH}).

To provide the both requirements we set

        $$\mu_{0j}(\xi)=\frac 1{\sqrt m}
        G^H_{\nul}(\xi)e^{-2\pi i (M^{-1}s_j,\xi)},$$
    where  $G^H_{\nul}(\xi)\in\Theta_{\nul,n}$ are trigonometric polynomials
    of semi-integer degrees associated with $\sigma= M^{-1}s_j-\lfloor M^{-1}s_j\rfloor$
    and
    $H$-symmetric with respect to the origin. For instance, they can be easily obtained by
    $$G^H_{\nul}(\xi)=\frac 1 {\#H}\sum\limits_{E\in H}G_{\nul}(E^*\xi),$$
    where  $G_{\nul}(\xi)\in\Theta_{\nul,n} $ are trigonometric polynomials of semi-integer
    degrees associated with $\sigma= M^{-1}s_j-\lfloor M^{-1}s_j\rfloor,$
    $\#H$ is the cardinality of the set $H.$

    It is clear that  condition~(\ref{19H}) is valid.
    Symmetry conditions~(\ref{fPolySymH}) are also satisfied since
    $$\mu_{0j}(M^*E^*M^{*-1}\xi)=\frac 1{\sqrt m}G^H_{\nul}(M^*E^*M^{*-1}\xi)e^{-2\pi i (M^{-1}s_j,M^*E^*M^{*-1}\xi)}=$$
    $$\frac 1{\sqrt m}G^H_{\nul}(\xi)e^{-2\pi i (M^{-1} E s_j,\xi)}=
    \mu_{0j}(\xi)%e^{-2\pi i (M^{-1} (Es_j -s_j),\xi)}
    e^{-2\pi i (r_j^E,\xi)}, \quad \forall E\in H, \quad \forall j\in\{0,\dots,m-1\}.$$
    Thus, we construct mask $m_0$ that is $H$-symmetric with respect to the origin and has sum rule of order $n$.

It is clear that
    $$D^{\beta}m_0(0)=\delta_{\beta\nul}, \quad \forall \beta \in \Delta_n.$$
    Therefore, according to condition~(\ref{20_new}) as appropriate dual mask $\w m_0$ we simply take $\w m_0\equiv 1$ or
          $$\w m_0(\xi)=
        G^H_{\nul}(\xi),$$
    where $G^H_{\nul}(\xi)\in\Theta_{\nul,n}$
    are $H$-symmetric trigonometric polynomials
     with respect to the origin.
    Again, they can be easily obtained by
    $$G^H_{\nul}(\xi)=\frac 1 {\#H}\sum\limits_{E\in H}G_{\nul}(E^*\xi),$$
    where $G_{\nul}(\xi)\in\Theta_{\nul,n}$ are trigonometric polynomials.
    $\Diamond$

 Now, we give some examples of symmetry groups with appropriate matrix dilations. Let $H$ be the axis symmetry group on $\zd$, i.e.
        $$H=\{\texttt{diag}(\varepsilon_1,\dots,\varepsilon_d): \varepsilon_r=\pm 1, \,\forall r=1,\dots,d\}.$$
     As matrix dilation, we take $M=2P,$ where $P$ is a permutation matrix, i.e.
    square matrix that has exactly one entry $1$ in each row and each column and zeros elsewhere.
    Due to the statement~\cite[Proposition 4.1]{Han3}, $H$ is a symmetry group with respect to these dilation matrices.
    Condition~(\ref{fR_j}) is also valid. Indeed, if digits $s_j,$ $j=0,\dots,2^d$ are chosen as coordinates of the vertices of $d$-dimensional unit cube, then
    vector $Es_j-s_j$ consists of numbers $-2$ and $0$. It is clear that $M^{-1}(Es_j-s_j)\in\zd,$ $\forall E \in H,$ $j=0,\dots,2^d.$
    %Condition $S$ is also valid.

    Another example, let $H$ be the 4-fold symmetry group on $\z^2$, i.e.
        $$ H=\left\{\pm I_2, \pm \left(\begin{matrix}
            -1 & 0\cr
            0 & 1\cr
        \end{matrix}\right), \pm \left(\begin{matrix}
            0 & 1\cr
            -1 & 0\cr
        \end{matrix}\right),
        \pm \left(\begin{matrix}
            0 & 1\cr
            1 & 0\cr
        \end{matrix}\right) \right\}.$$
         As matrix dilation, we take the quincunx matrix dilation
         $$M= \left(\begin{matrix} 1&1\cr  1&-1 \cr  \end{matrix}\right) \quad \texttt{or}\quad
         M= \left(\begin{matrix} 1&-1\cr  1&1 \cr  \end{matrix}\right). $$
         Not hard to check that $H$ is a symmetry group with respect to these dilation matrices.
         The set of digits is $D(M)=\{s_0=(0,0), s_1=(1,0)\}.$
         Condition~(\ref{fR_j}) is also valid. Assumption about trigonometric polynomials of
         semi-integer degrees is satisfied for the both case.

\section{Examples}

In this section we
    shall give several examples to illustrate the main results of this paper.
    All examples are based on the construction of the
    refinable mask $m_0$, dual refinable mask $\w m_0$, wavelet masks $m_{\nu}, $ $\w m_{\nu},$ $\nu=1,\dots,m$ by Theorem~\ref{theoMaskA},
    Theorem~\ref{DualTheo} and Theorem~\ref{theoWaveA}
    for point symmetry
    and by Theorem~\ref{theoAxisM}-\ref{theoWaveAxisS} for axis symmetry.

    1. Let $c=\nul,$ ${M}= \left(\begin{matrix} 1&-2\cr
    2&-1 \cr
    \end{matrix}\right).
    $
        The set of digits is $$D(M)=\{s_0=(0,0), s_1=(0,-1), s_2=(0,1)\},$$ $m=3.$
        Let us construct mask $m_0$ that is symmetric with respect to the origin
        and satisfies condition~(\ref{19mod}) with $n=4$ and $\lambda'_{\gamma}=\delta_{\gamma\nul},$ $[\gamma]<4.$
        According to Theorem~\ref{theoMaskA},
        %polyphase components can be take as follows $$\mu_{0j}=\frac 1 {\sqrt 3}, \mu_{0i}.$$
        mask $m_0$
        can be constructed as follows

            $$m_0: \frac 1 {243} \left(\begin{matrix}
       0      &        0       &       0  &            0       &      -3  &            0       &       0\cr
       0      &        0      &        0   &          -5      &        0   &           3        &     -4\cr
      -1       &       0     &         0    &         24     &        33    &          0         &    -5\cr
       0        &      3    &         36     &        \textbf{81 }   &         36     &         3          &    0\cr
      -5         &     0   &          33      &       24   &           0      &        0           &  -1\cr
      -4          &    3  &            0       &      -5  &            0       &       0            &  0\cr
       0           &   0 &            -3        &      0 &             0        &      0             & 0  \cr
    \end{matrix}\right)$$
    with support in $[-3,3]^2\bigcap\z^2.$
        The corresponding refinable function $\phi$ is in $L_2(\r^2)$ and $\nu_2(\phi)\ge 2.3477.$
        This value and the subsequent ones are calculated by~\cite[Theorem 7.1]{Han03Vect}.
        As for the dual mask we simply take $\w m_0 \equiv 1.$
        Condition~(\ref{20_new}) is obviously valid and $\w \phi (x)$ is the Dirac delta function.
        Due to Theorem~\ref{theoWaveA} we obtain wavelet masks

            $$m_1=\sqrt 3 (e^{2\pi i (s_1,\xi)}+e^{2\pi i (s_2,\xi)}),\quad
            m_2=\sqrt 3 (e^{2\pi i (s_1,\xi)}-e^{2\pi i (s_2,\xi)});$$

            $$\w m_1: \frac 1 {162\sqrt 3} \left(\begin{matrix}
       0    &          0       &       0&              5       &       0 &             0        &      4 \cr
       0      &        0      &        0 &             0      &      -30  &            0       &       0 \cr
       0      &        0     &       -36  &           81     &         0   &          -6      &        0 \cr
       6     &         0    &          0   &         \textbf{-48}    &          0    &          0     &         6 \cr
       0      &       -6   &           0    &         81   &         -36     &         0    &          0 \cr
       0       &       0  &          -30     &         0  &            0      &        0   &           0 \cr
       4        &      0 &             0      &        5 &             0       &       0  &            0    \cr
    \end{matrix}\right),\quad $$
            $$\w m_2: \frac 1 {162\sqrt 3} \left(\begin{matrix}
       0     &         0      &        0   &           5       &       0 &             0 &             4\cr
       0      &        0     &         0    &          0      &      -36  &            0  &            0\cr
       0       &       0    &        -36     &        81     &         0   &           0   &           0\cr
       4        &      0   &           0      &       \textbf{ 0 }   &          0    &          0    &         -4\cr
       0         &     0  &           0        &    -81    &         36      &        0      &        0\cr
       0          &    0 &            36        &      0  &            0      &        0      &        0\cr
      -4           &   0 &             0         &    -5 &             0       &       0       &       0   \cr    \end{matrix}\right)\quad$$
      with supports in $[-3,3]^2\bigcap\z^2.$ Wavelet masks $\w m_1, \w m_2$ have vanishing moments of order 4.
      Thus, we are in the conditions of Theorem~\ref{t7}. The corresponding symmetric/antisymmetric almost frame-like wavelet system
      provide approximation order 4 according to~(\ref{6}).

      2. Let $c=(1,1),$ ${M}= \left(\begin{matrix} 2&0\cr  0&2 \cr  \end{matrix}\right). $
      The set of digits is $$D(M)=\{s_0=(0,0), s_1=(0,1), s_2=(1,0), s_3=(1,1)\},$$ $m=4.$
      Let us construct mask $m_0$ that is axial symmetric with respect to the center $c$
        and satisfies condition~(\ref{19mod}) with $n=2.$
        According to the Theorem~\ref{theoAxisMS}
        %polyphase components can be take as follows $$\mu_{0j}=\frac 1 {\sqrt 3}, \mu_{0i}.$$
        mask $m_0$ can be constructed as follows

            $$m_0: \left(\begin{matrix}
       0       &       1/16     &      1/16 &          0\cr
       1/16     &      1/8     &       1/8   &         1/16\cr
       1/16      &    \textbf{ 1/8}    &        1/8    &        1/16\cr
       0          &    1/16  &         1/16    &       0      \cr
      \end{matrix}\right)\quad $$
      with support
        in $[-1,2]^2\bigcap\z^2$.
      The corresponding refinable function $\phi$ is in $L_2(\r^2)$ and $\nu_2(\phi)\ge 2.$
        As for the dual mask we take
            $$\w m_0: \left(\begin{matrix}
       1/4     &      1/4    \cr
        \textbf{1/4}&1/4 \cr
      \end{matrix}\right)$$
      with support in $[0,1]^2\bigcap\z^2.$
        Condition~(\ref{20_new}) is valid. The corresponding dual
        refinable function $\w\phi$ is in $L_2(\r^2)$ and $\nu_2(\phi)\ge \frac 12.$
        Due to Theorem~\ref{theoWaveAxisS} we obtain wavelet masks

        $$m_1: \frac 12\left(\begin{matrix}
   1     &      1    \cr
        \textbf{1}& 1  \cr
      \end{matrix}\right) \quad  m_2: \frac 12\left(\begin{matrix}
        -1     &      1    \cr
        \textbf{1}& -1  \cr
      \end{matrix}\right) \quad m_3: \frac 12\left(\begin{matrix}
   -1     &    -1    \cr
       \textbf{ 1}& 1  \cr
      \end{matrix}\right) \quad m_4: \frac 12\left(\begin{matrix}
     1     &      -1    \cr
       \textbf{ 1}& -1  \cr
      \end{matrix}\right)$$
      with supports in $[0,1]^2\bigcap\z^2.$

Dual wavelet masks are

    $$\w m_1: \frac 18\left(\begin{matrix}
       0 &0& -1/8 &-1/8  &0& 0 \cr
        0 &0 &-1/8 &-1/8 &0 & 0\cr
        -1/8 &-1/8& 1/2 &1/2 &-1/8& -1/8 \cr
        -1/8& -1/8 &\textbf{1/2} &1/2& -1/8 &-1/8 \cr
            0 &0 &-1/8 &-1/8& 0 &0 \cr
        0& 0& -1/8 &-1/8 & 0 &0\cr
       \end{matrix}\right) \quad
        \w m_2: \frac 18\left(\begin{matrix}
        -1     &      1    \cr
        \textbf{1}& -1  \cr
      \end{matrix}\right)  \quad$$
with supports in $[-2, 3]^2\bigcap\z^2$ and $[0,1]^2\bigcap\z^2$ accordingly and
    $$
    \w m_3: \frac 18\left(\begin{matrix}
    1/8 &1/8\cr
    1/8 &1/8 \cr
    -1 &-1  \cr
    \textbf{1} &1\cr
    -1/8 &-1/8 \cr
     -1/8 &-1/8\cr
      \end{matrix}\right) \quad
    \w m_4: \frac 18\left(\begin{matrix}
    -1/8 &-1/8& 1 &-1 &1/8& 1/8 \cr
    -1/8& -1/8 &\textbf{1} &-1& 1/8 &1/8 \cr
      \end{matrix}\right).
        $$
    with supports in $[0, 1]\times[-2,3]\bigcap\z^2 $ and $[-2, 3]\times[0,1] \bigcap\z^2$ accordingly.

The corresponding axial symmetric/antisymmetric frame-like wavelet system provide
    approximation order 2 according to~\cite[Theorem 18]{KrSk}.% for any $f\in S$ in $S'.$

    3. Let $c=\nul,$ ${M}= \left(\begin{matrix} 1&1\cr  1&-1 \cr  \end{matrix}\right).$
        The set of digits is $D(M)=\{s_0=(0,0), s_1=(1,0)\},$ $m=2.$
        Let us construct mask $m_0$ that is 4-fold symmetric with respect to the origin
        and satisfies condition~(\ref{19mod}) with $n=2$.
        According to Theorem~\ref{theoHSym}
        %polyphase components can be take as follows $$\mu_{0j}=\frac 1 {\sqrt 3}, \mu_{0i}.$$
        mask $m_0$ can be constructed as follows

        $$m_0: \left(\begin{matrix}
             0   &          -1/16 &          0 &            -1/16 &          0\cr
      -1/16       &    0      &        1/4      &      0      &       -1/16\cr
       0           &   1/4   &         \textbf{1/2   }    &     1/4   &         0\cr
      -1/16         &  0    &          1/4        &    0    &         -1/16\cr
       0             &-1/16&           0           &  -1/16&           0       \cr
      \end{matrix}\right)$$
      with support
        in $[-2,2]^2\bigcap\z^2.$
      The corresponding refinable function $\phi$ is in $L_2(\r^2)$ and $\nu_2(\phi)\ge 0.783.$
     As for the dual mask we take  $\w m_0=1.$
      Condition~(\ref{20_new}) is obviously valid and $\w \phi (x)= \delta (x).$

     Due to Theorem~\ref{theoWaveHSym} we obtain wavelet masks

     $$m_1 = e^{2\pi i \xi_1}, \quad \w m_1:  \frac 1 {\sqrt 2}\left(\begin{matrix}
     0          &    1/8      &      0    &          1/8       &     0\cr
       1/8       &     0     &        -1/2 &           0      &        1/8\cr
       0          &   \textbf{-1/2 } &          1    &         -1/2   &         0\cr
       1/8         &   0   &          -1/2   &         0    &          1/8\cr
       0            &  1/8&            0      &        1/8 &           0     \cr
      \end{matrix}\right),$$
with support in $[-1,3]\times[-2,2]\bigcap\z^2.$ Wavelet mask $\w m_1$ has vanishing moments of order 2.

% for any $f\in S$ in $S'.$
    Thus, we are in the conditions of Theorem~\ref{t7}. The corresponding 4-fold symmetric frame-like wavelet system provide
    approximation order 2 according to~(\ref{6}).

    \end{document}